\documentclass{article}
\pdfoutput=1

\usepackage{arxiv}

\usepackage[utf8]{inputenc} 
\usepackage[T1]{fontenc}    
\usepackage{hyperref}       
\usepackage{url}            
\usepackage{booktabs}       
\usepackage{amsfonts}       
\usepackage{nicefrac}       
\usepackage{microtype}      
\usepackage{lipsum}
\usepackage[english]{babel}
\usepackage[final]{graphicx}

\title{ANALYSIS OF THE RIEMANN ZETA FUNCTION}

\author{
 Kirill V.~Kapitonets \\
 BAUMAN MSTU, GRADUATE 1990\\
  MCC EuroChem\\
  Moscow\\
  Russian Federation\\
  \texttt{kkapitonets@live.com} \\
}

\begin{document}
\maketitle
\begin{abstract}
The paper uses a feature of calculating the Riemann Zeta function in the critical strip, where its approximate value is determined by partial sums of the Dirichlet series, which it is given.
\par
These expressions are called the first and second approximate equation of the Riemann Zeta function.
\par
The representation of the terms of the Dirichlet series by vectors allows us when analyzing the polyline formed by these vectors:
\par
1) explain the geometric meaning of the generalized summation of the Dirichlet series in the critical strip;
\par
2) obtain formula for calculating the Riemann Zeta function;
\par
3) obtain the functional equation of the Riemann Zeta function based on the geometric properties of vectors forming the polyline;
\par
4) explain the geometric meaning of the second approximate equation of the Riemann Zeta function;
\par
5) obtain the vector equation of non-trivial zeros of the Riemann Zeta function;
\par
6) determine why the Riemann Zeta function has non-trivial zeros on the critical line;
\par
7) understand why the Riemann Zeta function cannot have non-trivial zeros in the critical strip other than the critical line.
\par
The main result of the paper is a definition of possible ways of the confirmation of the Riemann hypothesis based on the properties of the vector system of the second approximate equation of the Riemann Zeta function.

\end{abstract}

\section{Introduction}
A more detailed title of the paper can be formulated as follows:
\par
\textit{Geometric analysis of expressions that determine a value of the Riemann zeta function in the critical strip.}
\par
However, this name will not fully reflect the content of the paper.
It is necessary to refer to the sequence of problems.
The starting point of the analysis of the Riemann zeta function was the question:
\par
\textit{Why the Riemann zeta function has non-trivial zeros.}
\par
A superficial study of the question led to the idea of geometric analysis of the Dirichlet series.
\par
Preliminary geometric analysis allowed to obtain the following results:
\par
1) definition of an expression for calculating a value of the Riemann zeta function in the critical strip;
\par
2) derivation of the functional equation of the Riemann zeta function based on the geometric properties of vectors;
\par
3) explanation of the geometric meaning of the second approximate equation of the Riemann zeta function.
\par
The most important result of the preliminary geometric analysis was the question why the Riemann zeta function cannot have non-trivial zeros in the critical strip except the critical line.
\par
Therefore, the extended title of the paper is follows:
\par
\textit{Geometric analysis of the expressions defining a value of the Riemann zeta function in the critical strip with the aim of answering the question of why the Riemann zeta function cannot have non-trivial zeros in the critical strip, except for the critical line.}
\par
At once, it is necessary to define the upper limit of the paper.
Although that the paper explores the main question of the Riemann hypothesis, we are not talking about its proof.
\par
The results of the geometric analysis of the expressions defining a value of the Riemann zeta function in the critical strip is the definition of the possible ways of the confirmation of the Riemann hypothesis based on the results of this analysis.
\par
In addition, if we are talking about the boundaries, then we will mark the lower limit of the paper.
We will not use the latest advances in the theory of the Riemann zeta function and improve anyone's result (we will state our position on this later).
\par
We turn to the origins of the theory of the Riemann zeta function and base on the very first results obtained by Riemann, Hardy, Littlewood and Titchmarsh, as well as by Adamar and Valle Poussin.
\par
In addition, of course, we cannot do without the official definition of the Riemann hypothesis:
\par
\textit{The Riemann zeta function is the function of the complex variable s, defined in the half-plane Re(s) > 1 by the absolutely convergent series:}
\begin{equation}\label{zeta_dirichlet}\zeta(s) = \sum\limits_{n=1}^{\infty}\frac{1}{n^s};\end{equation}
\textit{and in the whole complex plane C by analytic continuation. As shown by Riemann, $\zeta(s)$ extends to C as a meromorphic function with only a simple pole at s =1, with residue 1, and satisfies the functional equation:}
\begin{equation}\label{zeta_func_eq}\pi^{-s/2}\Gamma(\frac{s}{2})\zeta(s)=\pi^{-(1-s)/2}\Gamma(\frac{1-s}{2})\zeta(1-s);\end{equation}
\textit{Thus, in terms of the function $\zeta(s)$, we can state Riemann hypothesis: The non-trivial zeros of $\zeta(s)$ have real part equal to 1/2.}
\par
We define two important concepts of the Riemann zeta function theory: the critical strip and the critical line.
\par
Both concepts relate to the region where non-trivial zeros of the Riemann zeta function are located.
\par
Adamar and Valle-Poussin independently in the course of the proof of the theorem about the distribution of prime numbers, which is based on the behavior of non-trivial zeros of the Riemann zeta function , showed that all non-trivial zeros Riemann zeta function is located in a narrow strip $0<Re(s)<1$.
In the theory of the Riemann zeta function, this strip is called \textit {critical strip}.
\par
The line $Re(s)=1/2$ referred to in the Riemann hypothesis is called \textit {critical line}.
\par
It is known that the Dirichlet series (\ref{zeta_dirichlet}), which defines the Riemann zeta function, diverges in the critical strip, therefore, to perform an analytical continuation of the Riemann zeta function, \textit{the method of generalized summation} of divergent series, in particular the generalized Euler-Maclaurin summation formula is used.
\par
We should note that the generalized Euler-Maclaurin summation formula is used in the case where the partial sums of a divergent series are suitable for computing the generalized sum of this divergent series.
\par
From the theory of generalized summation of divergent series \cite{HA3}, it is known that any method of generalized summation, if it yields any result, yields the same result with other methods of generalized summation.
\par
On the one hand, this fact confirms the uniqueness of the analytical continuation of the function of a complex variable, and on the other hand, apparently, the principle of analytical continuation extends in the form of a generalized summation to the functions of a real variable, which are defined in divergent series. 
\par
Now we can formulate the principles of geometric analysis of expressions that define a value of the Riemann zeta function in the critical strip:
\par
1) a complex number can be represented geometrically using points or vectors on the plane;
\par
2) the terms of the Dirichlet series, which defines the Riemann zeta function, are complex numbers;
\par
3) a value of the Riemann zeta function in the region of analytic continuation can be represented by any method of generalized summation of divergent series.
\par
What exactly is the geometric analysis of values of the Riemann zeta function, we will see when we construct a polyline, which form a vector corresponding to the terms of the Dirichlet series, which defines the Riemann zeta function.
\par
All results in the paper are obtained empirically by calculations with a given precision (15 significant digits).
\section{Representation of the Riemann zeta function value by a vector system}
This section presents detailed results of geometric analysis of expressions that determine a value of the Riemann zeta function in the critical strip.
\par
The possibility of such an analysis is due to two important facts:
\par
1) representation of complex numbers by vectors;
\par
2) using of methods of generalized summation of divergent series.
\subsection{Representation of the partial sum of the Dirichlet series by a vector system}
Each term of the Dirichlet series (\ref{zeta_dirichlet}), which defines the Riemann zeta function, is a complex number $x_n+iy_n$, hence it can be represented by the vector $X_n=(x_n, y_n)$.
\par
To obtain the coordinates of the vector $X_n(s), s=\sigma+it$, we first present the expression $n^{-s}$ in the exponential form of a complex number, and then, using the Euler formula, in the trigonometric form of a complex number:
\begin{equation}\label{x_vect_direchlet}X_n(s)=\frac{1}{n^{s}}=\frac{1}{n^{\sigma}}e^{-it\log(n)}=\frac{1}{n^{\sigma}}(\cos(t\log(n))-i\sin(t\log(n)));\end{equation}
\par
Then the coordinates of the vector $X_n(s)$ can be calculated by formulas:
\begin{equation}\label{x_n_y_n}x_n(s)=\frac{1}{n^{\sigma}}\cos(t\log(n)); y_n(s)=-\frac{1}{n^{\sigma}}\sin(t\log(n));\end{equation}
\par
\textit{Using the rules of analytical geometry, we can obtain the coordinates corresponding to the partial sum $s_m(s)$ of the Dirichlet series (\ref{zeta_dirichlet}):}
\begin{equation}\label{s_m_x_s_m_y}s_m(s)_x=\sum_{n=1}^{m}{x_n(s)}; s_m(s)_y=\sum_{n=1}^{m}{y_n(s)};\end{equation}
\begin{figure}[ht]
\centering
\includegraphics[scale=0.6]{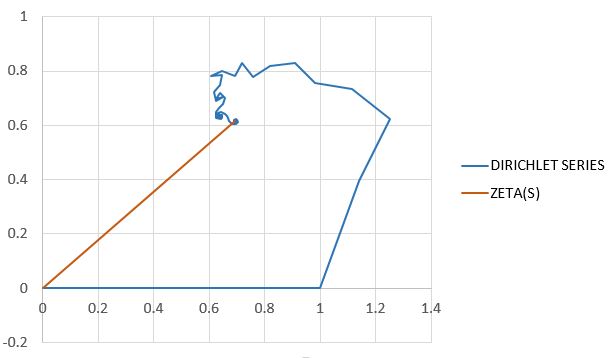}
\caption{Polyline, $s=1.25+279.229250928i$}
\label{fig:s1_1}
\end{figure}
\par
We construct a polyline corresponding to partial sums of $s_m(s)$.
Select a value \textit{in the convergence region} of the Dirichlet series (\ref{zeta_dirichlet}), for example, $s=1.25+279.229250928 i$.
\par
We also display the vector $(0.69444570272324, 0.61658346971775)$, which corresponds to a value of the Riemann zeta function at $s=1.25+279.229250928 i$.
\par
We will display the first $m=90$ vectors (later we will explain why we chose such a number) so that the vectors follow in ascending order of their numbers (fig. \ref{fig:s1_1}).
\par
Now we change a value of the real part $s=0.75+279.229250928 i$ and move to the region where the Dirichlet series (\ref{zeta_dirichlet}) \textit{diverges.}
\begin{figure}[ht]
\centering
\includegraphics[scale=0.6]{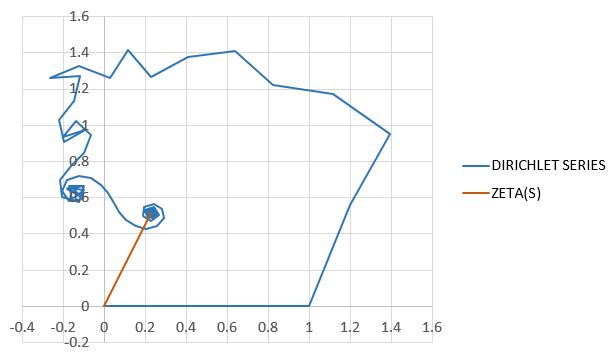}
\caption{Polyline, $s=0.75+279.229250928i$}
\label{fig:s3_1}
\end{figure}
\par
We also display the vector $(0.22903651233853, 0.51572970834588)$, which corresponds to a value of the Riemann zeta function at $s=0.75+279.229250928 i$.
\par
We observe (fig. \ref{fig:s3_1}) an increase of size of the polyline, but the qualitative behavior of the graph does not change.
\textit{The polyline twists around a point corresponding to a value of the Riemann zeta function.}
\par
To see what exactly is the difference, we need to consider the behavior of vectors with smaller numbers, for example, in the range $m=(300, 310)$ with relation to vectors with large numbers, for example, in the range $m=(500, 510)$.
\begin{figure}[ht]
\centering
\includegraphics[scale=0.6]{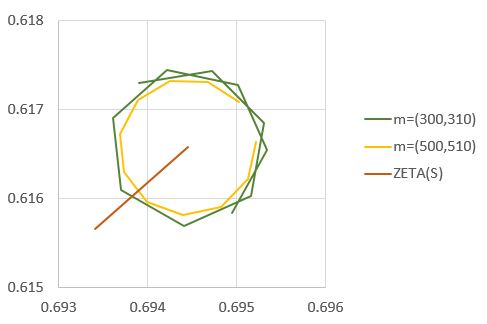}
\caption{Part of polyline, $s=1.25+279.229250928i$}
\label{fig:s1_2_1}
\end{figure}
\par
We see (fig. \ref{fig:s1_2_1}) that when $\sigma=1.25$ the polyline is \textit{a converging} spiral, as the radius of the spiral in the range of $m=(300, 310)$ is greater than the radius of the spiral in the range of $m=(500, 510)$.
\begin{figure}[ht]
\centering
\includegraphics[scale=0.6]{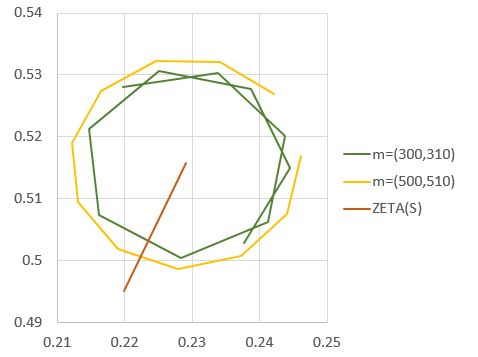}
\caption{Part of polyline, $s=0.75+279.229250928i$}
\label{fig:s3_2_1}
\end{figure}
\par
Unlike the first case, when $\sigma=0.75$ we see (fig. \ref{fig:s3_2_1}) that a polyline is \textit{a divergent} spiral, since the radius of the spiral in the range $m=(300, 310)$ is less than the radius of the spiral in the range $m=(500, 510)$.
\par
But regardless of whether the spiral converges or diverges, in both cases we see that the center of the spiral corresponds to the point corresponding to a value of the Riemann zeta function (we will later show that it is indeed).
\par
This fact can be considered a geometric explanation of the method of generalized summation:
\par
\textit{The coordinates of the partial sums $s_m(s)_x$ and $s_m(s)_y$ of the Dirichlet series (\ref{zeta_dirichlet}) vary with relation to some middle values of $x$ and $y$, which we take as a value of the Riemann zeta function, only in one case the coordinates of the partial sums converge infinitely to these values, and in the other case they diverge infinitely with relation to these values.}
\subsection{Properties of vector system of the partial sums of the Dirichlet series – the Riemann spiral}
In order not to use in the further analysis of the long definitions of the vector system, which we are going to study in detail, we introduce the following definition:
\par
\textit{The Riemann spiral is a polyline formed by vectors corresponding to the terms of the Dirichlet series (\ref{zeta_dirichlet}) defining the Riemann zeta function, in ascending order of their numbers.}
\par
A preliminary analysis of \textit{Riemann spiral} showed that vectors arranged in ascending order of their numbers form a converging or divergent spiral with a center at a point corresponding to a value of the Riemann zeta function.
\par
We also see (fig. \ref{fig:s3_1}) that the Riemann spiral has several more centers where the vectors first form a converging spiral and then a diverging spiral.
\par
Comparing (fig. \ref{fig:s1_1}) and (fig. \ref{fig:s3_1}) we see that the number of such centers does not depend on the real part of a complex number.
\par
\begin{figure}[ht]
\centering
\includegraphics[scale=0.6]{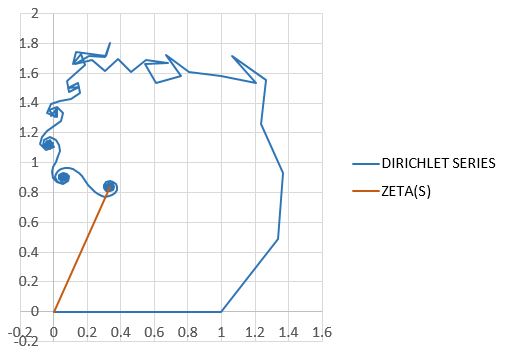}
\caption{The Riemann spiral, $s=0.75+959.459168807i$}
\label{fig:s7_1}
\end{figure}

While comparing (fig. \ref{fig:s3_1}) and (fig. \ref{fig:s7_1}), we see that the number of such centers increases when the imaginary part of a complex number increases.
\par
We can explain this behavior of Riemann spiral vectors if we return to the exponential form (\ref{x_vect_direchlet})of the record of terms of the Dirichlet series (\ref{zeta_dirichlet}).
\par
\begin{figure}[ht]
\centering
\includegraphics[scale=0.6]{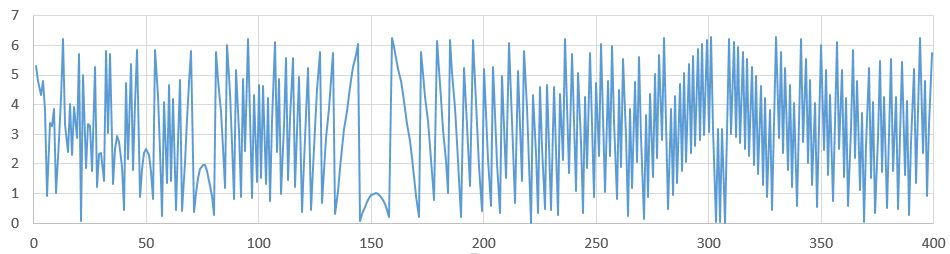}
\caption{Absolute angles $\varphi_n(t)$ of vectors of the Riemann spiral, rad, $s=0.75+959.459168807i$}
\label{fig:absolute_angles}
\end{figure}

The paradox of the Riemann spiral is that with unlimited growth of the function $t\log(n)$, the absolute angles $\varphi_n(t)$ of its vectors behave in a pseudo-random way (fig. \ref{fig:absolute_angles}), because we can recognize angles only in the range $[0, 2\pi]$.
\begin{equation}\label{phi_n}\varphi_n(t)=t\log(n)mod\ 2\pi;\end{equation}
\par
\begin{figure}[ht]
\centering
\includegraphics[scale=0.6]{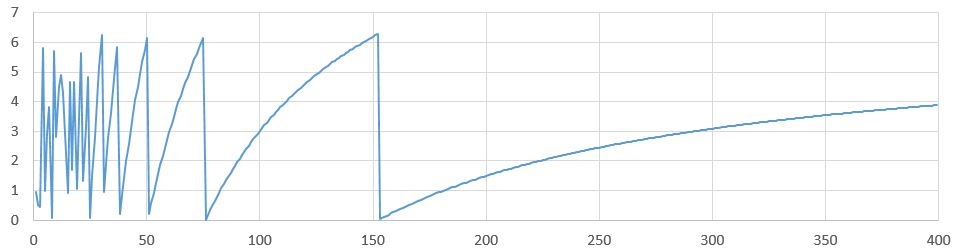}
\caption{Relative angles $\Delta\varphi_n(t)$ of vectors of the Riemann spiral, rad, $s=0.75+959.459168807i$}
\label{fig:relative_angles}
\end{figure}

The angles between the vectors $\Delta\varphi_n(t)$, if they are measured not as visible angles between segments, but as angles between directions, permanently grow (fig. \ref{fig:relative_angles}) to a value $2\pi$, then they sharply decrease to a value $0$ and again grow to a value $2\pi$.
\begin{equation}\label{delta_phi_n}\Delta\varphi_n(t)=\varphi_n(t)-\varphi_{n-1}(t);\end{equation}
\par
On the last part this growth is \textit{asymptotic,} i.e. no matter how large the vector number $n$ is, a value $\log(n)$ will never be equal to a value $\log(n+1)$.
\par
In consequence of the revealed properties of the Riemann spiral vectors, we observe two types of singular points:
\par
1) \textit{reverse points} where the visible twisting of the vectors is replaced by untwisting, these points are multiples of $(2k-1)\pi$;
\par
2) \textit{inflection points}, in which the visible untwisting of the vectors is replaced by twisting, these points are multiples of $2k\pi$;
\par
Now we can explain why for a value of $s=1.25+279.229250928 i$ we took 90 vectors to construct the Riemann spiral:
\par
\begin{equation}\label{m1}279.229250928/\pi=88.881431082;\end{equation}
\par
Therefore, the first reverse point is between the 88th and 89th vectors, we just rounded the number of vectors to a multiple.
\par
\textit{This is the number of vectors we must use to build the Riemann spiral in order to vectors fully twisted around the point corresponding to a value of the Riemann zeta function.}
\par
In addition, we can determine the number of reverse points $m$, as we remember, this number determines the range in which (fig. \ref{fig:relative_angles}) the periodic monotonous increase of angles between the Riemann spiral vectors is observed, moreover this number plays an important role (as we will see later) in the representation of the Riemann zeta function values by the vector system.
\par
We can determine the number of reverse points $m$ from the condition that between two reverse points there is at least one vector:
\par
\begin{equation}\label{m2_eq}\frac{t}{(2m-1)\pi}-\frac{t}{(2m+1)\pi}=1;\end{equation}
\par
From this equation we find:
\par
\begin{equation}\label{m2}m=\sqrt{\frac{t}{2\pi}+\frac{1}{4}};\end{equation}
\par
At the end of the consideration of the static parameters of the Riemann spiral, we perform an analysis of its radius of curvature:
\par
\begin{equation}\label{curvature_radius_eq}r_n=\frac{|X_n|\cos(\Delta\varphi_n)}{\sqrt{1-\cos(\Delta\varphi_n)^2}};\end{equation}
\par
where
\begin{equation}\label{cos_delta_varphi_n}\cos(\Delta\varphi_n)=\frac{(X_n,X_{n-1})}{|X_n||X_{n-1}|};\end{equation}
\begin{equation}\label{module}|X_n|=\sqrt{x_n^2+y_n^2};\end{equation}
\begin{equation}\label{scalar_product}(X_n,X_{n-1})=x_nx_{n-1}+y_ny_{n-1};\end{equation}
\begin{figure}[ht]
\centering
\includegraphics[scale=0.6]{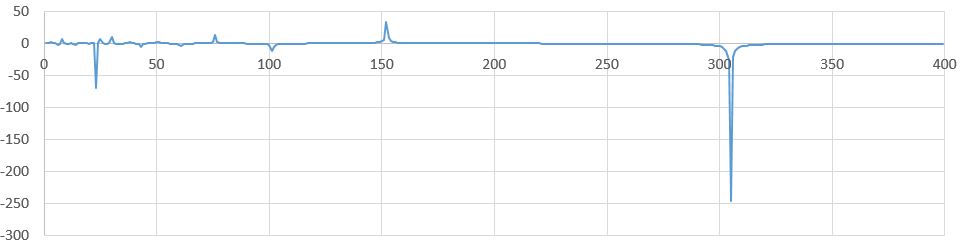}
\caption{Curvature radius of the Riemann spiral, $s=0.75+959.459168807i$}
\label{fig:curvature_radius}
\end{figure}
\par
We see (fig. \ref{fig:curvature_radius}) that the Riemann spiral has \textit{an alternating sign} of radius of curvature.
This is the only spiral that has this property.
\par
The maximum negative value of the radius of curvature takes at the reverse point, and the maximum positive is at the inflection point. 
\subsection{Derivation of an empirical expression for the Riemann zeta function}
We study in detail the behavior of Riemann spiral vectors after the first reverse point: $$m=\frac{t}{\pi};$$
\par
As we know, the Riemann spiral vectors in the critical strip form a divergent spiral (fig. \ref{fig:s3_2_1}).
\par
\begin{figure}[ht]
\centering
\includegraphics[scale=0.6]{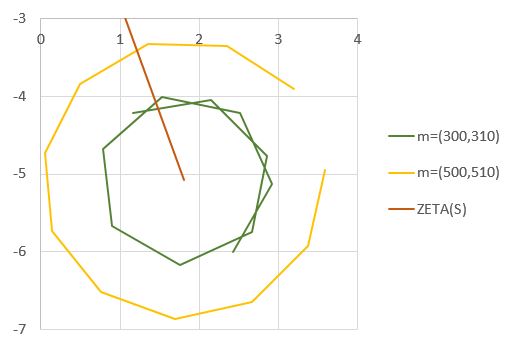}
\caption{Part of polyline, $s=0+279.229250928i$}
\label{fig:s8_2_1}
\end{figure}

Consider the Riemann spiral after the first reverse point on the left boundary of the critical strip (Fig. \ref{fig:s8_2_1}) i.e. when $\sigma=0$.
\par
We see that the size of the polyline has increased in comparison with (fig. \ref{fig:s3_2_1}), but the behavior of the vectors has not changed, as their numbers increase, the vectors tend to form a regular polygon and the quantity of edges grows with the growth of vector numbers.
\par
And as we understand, as a result of the vectors tend to form a circle of the greater radius, than the greater the number of vectors, and the center of the circle tends to the center of the spiral, consequently, to the point, which corresponds to a value of the Riemann zeta function.
\par
In analytic number theory, this fact\footnote{The first is the correspondence to the graph of partial sums of the Dirichlet series considered by Erickson in his paper \cite{ER}} is known as the first approximate equation of the Riemann zeta function:
\begin{equation}\label{zeta_eq_1}\zeta(s)=\sum_{n\le x}{\frac{1}{n^s}}-\frac{x^{1-s}}{1-s}+\mathcal{O}(x^{-\sigma}); 
\sigma>0; |t|< 2\pi x/C; C > 1\end{equation}
\par
As we know, Hardy and Littlewood obtained this approximate equation based on the generalized Euler-Maclaurin summation method.
\par
It is known from the theory of generalized summation of divergent series \cite{HA3} that we can apply \textit{any other method} of generalized summation (if it gives any value) and get the same result.
\par
We use the Riemann spiral to determine another method of generalized summation.
\par
\begin{figure}[ht]
\centering
\includegraphics[scale=0.6]{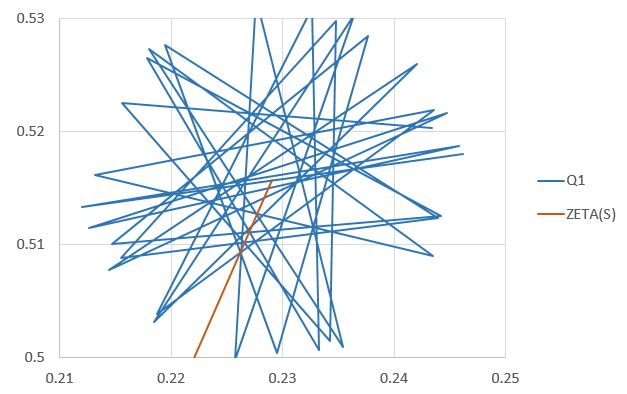}
\caption{First 30 vectors after the first reverse point, $s=0.75+279.229250928i$}
\label{fig:q_1}
\end{figure}

Consider the first 30 vectors of the Riemann spiral after the first reverse point (fig. \ref{fig:q_1}), here the vectors form a star-shaped polygon, with the point that corresponds to a value of the Riemann zeta function at the center of this polygon, as in the case of a divergent spiral (fig. \ref{fig:s3_2_1}).
\par
\begin{figure}[ht]
\centering
\includegraphics[scale=0.6]{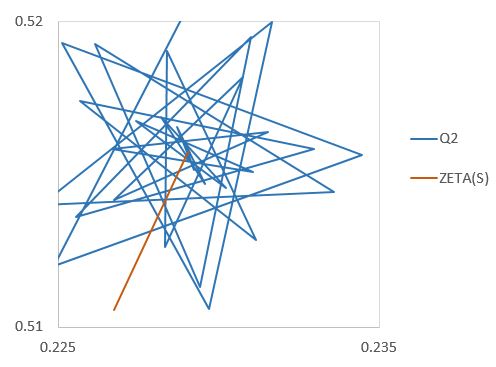}
\caption{The first step of generalized summation, $s=0.75+279.229250928i$}
\label{fig:q_2}
\end{figure}

Connect the middle of the segments formed by vectors and get 29 segments (fig. \ref{fig:q_2}).
\par
We see that the star-shaped polygon has decreased in size, and the point that corresponds to a value of the Riemann zeta function is again at the center of this polygon.
\begin{figure}[ht]
\centering
\includegraphics[scale=0.6]{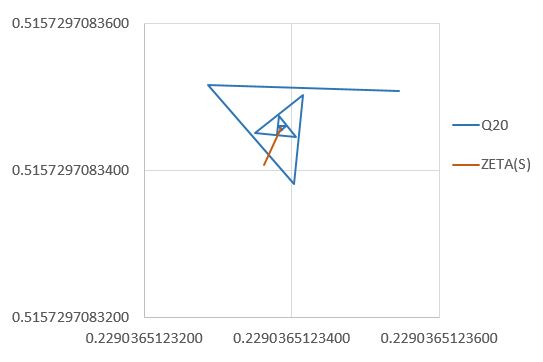}
\caption{The twentieth step of generalized summation, $s=0.75+279.229250928i$}
\label{fig:q_20}
\end{figure}
\par
We will repeat the operation of reducing the polygon (fig. \ref{fig:q_20}) until there is one segment left.
\par
 \textit{Calculations show that when calculating the coordinates of vectors with an accuracy of 15 characters, the coordinates of the middle of this segment with an accuracy of not less than 13 digits match a value of the zeta function of Riemann, for example, at point $0.75+279.229250928i$ the exact value of \cite{ZF} equals $0.22903651233853+0.51572970834588i$, and in the calculation of midpoints of the segments we get a value $0.22903651233856+0.51572970834589i$.}
\par
To calculate values of the Riemann zeta function, a reduced formula can be obtained by the described method.
Write sequentially the expressions to calculate the midpoints of the segments and substitute successively the obtained formulas one another:
\begin{equation}\label{a}a_i = \frac{x_i+x_{i+1}}{2};\end{equation}
\begin{equation}\label{b}b_i =\frac{a_i+a_{i+1}}{2}=\frac{x_i+2x_{i+1}+x_{i+2}}{4};\end{equation}
\begin{equation}\label{c}c_i =\frac{b_i+b_{i+1}}{2}= \frac{x_i+3x_{i+1}+3x_{i+2}+x_{i+3}}{8};\end{equation}
\begin{equation}\label{d}d_i =\frac{c_i+c_{i+1}}{2}=\frac{x_i+4x_{i+1}+6x_{i+2}+4x_{i+3}+x_{i+4}}{16};\end{equation}
\par
We write the same formulas for the coordinates $y_k$.
\par
In the numerator of each formula (\ref{a}-\ref{d}) we see the sum of vectors multiplied by binomial coefficients, which correspond to the degree of Newton's binomial, equal to the number of vectors minus one, and in the denominator the degree of two, equal to the number of vectors minus one.
\par
Now we can write down the abbreviated formula:
\par
\begin{equation}\label{s_x_s_y}s_x=\frac{1}{2^m}\sum_{k=0}^{m}\Big(_m^k\Big)x_k; s_y=\frac{1}{2^m}\sum_{k=0}^{m}\Big(_m^k\Big)y_k;\end{equation}
\par
where $x_k$ and $y_k$ are coordinates of partial sums (\ref{s_m_x_s_m_y}) of Dirichlet series.
\par
\textit{We obtained the formula of the generalized Cesaro summation method $(C, k)$ \cite{HA3}.}
\par
It should be noted that to calculate the coordinates of the center of the star-shaped polygon when a value of the imaginary part of the complex number is equal to 279.229250928, we use 30 vectors after the first reverse point, while starting with a value of the imaginary part of a complex number equal to 1000, it is enough to use only 10 vectors after the first reverse point.
\par
We obtained a result that applies not only to the Riemann zeta function, but also to all functions of a complex variable that have an analytic continuation.
\par
A result that relates not only to the analytic continuation of the Riemann zeta function, but to the analytic continuation as the essence of any function of a complex variable.
\par
A result that refers to any function and any physical process where such functions are applied, and hence a result that defines the essence of that physical process.
\par
Euler's intuitive belief\footnote{These inconveniences and apparent contradictions can be avoided if we give the word \glqq sum\grqq \ a meaning different from the usual. Let us say that the sum of any infinite series is a finite expression from which the series can be derived.} that a divergent series can be matched with a certain value \cite{EL}, has evolved into the fundamental theory\footnote{It is impossible to state Euler's principle accurately without clear ideas about functions of a complex variable and analytic continuation.} of generalized summation of divergent series, which Hardy systematically laid out in his book \cite{HA3}.
\par
The essence of our conclusions is as follows:
\par
1) When it comes to analytical continuation of some function of a complex variable, it automatically means that there are at least two regions of definition of this function:
\par
a) the region where the series by which this function is defined converges;
\par
b) the region where the same series diverges.
\par
\textit{Actually, the question of analytical continuation arises because there is an region where the series that defines the function of a complex variable diverges.}
\par
Thus, the analytical continuation of any function of a complex variable is inseparably linked to the fact that there is a divergent series.
\par
We can say that this is the very essence of the analytical continuation.
\par
2) Analytical continuation is possible only if there is some method of generalized summation that will give a result, i.e. some value other than infinity.
\par
This value will be a value of the function, in the region where the series by which this function is defined diverges.
\par
3) The most important thing in this question is that the series by which the function is defined must behave asymptotically in the region where this series infinitely converges and in the region where it infinitely diverges.
\par
And then we come to an important point:
\par
\textit{The function of a complex variable, if it has an analytical continuation (and therefore is given by a series that converges in one region and diverges in another, i.e. has no limit of partial sums of this series) is determined by the asymptotic law of behavior of the series by which this function is given and it is no matter whether this series converges or diverges, a value of this function will be the asymptotic value with relation to which this series converges or diverges.}
\par
Analytical continuation is possible only if the series with which the function is given has asymptotic behavior, i.e. its values oscillate with relation to the asymptote, which is a value of the function.
\par
In the case of a function of a complex variable, there are two such asymptotes, and in the case of a function of a real variable, such an asymptote is one.
\par
If the Riemann zeta function is given asymptotic values, then it also has \textit{asymptotic value of zero,} hence it may seem that the Riemann hypothesis cannot be proved.
\par
As we will show later, a value of the Riemann zeta function can be given by a finite vector system, the sum of which gives \textit{the exact value of zero} if these vectors form a polygon.
\par
\textit{In this regard, we can conclude that the Riemann hypothesis can be confirmed only if such a finite vector system exists and only using the properties of the vectors of this system.}
\par
One can disagree with the conclusion that it is possible to confirm the Riemann hypothesis using a finite vector system, but we will go this way, because the chosen method of geometric analysis of the Riemann zeta function allows us to penetrate into the essence of the phenomenon.
\subsection{Derivation of an empirical expression for the functional equation of the Riemann zeta function}
We already know one dynamic property of the Riemann spiral, which is that the number of reverse points increases with the growth of the imaginary part of a complex number.
\par
We will study \textit{the reverse points} and \textit{inflection points,} which, as we will see later, are not just special points of the Riemann spiral, they are its essential points that define the essence of the Riemann spiral and the Riemann zeta function.
\par
We will use our empirical formula for calculating a value of the Riemann zeta function, which corresponds to the Cesaro generalized summation method.
\par
As we saw (fig. \ref{fig:s7_1}), the Riemann spiral vectors at any reverse point where the Riemann spiral has a negative radius of curvature, up to the reverse point form a converging spiral, and after the reverse point form a divergent spiral.
\par
This fact allows us to apply the formula for calculating a value of the Riemann zeta function at any reverse point, because this formula allows us to calculate the coordinates of any center of the divergent spiral, which form the Riemann spiral vectors.
\par
For the convenience of further presentation, we introduce the following definition:
\par
\begin{figure}[ht]
\centering
\includegraphics[scale=0.6]{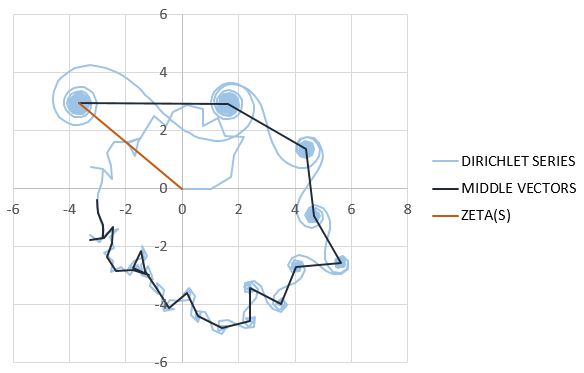}
\caption{The middle vectors of the Riemann spiral , $s=0.25+5002.981i$}
\label{fig:average_vectors}
\end{figure}

\textit{Middle vector of the Riemann spiral is a directed segment connecting two adjacent centers of the Riemann spiral, drawn in the direction from the center with a smaller number to the center with a larger number (fig. \ref{fig:average_vectors}) when numbering from the center corresponding to value of the Riemann zeta function.}
\par
To identify modulus of middle vectors of the Riemann spiral first using the formula (\ref{s_x_s_y}) we define the coordinates of the centers of the Riemann spiral.
\par
To calculate the coordinates of the first center of the Riemann spiral (a value of the Riemann zeta function) we will use 30 vectors, to calculate the coordinates of the second center - 20 vectors.
\par
To increase the accuracy, we can choose a different number of vectors for each center of the Riemann spiral, but since we choose a sufficiently large value of the imaginary part of a complex number and a sufficiently small number of vectors, it will be enough to use 5 vectors to calculate the coordinates starting from the third center of the Riemann spiral.
\par
Then, using the formula to determine the distance between two points, we find the modulus of the middle vectors:
\begin{equation}\label{dist}|Y_n|=\sqrt{(x_{n+1}-x_n)^2+(y_{n+1}-y_n)^2};\end{equation}
\par
\begin{figure}[ht]
\centering
\includegraphics[scale=0.6]{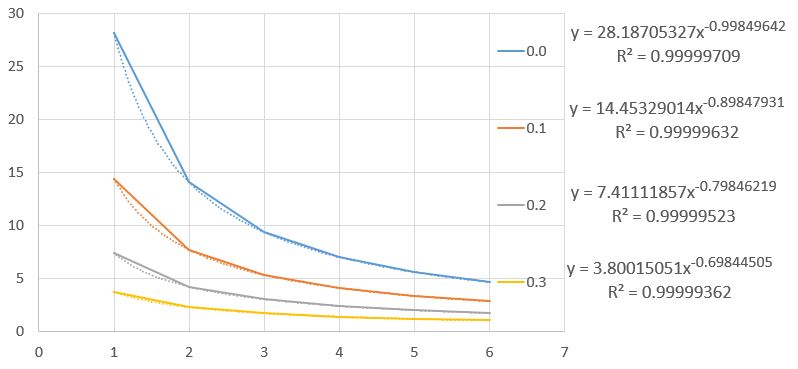}
\caption{The dependence of the modulus of the middle vector from the sequence number, $Re(s)=(0.0, 0.1, 0.2, 0.3)$, $Im(s)=5000$}
\label{fig:v1_5000_6}
\end{figure}

We calculate the modulus of the first six middle vectors of the Riemann spiral for values of the real part in the range from 0 to 1 in increments of 0.1 and a fixed value of the imaginary part of a complex number - 5000 and then we approximate the dependence of the modulus of the middle vector of the Riemann spiral on its sequence number (fig. \ref{fig:v1_5000_6}).
\par
The best method of approximation of dependence of the modulus of the middle vector from the sequence number (fig. \ref{fig:v1_5000_6}, where $x=n$) for different values of the real part of a complex number (and the fixed imaginary part of a complex number) is a power function
\begin{equation}\label{depence_1}|Y_n|=An^B;\end{equation}
\par
\begin{figure}[ht]
\centering
\includegraphics[scale=0.6]{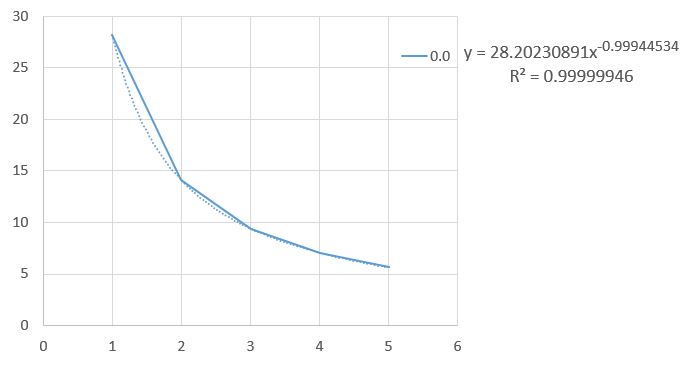}
\caption{The dependence of the modulus of the middle vector for the five vectors , $s=0+5000i$}
\label{fig:v0_5000_5}
\end{figure}
\begin{figure}[ht]
\centering
\includegraphics[scale=0.6]{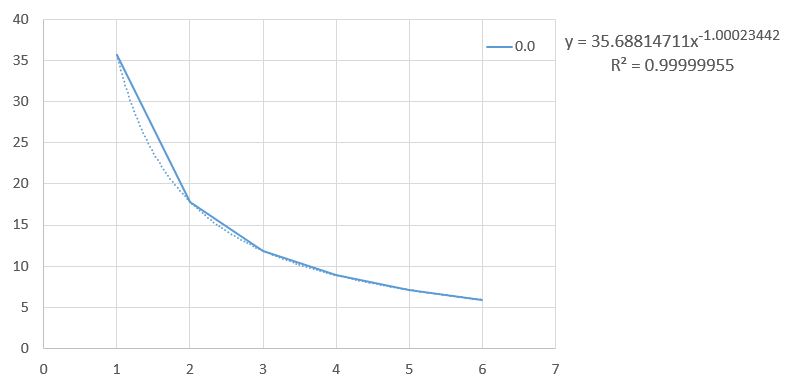}
\caption{Dependence of the middle vector modulus for a larger value of the imaginary part of a complex number , $s=0+8000i$}
\label{fig:v0_8000_6}
\end{figure}

It should be noted that the accuracy of the approximation increases or with a decrease in the number of vectors (fig. \ref{fig:v0_5000_5}), or by increasing a value of the imaginary part of a complex number (fig. \ref{fig:v0_8000_6}), this fact indicates the \textit{asymptotic} dependence of the obtained expressions.
\par
Now approximate the dependence of the coefficients $A$ and $B$ on a value of the real part of a complex number.
\par
We start with the coefficient $B$, because at this stage we will complete its analysis, and for the analysis of the coefficient $A$ additional calculations will be required.
\begin{figure}[ht]
\centering
\includegraphics[scale=0.6]{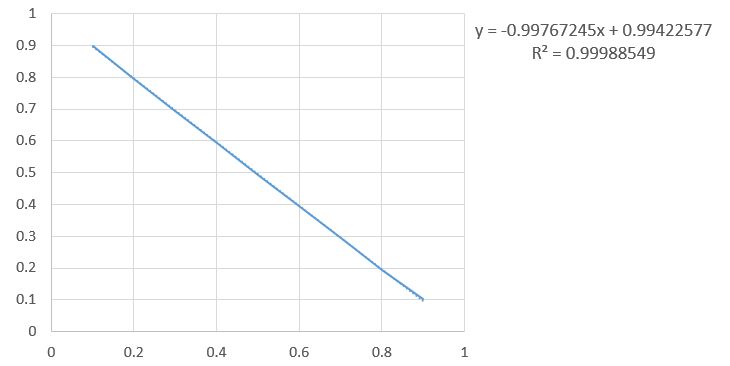}
\caption{Dependence of the coefficient $B$ on the real part of a complex number, $Im(s)=5000$}
\label{fig:factor2_5000}
\end{figure}
\par
The coefficient $B$ has a linear dependence (fig. \ref{fig:factor2_5000}, where $x=\sigma$) from a value of the real part of a complex number.
Therefore, taking into account the identified asymptotic dependence, we can rewrite the expression (\ref{depence_1}) for the modulus of the middle vectors of the Riemann spiral:
\begin{equation}\label{depence_2}|Y_n|=A\frac{1}{n^{1-\sigma}};\end{equation}
\par
\begin{figure}[ht]
\centering
\includegraphics[scale=0.6]{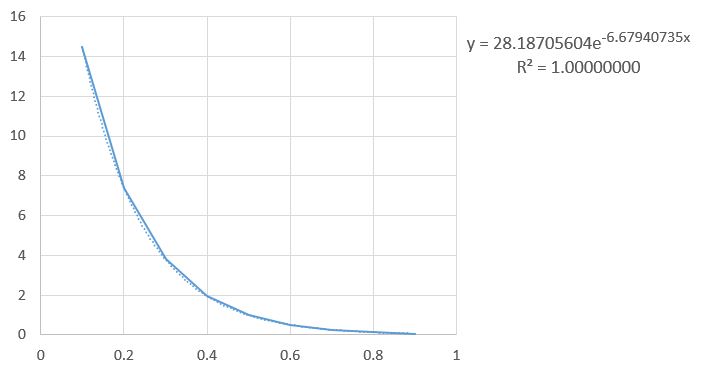}
\caption{Dependence of the coefficient $A$ on the real part of a complex number, $Im(s)=5000$}
\label{fig:factor1_5000}
\end{figure}

The best way to approximate the dependence of the coefficient $A$ (fig. \ref{fig:factor1_5000}, where $x=\sigma$) from a value of the real part of a complex number is the exponent:
\begin{equation}\label{depence_3}A=Ce^{D\sigma};\end{equation}
\par
We first find the ratio of the coefficients $C$ and $D$, for this we calculate $\log(C)$:
\begin{equation}\label{depence_4}2\log(C)=6.67772573;\end{equation}
\par
Taking into account the revealed asymptotic dependence $2\log(C)=D$, we can rewrite the expression (\ref{depence_3}) for the coefficient $A$:
\begin{equation}\label{depence_5}A=Ce^{-2\log(C)\sigma}=e^{log(C)-2\log(C)\sigma}=e^{2log(C)(\frac{1}{2}-\sigma)}=(C^2)^{\frac{1}{2}-\sigma};\end{equation}
\par
Now, taking into account the identified asymptotic dependence $|Y_n|=A=C$ when $\sigma=0$, we calculate the modulus of the first middle vector when $\sigma=0$ for different values of the imaginary part of a complex number in the range from 1000 to 9000 in increments of 1000 and then approximate the dependence of the coefficients $C$ on a value of the imaginary part of a complex number.
\begin{figure}[ht]
\centering
\includegraphics[scale=0.6]{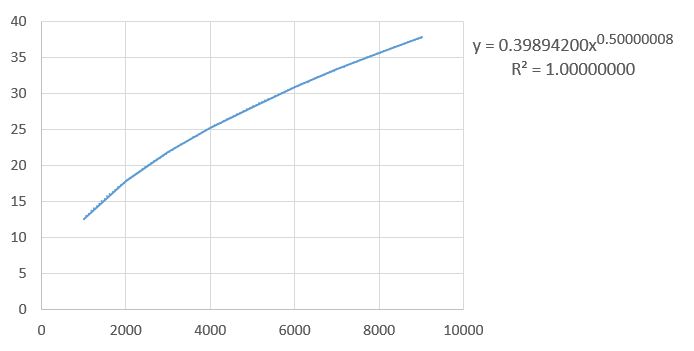}
\caption{Dependence of the coefficient $C$ on the imaginary part of a complex number, $Re(s)=0$}
\label{fig:factor1}
\end{figure}
\par
The best way to approximate the dependence of the coefficient $C$ (Fig. \ref{fig:factor1}, where $x=t$) from a value of the imaginary part of a complex number is a power function:
\begin{equation}\label{depence_6}C=Et^F;\end{equation}
\par
Taking into account the revealed asymptotic dependence consider $$F=\frac{1}{2};$$ \par then we can write the final expression for the modulus of the middle vectors of the Riemann spiral:
\begin{equation}\label{depence_8}|Y_n|=(E^2t)^{\frac{1}{2}-\sigma}\frac{1}{n^{1-\sigma}};\end{equation}
\par
where $E^2=0.159154719364$ some constant, the meaning of which we learn later.
\par
We obtained an asymptotic expression (\ref{depence_8}) for the modulus of the middle vectors of the Riemann spiral, which becomes, as we found out, more precisely when the imaginary part of a complex number increases.
\par
\textit{As a consequence of the asymptotic form of the resulting expression, we can apply it to any middle vector, even if we can no longer calculate its coordinates or the calculated coordinates give an inaccurate value of the middle vector modulus.}
\par
Therefore, we can obtain any quantity of middle vectors necessary to construct \textit{the inverse Riemann spiral.}
\par
So we can get an infinite series, which is given by the middle vectors of the Riemann spiral.
\par
By comparing the expression for the modulus of vectors (\ref{x_n_y_n}) of the Riemann spiral and the expression for the modulus of its middle vectors (\ref{depence_8}), we can assume that the infinite series formed by the middle vectors of the Riemann spiral sets a value of the Riemann zeta function $\zeta(1-s)$.
\par
\textit{So we can assume that values of the Riemann zeta function $\zeta(s)$ and $\zeta(1-s)$ are related through an expression for the middle vectors of the Riemann spiral.}
\par
To complete the derivation of the dependence of the Riemann zeta function $\zeta(s)$ and $\zeta(1-s)$, it is necessary to determine the dependence of the angles between the middle vectors of the Riemann spiral and construct an inverse Riemann spiral whose vectors, as we show further, asymptotically twist around the zero of of the complex plane.
\par
If to determine the coordinates of the centers of Riemann spiral, we used the \textit{reverse points,} to determine the angles between middle vectors of the Riemann spiral, we will use \textit{inflection points.}
We see (fig. \ref{fig:average_vectors}) that the inflection points are not only the points at which the visible untwisting of the vectors is replaced by the twisting, but also the points at which the middle vectors intersect the Riemann spiral.
\par
One can show by computing that the angles between the middle vectors and the Riemann spiral at the intersection points are asymptotically equal to $\pi/4$, then the angles between the middle vectors can be equated to the angles between the Riemann spiral vectors at the inflection points.
\par
As we remember, the inflection points are multiples of $2k\pi$, then using the properties of the logarithm, we can find the angle between the first and any other middle vector of the Riemann spiral:
\begin{equation}\label{depence_9}\beta_k=\alpha_k-\alpha_1=t\log\Big(\frac{t}{2\pi}\Big)-t\log(k)-t\log\Big(\frac{t}{2\pi}\Big)=-t\log(k);\end{equation}
\par
We have obtained an expression that shows that the angles between the first middle vector of the Riemann spiral and any other middle vector are equal in modulus and opposite in sign to the angles between the first vector and corresponding vector of the Riemann spiral, the negative sign shows that the middle vectors have a special kind of symmetry (which we will consider later) with relation to the vectors of value of the Riemann spiral.
\par
\begin{figure}[ht]
\centering
\includegraphics[scale=0.6]{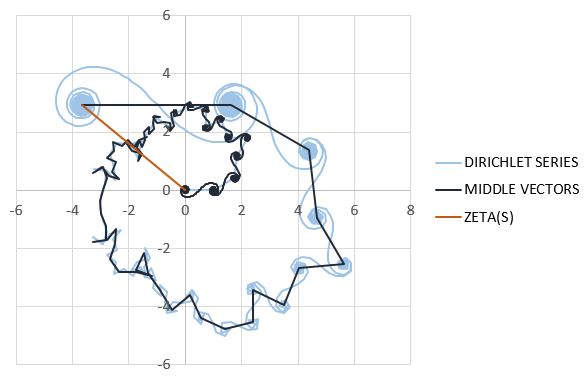}
\caption{Inverse Riemann spiral, $s=0.25+5002.981i$}
\label{fig:reverse_spiral_approx}
\end{figure}

Knowing the coordinates of the first middle vector, which we calculated with sufficient accuracy, we can now find the angles and modulus of the remaining middle vectors using the obtained asymptotic expressions and construct the inverse Riemann spiral (Fig. \ref{fig:reverse_spiral_approx}).
\par
Now at $n=k$ we can write the final formula to calculate the coordinates of the middle vectors $Y_n$:
\begin{equation}\label{h_x_n_h_y_n}\tilde{x}_n(s)=(E^2t)^{\frac{1}{2}-\sigma}\frac{1}{n^{1-\sigma}}\cos(\alpha_1-t\log(n)); \tilde{y}_n(s)=(E^2t)^{\frac{1}{2}-\sigma}\frac{1}{n^{1-\sigma}}\sin(\alpha_1-t\log(n));\end{equation}
\par
Using Euler's formula for complex numbers we write the expression for the middle vectors of the Riemann spiral in exponential form:
\begin{equation}\label{y_n_exp}Y_n(s)=(E^2t)^{\frac{1}{2}-\sigma}\frac{1}{n^{1-\sigma}}e^{-i(\alpha_1-t\log(n))};\end{equation}
\par
Using the rules of analytical geometry, we can obtain the coordinates of the vector corresponding to the partial sum of $\hat{s}_m(s)$ inverse Riemann spiral:
\begin{equation}\label{h_s_m_x_h_s_m_y}\tilde{s}_m(s)_x=\zeta(s)_x-\sum_{n=1}^{m}{\tilde{x}_n(s)}; \tilde{s}_m(s)_y=\zeta(s)_y-\sum_{n=1}^{m}{\tilde{y}_n(s)};\end{equation}
\par
\begin{figure}[ht]
\centering
\includegraphics[scale=0.6]{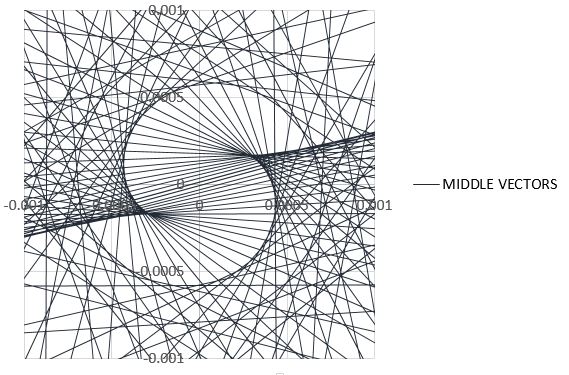}
\caption{Middle vectors of the inverse Riemann spiral near the point of zero, $s=0.25+5002.981i$}
\label{fig:reverse_spiral_zero}
\end{figure}

To verify that the obtained expression for the middle vectors of the Riemann spiral defines the relation values of the Riemann zeta function $\zeta(s)$ and $\zeta(1-s)$, consider the middle vectors of the Riemann spiral near the point of zero (fig. \ref{fig:reverse_spiral_zero}).
\par
We see that the middle of the vectors twist around the zero point, in the same way as the vectors of the Riemann spiral twist around the point with the coordinates of a value of the Riemann zeta function\footnote{This fact was already considered by the forum Stack Exchange \cite{MA}, but nobody try to compute the Riemann zeta function with geometric method, which corresponds to the method of generalized summation Cesaro.}.
\par
Now we can write the final equation that relates values of the Riemann zeta function $\zeta(s)$ and $\zeta(1-s)$ \textit{taking into account the rules of generalized summation of divergent series:}
\begin{equation}\label{zeta_func_approx}\sum_{n=1}^{\infty}{\frac{1}{n^s}}-(E^2t)^{\frac{1}{2}-\sigma}e^{-i\alpha_1}\sum_{n=1}^{\infty}{\frac{1}{n^{1-s}}}=0;\end{equation}
\par
where $\alpha_1$ is the angle of the first middle vector of the Riemann spiral.
\par
\textit{This fact shows the asymptotic form of functional equation and the geometric nature of the Riemann zeta function, based on the significant role of turning points and inflection points of the Riemann spiral, which determine the middle vectors and the inverse Riemann spiral.}
\par
We can try to find an asymptotic expression for the angle $\alpha_1$ of the first middle vector of the Riemann spiral in a geometric way, but this will not significantly improve the result.
\par
We use the results of the analytical theory of numbers to an arithmetic way to find the exact expression for the angle $\alpha_1$ of the first middle vector of the Riemann spiral and to determine the meaning of the constant $E^2=0.159154719364$.
\subsection{Derivation of empirical expression for CHI function}
The functional equation of the Riemann zeta function \cite{TI} has several equivalent entries:
\begin{equation}\label{zeta_func_eq2}\zeta(s)=\chi(s)\zeta(1-s); \end{equation}
\par
where
\begin{equation}\label{chi_eq}\chi(s)=\frac{(2\pi)^s}{2\Gamma(s)\cos(\large\frac{\pi s}{2})}=2^s\pi^{s-1}\sin(\frac{\pi s}{2})\Gamma(1-s)=\pi^{s-\frac{1}{2}}\frac{\Gamma(\frac{1-s}{2})}{\Gamma(\frac{s}{2})};\end{equation}
\par
Comparing (\ref{zeta_func_approx}) and (\ref{zeta_func_eq2}), we get:
\begin{equation}\label{chi_eq_app}\chi(s)=(E^2t)^{\frac{1}{2}-\sigma}e^{-i\alpha_1};\end{equation}
\par
We find the same expression in Titchmarsh \cite{TI}:
\par
in any fixed strip $\alpha\le\sigma\le\beta$, when $t\to \infty$:
\begin{equation}\label{chi_eq_app2}\chi(s)= \Big(\frac{2\pi}{t}\Big)^{(\sigma+it-\frac{1}{2})}e^{i(t+\frac{\pi}{4})}\Big\{1+\mathcal{O}\Big(\frac{1}{t}\Big)\Big\};\end{equation}
\par
We write the expression (\ref{chi_eq_app2}) in the exponential form of a complex number:
\begin{equation}\label{chi_eq_ex}\chi(s)= \Big(\frac{t}{2\pi}\Big)^{(\frac{1}{2}-\sigma)}e^{-i(t(\log{\frac{t}{2\pi}}-1)-\frac{\pi}{4}+\tau(s))};\end{equation}
\par
where
\begin{equation}\label{tau}\tau(s)= \mathcal{O}\Big(\frac{1}{t}\Big);\end{equation}
\par
By matching (\ref{chi_eq_app}) and (\ref{chi_eq_ex}), we obtain the angle of the first middle vector of the Riemann spiral:
\begin{equation}\label{alpha1}\alpha_1=t(\log{\frac{t}{2\pi}}-1)-\frac{\pi}{4}+\tau(s);\end{equation}
\par
And define the meaning of the constant $E^2=0.159154719364$:
\begin{equation}\label{e_2}E^2=\frac{1}{2\pi}=0.159154943091;\end{equation}
\par
Now we can write the asymptotic equation (\ref{zeta_func_approx}) in its final form:
\begin{equation}\label{zeta_func_approx2}\sum_{n=1}^{\infty}{\frac{1}{n^s}}-\Big(\frac{t}{2\pi}\Big)^{\frac{1}{2}-\sigma}e^{-i\alpha_1}\sum_{n=1}^{\infty}{\frac{1}{n^{1-s}}}=0;\end{equation}
\par
We will find the empirical expression for remainder term $\tau(s)$ of the CHI function, which defines the ratio of the modulus and the argument of the exact $\chi(s)$ and the approximate $\tilde\chi(s)$ value of the CHI function:
\begin{equation}\label{chi_eq_rem}\tau(s)=\Delta\varphi_{\chi}+\log\Big(\frac{|\chi(s)|}{|\tilde\chi(s)|}\Big)i;\end{equation}
\par
where
\begin{equation}\label{delta_varphi_chi}\Delta\varphi_{\chi}=Arg(\chi(s))-Arg(\tilde\chi(s));\end{equation}
\par
The exact values\footnote{The exact value will be understood as a value obtained with a given accuracy.} CHI functions we find from the functional equation of the Riemann zeta function, substituting the exact values of the Riemann zeta function:
\begin{equation}\label{chi_eq_ex2}\chi(s)=\frac{\zeta(s)}{\zeta(1-s)};\end{equation}
\par
Approximate values of the CHI function we find from the expression (\ref{chi_eq_ex}), dropping the function $\tau(s)$.
\begin{equation}\label{chi_eq_app3}\tilde\chi(s)= \Big(\frac{t}{2\pi}\Big)^{(\frac{1}{2}-\sigma)}e^{-i(t(\log{\frac{t}{2\pi}}-1)-\frac{\pi}{4})};\end{equation}
\par
The calculation of values of the Riemann zeta function is currently available in different mathematical packages.
To calculate the exact values of the Riemann zeta function we will use the Internet service \cite{ZF}.
\par
We will use 15 significant digits, because this accuracy is enough to analyze the CHI function.
\par
We will calculate values of the Riemann zeta function in the numerator (\ref{chi_eq_ex2}) for values of the real part of a complex number in the range from 0 to 1 in increments of 0.1 and values of the imaginary part of a complex number in the range from 1000 to 9000 in increments of 1000.
Note that the range of values of the real part of the a complex number contains 2m+1 value, which are related by the following relation:
\begin{equation}\label{sigma} 1-\sigma_{k+1}=\sigma_{2m-k+1};\end{equation}
\par
where k varies between 0 and 2m, hence, for k=m:
\begin{equation}\label{sigma_m}\sigma_{m+1} = 0.5;\end{equation}
We use, on the one hand, the property of the Riemann zeta function as a function of a complex variable:
\begin{equation}\label{zeta_conj}\overline{\zeta(1-\sigma+it)}=\zeta(\overline{1-\sigma+it})=\zeta(1-\sigma-it)=\zeta(1-s);\end{equation}
\par
From other side:
\begin{equation}\label{zeta_conj2}\overline{\zeta(1-\sigma+it)}=Re(\zeta(1-\sigma+it))-Im(\zeta(1-\sigma+it))i;\end{equation}
\par
Equate (\ref{zeta_conj}) and (\ref{zeta_conj2}):
\begin{equation}\label{zeta_right}\zeta(1-s)=\zeta(1-\sigma-it)=Re(\zeta(1-\sigma+it))-Im(\zeta(1-\sigma+it))i;\end{equation}
\par
Now we use the relation (\ref{sigma}):
\begin{equation}\label{zeta_right_req}\zeta(1-s_{k+1})=Re(\zeta(\sigma_{2m-k+1}+it))-Im(\zeta(\sigma_{2m-k+1}+it))i=Re(\zeta(s_{2m-k+1})-Im(s_{2m-k+1})i;\end{equation}
\par
We use the resulting formula to compute the Riemann zeta function in the denominator (\ref{chi_eq_ex2}) based on values computed for the numerator.
Then we get an expression for CHI function:
\begin{equation}\label{chi_req}\chi(s_{k+1})=\frac{(Re(\zeta(s_{k+1}))+Im(\zeta(s_{k+1}))i)(Re(\zeta(s_{2m-k+1}))+Im(\zeta(s_{2m-k+1}))i)}{Re(\zeta(s_{2m-k+1}))^2+Im(\zeta(s_{2m-k+1}))^2};\end{equation}
\par
We will open brackets and write separate expressions for the real part of the CHI function:
\begin{equation}\label{chi_req_re}Re(\chi(s_{k+1}))=\frac{Re(\zeta(s_{k+1}))Re(\zeta(s_{2m-k+1}))-Im(\zeta(s_{k+1}))Im(\zeta(s_{2m-k+1}))}{Re(\zeta(s_{2m-k+1}))^2+Im(\zeta(s_{2m-k+1}))^2};\end{equation}
\par
and for imaginary:
\begin{equation}\label{chi_req_im}Im(\chi(s_{k+1}))=\frac{Re(\zeta(s_{k+1}))Im(\zeta(s_{2m-k+1}))+Im(\zeta(s_{k+1}))Re(\zeta(s_{2m-k+1}))}{Re(\zeta(s_{2m-k+1}))^2+Im(\zeta(s_{2m-k+1}))^2};\end{equation}
\par
We will consider this to be the exact value of CHI function, because we can calculate it with a given accuracy.
\par
Now write the expression for the real part of the approximate value of the CHI function:
\begin{equation}\label{chi_req_re_app}Re(\tilde\chi(s_{k+1}))=(\frac{2\pi}{t})^{(\sigma-\frac{1}{2})}\cos(t(\log(\frac{2\pi}{t})+1)+\frac{\pi}{4});\end{equation}
\par
and for imaginary:
\begin{equation}\label{chi_req_im_app}Im(\tilde\chi(s_{k+1}))=(\frac{2\pi}{t})^{(\sigma-\frac{1}{2})}\sin(t(\log(\frac{2\pi}{t})+1)+\frac{\pi}{4});\end{equation}
\begin{figure}[ht]
\centering
\includegraphics[scale=0.6]{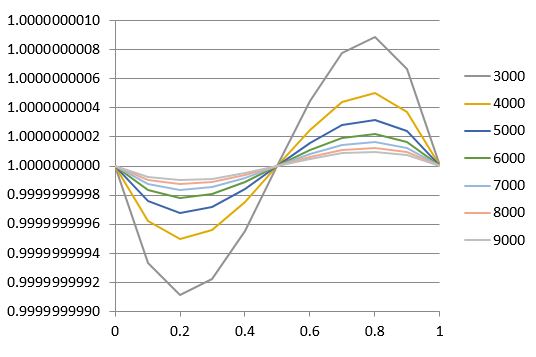}
\caption{The ratio of the exact $|\chi(s)|$ and approximate $|\tilde\chi(s)|$ module of the CHI function}
\label{fig:ratio_chi}
\end{figure}
\par
Calculate the modulus for the exact $|\chi(s)|$ value of CHI function:
\begin{equation}\label{chi_abs}|\chi(s)|=\sqrt{Re(\chi(s))^2+Im(\chi(s))^2};\end{equation}
\par
and for the approximate $|\tilde\chi(s)|$ value of CHI function:
\par
\begin{equation}\label{chi_abs_app}|\tilde\chi(s)|=\sqrt{Re(\tilde\chi(s))^2+Im(\tilde\chi(s))^2};\end{equation}
\par
and the angle between the exact $\chi(s)$ and approximate $\tilde\chi(s)$ value of CHI function:
\begin{equation}\label{chi_angle}\Delta\varphi_{\chi}=Arg(\chi(s))-Arg(\tilde\chi(s))= \arccos(\frac{Re(\chi(s))}{|\chi(s)|})-\arccos(\frac{Re(\tilde\chi(s))}{|\tilde\chi(s)|});\end{equation}
\par
We construct graphs of the ratio of the modulus of the exact $|\chi(s)|$ and approximate $|\tilde\chi(s)|$ values of the CHI function of the real part of a complex number of numbers (fig. \ref{fig:ratio_chi}).
\par
\textit{The graph shows (fig. \ref{fig:ratio_chi}) that the ratio of the modulus of the exact $|\chi(s)|$ and approximate $|\tilde\chi(s)|$ values CHI functions can be taken as 1, therefore, we can say that:}
\begin{equation}\label{tau2}\tau(s)=\Delta\varphi_{\chi};\end{equation}
\par
\begin{figure}[ht]
\centering
\includegraphics[scale=0.6]{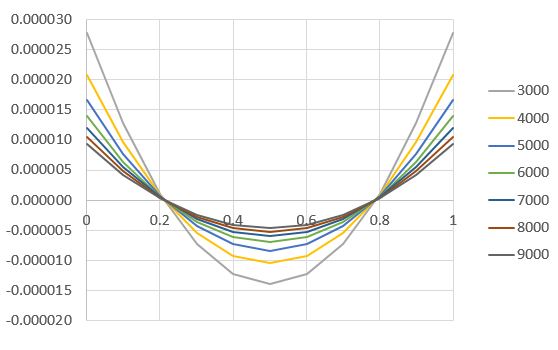}
\caption{The angle $\Delta\varphi_{\chi}$ between the exact $\chi(s)$ and approximate $\tilde\chi(s)$ value of the CHI function (on the real part)}
\label{fig:angle_chi_real}
\end{figure}
\begin{figure}[ht]
\centering
\includegraphics[scale=0.6]{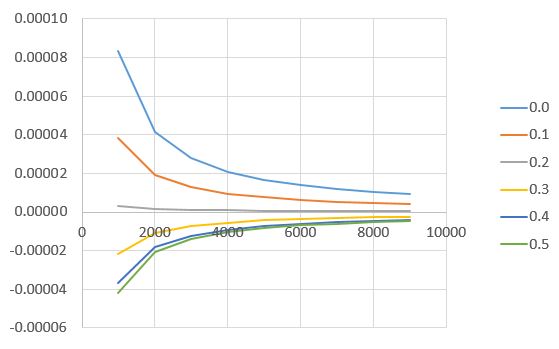}
\caption{The angle $\Delta\varphi_{\chi}$ between the exact $\chi(s)$ and approximate $\tilde\chi(s)$ value of the CHI function (on the imaginary part)}
\label{fig:angle_chi_complex}
\end{figure}

We construct graphs of the dependence of $\Delta\varphi_{\chi}$ on the real part of a complex number (fig. \ref{fig:angle_chi_real}) and from the imaginary part of a complex number (fig. \ref{fig:angle_chi_complex}).
\par
$\tau(s)$ provided (\ref{tau2}) is the argument of remainder term  of the CHI function, hence it is a combination of the arguments of remainder terms of the product $\Gamma(s)\cos(\pi s/2)$ in (\ref{chi_eq}).
\par
The argument of the remainder term $\mu(s)$ of the gamma function can be obtained from an expression we can find in Titchmarsh \cite{TI}:
\begin{equation}\label{gamma_app}\log(\Gamma(\sigma+it))=(\sigma+it-\frac{1}{2})\log{it}-it+\frac{1}{2}\log{2\pi}+\mu(s);\end{equation}
\par
The most significant researches of remainder term of the gamma function can be found in the paper of Riemann \cite{SI} and Gabcke \cite{GA}, they independently and in different ways obtain an expression for the argument of the remainder term of the gamma function when $\sigma=1/2$.
\par
We use the expression explicitly written by Gabcke \cite{GA}:
\begin{equation}\label{mu}\mu(t)=\frac{1}{48t}+\frac{1}{5760t^3}+\frac{1}{80640t^5}+\mathcal{O}(t^{-7});\end{equation}
\par
\begin{figure}[ht]
\centering
\includegraphics[scale=0.6]{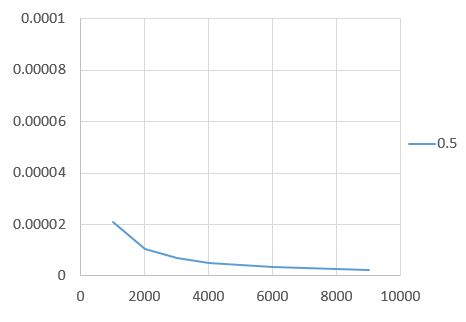}
\caption{The remainder term of the gamma function, $\sigma=1/2$}
\label{fig:remainder_gamma}
\end{figure}
\begin{figure}[ht]
\centering
\includegraphics[scale=0.6]{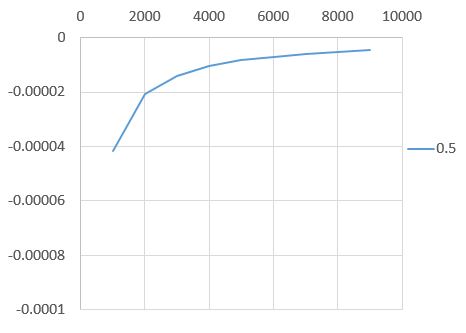}
\caption{The angle $\Delta\varphi_{\chi}$ between the exact $\tilde(s)$ and approximate $\hat\chi(s)$ value of the CHI function, $\sigma=1/2$}
\label{fig:angle_chi_complex_0_5}
\end{figure}

We construct a graph of the dependence of the argument of the remainder term $\mu(t)$ of the gamma function from the imaginary part of a complex number when $\sigma=1/2$ (fig. \ref{fig:remainder_gamma}).
We compare the obtained graph with the graph of dependence $\Delta\varphi_{\chi}$ on the imaginary part of a complex number when $\sigma=1/2$ (fig. \ref{fig:angle_chi_complex_0_5}).
\par
These graphs correspond to each other up to sign and constant value, i.e. the absolute value of the angle between the exact $\chi(s)$ and approximate $\tilde\chi(s)$ value of the CHI function is exactly two times greater than a value of the argument of the remainder term of the gamma function when $\sigma=1/2$.
\par
\textit{But a more significant result is obtained by dividing the angle values $\Delta\varphi_{\chi}$ between the exact $\chi(s)$ and approximate $\tilde\chi(s)$ value of the CHI function by a value of  the argument of remainder term $\mu(t)$ of the gamma function when $\sigma=1/2$.}
\par
The result of this operation we get the functional dependence $\lambda(\sigma)$ values of the argument of the remainder term $\tau(s)$ of CHI function (with different values of the real part of a complex number) from a value of the argument of the remainder term $\mu(t)$ of the gamma function when $\sigma=1/2$.
\begin{figure}[ht]
\centering
\includegraphics[scale=0.6]{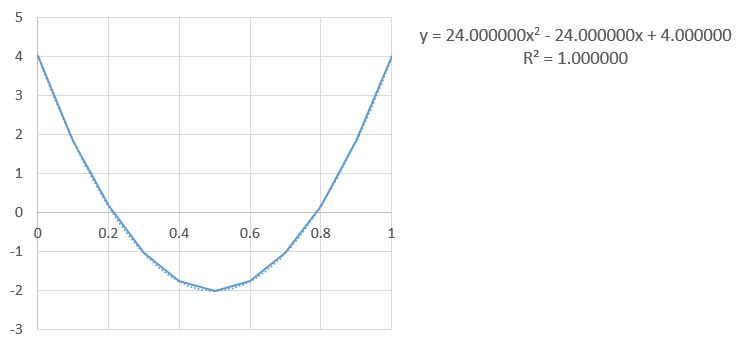}
\caption{Factor $\lambda(\sigma)$ (from the real part of a complex number)}
\label{fig:lambda_func}
\end{figure}
\begin{equation}\label{tau3}\tau(s)=\Delta\varphi_{\chi}=\lambda(\sigma)\mu(t);\end{equation}
\par
The coefficient of $\lambda(\sigma)$ shows which number to multiply a value of the argument of the remainder term $\mu(t)$ of the gamma function in the form of the Riemann-Gabcke to get a value of the argument of the remainder term $\tau(s)$ of CHI function.
\par
We will find an explanation for this paradoxical identity in the further study of expressions that determine a value of the Riemann zeta function.
\par
We obtained an exact expression for the angle $\alpha_1$ of the first middle vector of the Riemann spiral:
\par
\begin{equation}\label{alpha1_ex}\alpha_1=t(\log{\frac{t}{2\pi}}-1)-\frac{\pi}{4}+\lambda(\sigma)\mu(t);\end{equation}
\par
As well as the exact expression for the CHI function:
\begin{equation}\label{chi_eq_ex3}\chi(s)= \Big(\frac{t}{2\pi}\Big)^{(\frac{1}{2}-\sigma)}e^{-i(t(\log{\frac{t}{2\pi}}-1)-\frac{\pi}{4}+\lambda(\sigma)\mu(t))};\end{equation}
\par
Now we can perform the exact construction of the inverse Riemann spiral.
\subsection{Representation of the second approximate equation of the Riemann zeta function by a vector system}
We write the second approximate equation of the Riemann zeta function \cite{HA1} in vector form.
\begin{equation}\label{zeta_eq_2}\zeta(s)=\sum_{n\le x}{\frac{1}{n^s}}+\chi(s)\sum_{n\le y}{\frac{1}{n^{1-s}}}+\mathcal{O}(x^{-\sigma})+\mathcal{O}(|t|^{1/2-\sigma}y^{\sigma-1}); \end{equation}
\begin{equation}\label{zeta_eq_2_cond}0<\sigma <1; 2\pi xy=|t|;\end{equation}
\begin{equation}\label{chi_eq2}\chi(s)=\frac{(2\pi)^s}{2\Gamma(s)\cos(\large\frac{\pi s}{2})}=2^s\pi^{s-1}\sin(\frac{\pi s}{2})\Gamma(1-s);\end{equation}
\par
We use the exponential form of a complex number, then , using Euler's formula, go to the trigonometric form of a complex number and get the coordinates of the vectors.
Put $x=y$, then for $m=\Big[\sqrt{\frac{t}{2\pi}}\Big]$ we get:
\begin{equation}\label{zeta_eq_2_vect}\zeta(s)=\sum_{n=1}^{m}{X_n(s)}+\sum_{n=1}^{m}{Y_n(s)}+R(s);\end{equation}
\begin{description}
\item where
\begin{equation}\label{x_vect}X_n(s)=\frac{1}{n^{s}}=\frac{1}{n^{\sigma}}e^{-it\log(n)}=\frac{1}{n^{\sigma}}(\cos(t\log(n))-i\sin(t\log(n)));\end{equation}
\begin{equation}\label{y_vect}Y_n(s)=\chi(s)\frac{1}{n^{1-s}}=\chi(s)\frac{1}{n^{1-\sigma}}e^{it\log(n)}
=\chi(s)\frac{1}{n^{1-\sigma}}(\cos(t\log(n))+i\sin(t\log(n));\end{equation}
\item $R(s)$ - some function of the complex variable, which we will estimate later using the exact values of $\zeta(s)$ and $\chi(s)$.
\end{description}
\par
The vector system (\ref{zeta_eq_2_vect}) determines a value of $\zeta(s)$ at each interval:
\begin{equation}\label{m_interval}t\in[2\pi m^2,2\pi (m+1)^2);
m=1,2,3...\end{equation}
\begin{figure}[ht]
\centering
\includegraphics[scale=0.6]{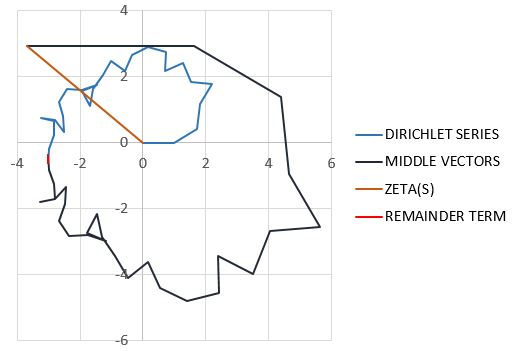}
\caption{Finite vector system, $s=0.25+5002.981i$}
\label{fig:finite_vector_system}
\end{figure}
\par
We see that the first sum (\ref{zeta_eq_2_vect}) corresponds to the vectors of the Riemann spiral, and the second is the middle vector of the Riemann spiral.
\par
Thus, we can explain the geometric meaning of the second approximate equation of the Riemann zeta function.
\begin{figure}[ht]
\centering
\includegraphics[scale=0.6]{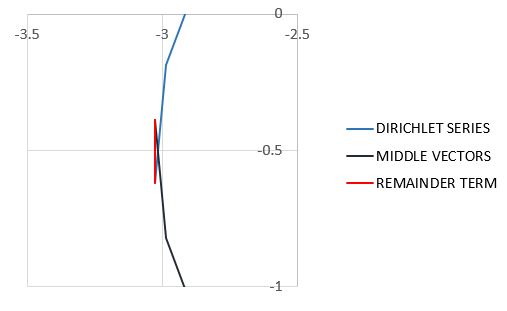}
\caption{Gap of vector system, $s=0.25+5002.981i$}
\label{fig:remainder_gap}
\end{figure}
\par
In the analysis of Riemann spiral we received the quantity of reverse points: $$m=\sqrt{\frac{t}{2\pi}+\frac{1}{4}}$$ of the conditions that between two reverse points is at least one vector of the Riemann spiral.
\par
If we look at the inverse Riemann spiral (fig. \ref{fig:reverse_spiral_approx}), we note that the number of reverse points, which corresponds to one side of the middle vectors of the Riemann spiral, and on the other hand, the vectors of value of the Riemann spiral, which had not yet twist in a convergent and then a divergent spiral.
\par
\textit{If we remove the vectors that twist into spirals, we get the vector system of (fig. \ref{fig:finite_vector_system}), which corresponds to the second approximate equation of the Riemann zeta function.}
\par
We can also explain the geometric meaning of the remainder term of the second approximate equation of the Riemann zeta function.
\par
As the scale of the picture of the vector system increases (fig. \ref{fig:finite_vector_system}), we see a gap between the vectors and the middle vectors of the Riemann spiral (fig. \ref{fig:remainder_gap}), this gap is the remainder term of the second approximate equation of the Riemann zeta function.
\subsection{The axis of symmetry of the vector system of the second approximate equation of the Riemann zeta function - conformal symmetry}
\par
In the process of geometric derivation of the functional equation of the Riemann zeta function we found that the angles (\ref{depence_9}) between the first middle vector of the Riemann spiral and any other middle vector are equal in modulus and opposite in sign to the angles between the first vector and appropriate vector of the Riemann spiral.
\par
The same result (\ref{x_vect}) and (\ref{y_vect}) we obtained when writing the second approximate equation of the Riemann zeta function in vector form.
\par
We show 1) that the angles between any two middle vectors $Y_i$ and $Y_j$ of the Riemann spiral are equal in modulus and opposite in sign, if we measure the angle from the respective vectors, to the angles between the corresponding vectors $X_i$ and $X_j$ of the Riemann spiral and 2) there is a line that has angles equaled modulus and opposite in sign, if we measure the angle from this line, with any pair of corresponding vectors $Y_i$ and $X_i$ of the Riemann spiral (Lemma 1).
\par
Put in accordance with the first middle vector $Y_1$, the two random middle vectors $Y_i$ and $Y_j$ of the Riemann spiral, the vector $X_1$ and two relevant vectors $X_i$ and $X_j$ of the Riemann spiral the segments $A_1A_2$, $A_2A_3$, $A_3A_4$, $A_1'A_2'$ , $A_2'A_3'$ and $A_3'A_4'$ respectively.
\par
\begin{figure}[ht]
\centering
\includegraphics[scale=0.6]{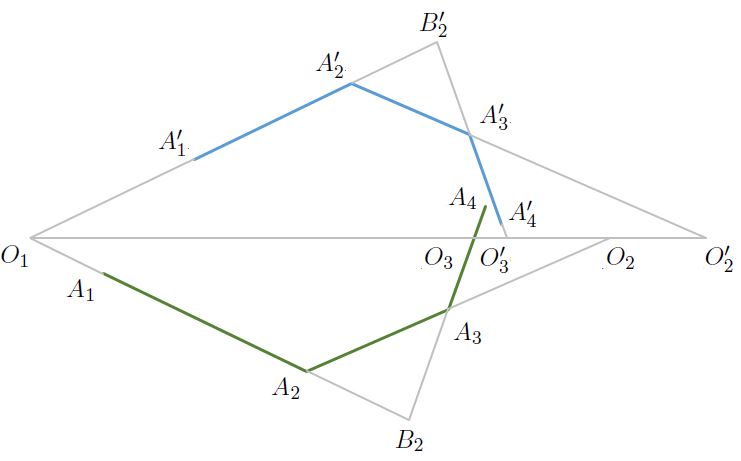}
\caption{Polylines with equal angles}
\label{fig:conformal_symmetry_lemma}
\end{figure}

Consider (fig. \ref{fig:conformal_symmetry_lemma}) two polyline formed by vertices $A_1A_2A_3A_4$ and $A_1'A_2'A_3'A_4 '$ respectively, and oriented arbitrarily. Then edges $A_2A_3$ and $A_3A_4$ have angles with the edge $A_1A_2$ is equal in modulus and opposite in sign to the angles, which have edge $A_2'A_3'$ and $A_3'A_4'$ respectively, with the edge $A_1'A_2'$, if we measure angles from the edges $A_1A_2$ and $A_1'A_2'$ respectively.
\par
1) We show that the angle between edges $A_2A_3$ and $A_3A_4$ is equal in modulus and opposite in sign to the angle between the edges $A_2'A_3'$ and $A_3'A_4 '$, if we measure the angle from the corresponding edges, for example, from $A_3A_4$ and $A_3'A_4'$.
\par
A) Consider first (fig. \ref{fig:conformal_symmetry_lemma}) case when edges $A_1A_2$, $A_3A_4$ and edges $A_1'A_2'$, $A_3'A_4'$ respectively are not parallel to each other.
\par
Continue edges $A_1A_2$, $A_3A_4$ and edges $A_1'A_2'$, $A_3'A_4'$ to the intersection, we obtain the triangles $A_2B_2A_3$ and $A_2'B_2'A_3'$ respectively (fig. \ref{fig:conformal_symmetry_lemma}).
\par
These triangles are congruent by two angles, because the angles $A_1B_2A_4$ and $A_1'B_2'A_4'$ are equal in modulus, as the angle of the edges $A_3A_4$ and $A_3'A_4'$ with edges $A_1A_2$ and $A_1'A_2'$ respectively, and the angles $B_2A_2A_3$ and $B_2'A_2'A_3'$ are adjacent angles $A_1A_2A_3$ and $A_1'A_2'A_3'$ respectively which are also equal in modulus, as the angle of the edges $A_2A_3$ and $A_2'A_3'$ with edges $A_1A_2$ and $A_1'A_2'$ respectively.
\par
Hence, the angles $A_2A_3A_4$ and $A_2'A_3'A_4'$ are equal in modulus as the angles adjacent to the angles $A_2A_3B_2$ and $A_2'A_3'B_2'$ respectively, which are equal as corresponding angles of congruent triangles.
\par
\begin{figure}[ht]
\centering
\includegraphics[scale=0.6]{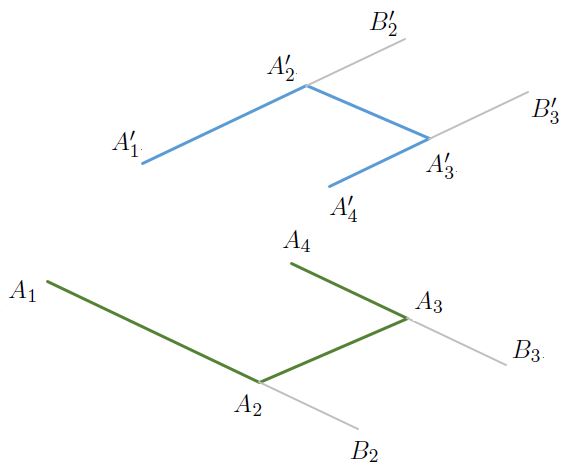}
\caption{Polylines with equal angles, parallel edges}
\label{fig:conformal_symmetry_2_lemma}
\end{figure}

B) Now consider (fig. \ref{fig:conformal_symmetry_2_lemma}) case when edges $A_1A_2$, $A_3A_4$ and edges $A_1'A_2'$, $A_3'A_4'$ respectively parallel to each other.
\par
The continue of edges $A_1A_2$, $A_3A_4$ and edges $A_1'A_2'$, $A_3'A_4'$ respectively do not intersect, because they form parallel lines $A_1B_2$, $A_4B_3$ and $A_1'B_2'$, $A_4'B_3'$ respectively.
\par
Thus edges $A_2A_3$ and $A_2'A_3'$ are intersecting lines of those parallel lines, respectively.
\par
Hence, the angles $A_2A_3A_4$ and $A_2'A_3'A_4'$ are equal in modulus as corresponding angles at intersecting lines of two parallel lines because the angles $A_1A_2A_3$ and $A_1'A_2'A_3'$ are equal in modulus, as the angle of the edges $A_2A_3$ and $A_2'A_3'$ with edges $A_1A_2$ and $A_1'A_2'$ respectively.
\par
If we measure the angle $A_2A_3A_4$ from edge $A_3A_4$, it is necessary to count its anti-clockwise, then it has a positive sign.
\par
If we measure the angle $A_2'A_3'A_4'$ from edge $A_3'A_4'$, it is necessary to count its clockwise, then it has a negative sign.
\par
Therefore, the angle between the edges $A_2A_3$ and $A_3A_4$ equal in modulus and opposite in sign to the angle between the edges $A_2'A_3'$ and $A_3'A_4'$, if we measure the angle from the edges $A_3A_4$ and $A_3'A_4'$ respectively.
\par
2) Now we continue edges $A_1A_2$ and $A_1'A_2'$ to the intersection and divide the angle $A_2O_1A_2'$ into two equal angles, we will get a line $O_1O_3$, which has equal in modulus and opposite in sign angles, if we measure the angle from this line, to edges $A_1A_2$ and $A_1'A_2'$ (fig. \ref{fig:conformal_symmetry_lemma}).
\par
We show that line $O_1O_3$ also is equal in modulus and opposite in sign to angles, if we measure the angle from this line, with edges $A_2A_3$, $A_3A_4$ and $A_2'A_3'$, $A_3'A_4'$ respectively.
\par
A) Consider first (fig. \ref{fig:conformal_symmetry_lemma}) case when edges $A_2A_3$ and $A_2'A_3'$ are not parallel to each other (sum of angles$A_2O_1O_3$, $A_1A_2A_3$ and $A_2'O_1O_3$, $A_1'A_2'A_3'$ is not equal to $\pi$).
\par
Continue edges $A_2 A_3$, $A_2'A_3'$ to the intersection with the line $O_1 O_3$, we get the triangles $O_1A_2 O_2$ and $O_1A_2'O_2'$ respectively.
\par
These triangles are congruent by two angles, because the angles $A_2O_1O_3$ and $A_2'O_1O_3$ are equal in modulus by build, and the angles $A_1A_2A_3$ and $A_1'A_2'A_3'$ are equal in modulus, as the angle of the edges $A_1A_2$ and $A_1'A_2'$ with edges $A_2A_3$ and $A_2'A_3'$ respectively.
\par
Hence, angles $A_2O_2O_1$ and $A_2'O_2'O_1$ are equal in modulus as the corresponding angles of the congruent triangles.
\par
\begin{figure}[ht]
\centering
\includegraphics[scale=0.6]{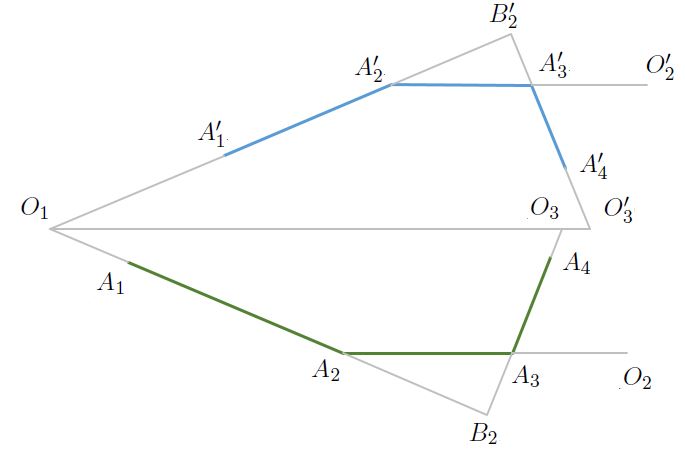}
\caption{Polylines with equal angles, sum of angles $\pi$}
\label{fig:conformal_symmetry_3_lemma}
\end{figure}

B) Now consider (fig. \ref{fig:conformal_symmetry_3_lemma}) case when edges $A_2A_3$ and $A_2'A_3'$ are parallel to each other (sum of angles$A_2O_1O_3$, $A_1A_2A_3$ and $A_2'O_1O_3$, $A_1'A_2'A_3'$ is equal to $\pi$).
\par
In this case, lines $O_1A_2$ and $O_1A_2'$ are intersecting two parallel lines $O_1O_3$, $A_2O_2$ and $O_1O_3$, $A_2'O_2'$ respectively, since the angles $A_2O_1O_3$ and $A_2'O_1O_3$ are equal in modulus by build, and the angles $A_1A_2A_3$ and $A_1'A_2'A_3'$ are equal in modulus as the angle of the edges $A_1A_2$ and $A_1'A_2'$ with edges $A_2A_3$ and $A_2'A_3'$ respectively, and the sum of the respective angles at intersecting lines is equal to $\pi$.
\par
Therefore, the angles between the line segments $A_2A_3$ and $A_2'A_3'$ and line $O_1O_3$ is equal to zero because $A_2A_3$ and $A_2'A_3'$ are parallel to line $O_1O_3$.
\par
Continue edges $A_3A_4$, $A_3'A_4'$ to the intersection with the straight line $O_1O_3$, get the triangles $O1B_2O_3$ and $O_1B_2'O_3'$ respectively (fig. \ref{fig:conformal_symmetry_lemma}).
\par
These triangles are congruent by two angles, because the angles $A_2O_1O_3$ and $A_2'O_1O_3$ are equal in modulus by build, and the angles $A_1B_2A_4$ and $A_1'B_2'A_4'$ are equal in modulus, as the angle of the edges $A_1A_2$ and $A_1'A_2'$ with edges $A_3A_4$ and $A_3'A_4'$ respectively.
\par
Hence, angles $B_2O_3O_1$ and $B_2'O_3'O_1$ are equal in modulus as the corresponding angles of the congruent triangles.
\par
If we measure the angles $A_2O_2O_1$ and $B_2O_3O_1$ from a line $O_1O_3$, they need to count anti-clockwise, then either have a positive sign.
\par
If we measure the angles $A_2'O_2'O_1$ and $B_2'O_3'O_1$ from a line $O_1O_3$, they need to count clockwise, then either have a negative sign.
\par
Therefore, line $O_1O_3$ has equal in modulus and opposite in sign angles, if we measure the angle from this line, with edges $A_1A_2$, $A_2A_3$, $A_3A_4$ and $A_1'A_2'$, $A_2'A_3'$, $A_3'A_4'$ respectively. $\square$
\par
According to the Lemma 1, the vector system of the second approximate equation of the Riemann zeta function has \textit{a special kind of symmetry} when there is a line that has angles equal in modulus and opposite in sign, if we measure its from this line, with any pair of corresponding vectors $Y_i$ and $X_i$ of the Riemann spiral.
\par
It should be noted that this symmetry of angles is kept when $\sigma=1/2$ when, as we will show later, the vector system of the second approximate equation of the Riemann zeta function has \textit{mirror symmetry.}
\par
To distinguish these two types of symmetry, we give a name to a special kind of symmetry of the vector system of the second approximate equation of the Riemann zeta function by analogy with the conformal transformation in which the angles are kept.
\par
\textit{Conformal symmetry - a special kind of symmetry in which there is a line that has equal modulus and opposite sign angles, if we measure the angle from this line, with any pair of corresponding segments.}
\par
The angle $\hat\varphi_M$ of the axis of mirror symmetry is equal to the angle $\varphi_M$ of the axis of conformal symmetry
\begin{equation}\label{phi_m}\hat\varphi_M=\varphi_M=\frac{Arg(\chi(s))}{2}+\frac{\pi}{2};\end{equation}
\par
In other words, it is the same line if we draw the axis of conformal symmetry at the same distance from the end of the first middle vector $Y_1$ and from the end of the first vector $X_1$ of the Riemann spiral.
\subsection{Mirror symmetry of the vector system of the second approximate equation of the Riemann zeta function}
In 1932 Siegel published notes of Riemann \cite{SI} in which Riemann, unlike Hardy and Littlewood, represented the remainder term of the second approximate equation of the Riemann zeta function explicitly:
\begin{equation}\label{zeta_eq_2_zi}\zeta(s)=\sum_{l=1}^{m}{l^{-s}}+\frac{(2\pi)^s}{2\Gamma(s)\cos(\frac{\pi s}{2})}\sum_{l=1}^{m}{l^{s-1}}+(-1)^{m-1}\frac{(2\pi) ^{\frac{s+1}{2}}}{\Gamma(s)}t^{\frac{s-1}{2}}e^{ \frac{\pi is}{2}- \frac{ti}{2}- \frac{\pi i}{8}}\mathcal{S}; 
\end{equation}
\begin{equation}\label{rem_sum}\mathcal{S}=\sum_{0\le 2r\le k\le n-1}{\frac{2^{-k}i^{r-k}k!}{r!(k-2r)!}a_kF^{(k-2r)}(\delta)}+\mathcal{O}\Big(\big(\frac{3n}{t}\big)^{\frac{n}{6}}\Big);\end{equation}
\begin{equation}\label{rem_param}n\le 2\cdot 10^{-8}t; 
m=\Big[\sqrt{\frac{t}{2\pi}}\Big]; \\
\delta=\sqrt{t}-(m+\frac{1}{2})\sqrt{2\pi};\end{equation}
\begin{equation}\label{rem_func}F(u) =\frac{\cos{(u^2+\frac{3\pi}{8})}}{\cos{(\sqrt{2\pi}u)}};\end{equation}
\par
We use an approximate expression (\ref{gamma_app}) for the gamma function to write the expression for the remainder term of the second approximate equation of the Riemann zeta function in exponential form.
\begin{equation}\label{rem_exp}R(s) = (-1)^{m-1}\Big(\frac{t}{2\pi}\Big)^{-\frac{\sigma}{2}}e^{-i[\frac{t}{2}(\log{\frac{t}{2\pi}}-1)-\frac{\pi}{8}+\mu(s)]}\mathcal{S},\end{equation}
\par
Given an estimate for the sum of $\mathcal{S}$, which can be found, for example, Titchmarsh \cite{TI}:
\begin{equation}\label{sum_app}\mathcal{S}=\frac{\cos{(\delta^2+\frac{3\pi}{8})}}{\cos{(\sqrt{2\pi}\delta)}}+\mathcal{O}(t^{-\frac{1}{2}});\end{equation}
we can identify the expression for the argument of the remainder term of the second approximate equation of the Riemann zeta function:
\begin{equation}\label{arg_rem}Arg(R(s))=-(\frac{t}{2}(\log{\frac{t}{2\pi}}-1)-\frac{\pi}{8}+\mu(s));\end{equation}
\par
Taking into account the expression (\ref{mu}) for the argument of the remainder term of the gamma function when $\sigma=1/2$, we obtain an expression for the argument of the remainder term of the second approximate equation of the Riemann zeta function on the critical line.
\begin{equation}\label{arg_rem_2}Arg(\hat R(s))=-(\frac{t}{2}(\log{\frac{t}{2\pi}}-1)-\frac{\pi}{8}+\mu(t));\end{equation}
\par
In the derivation of exact expressions for the CHI function, we found the expression (\ref{tau3}), which shows that when $\sigma=1/2$ value of the argument of remainder term $\tau(s)$ of CHI function is exactly twice the argument of the remainder term $\mu(t)$ of gamma functions (\ref{mu}) in the form of the Riemann-Gabcke.
\begin{equation}\label{arg_chi}Arg(\hat\chi(s))=-(t(\log{\frac{t}{2\pi}}-1)-\frac{\pi}{4}+2\mu(t));\end{equation}
\par
Comparing the expression (\ref{arg_rem_2}) for the argument, the remainder term of the second approximate equation of the Riemann zeta function of and the expression (\ref{arg_chi}) for the argument CHI-function on the critical line, we find a fundamental property of the vector of the remainder term of the second approximate equation of the Riemann zeta function:
\begin{equation}\label{arg_rem_3}Arg(\hat R(s))=\frac{Agr(\hat\chi(s))}{2};\end{equation}
\par
\textit{On the critical line, the argument of the remainder term of the second approximate equation of the Riemann zeta function is exactly half the argument of the CHI function.}
\par
Comparing (\ref{phi_m}) and (\ref{arg_rem_3}) we see that when $\sigma=1/2$, the vector of the remainder term $\hat R(s)$ is perpendicular to the axis of symmetry of the vector system of the second approximate equation of the zeta function of Riemann:
\begin{equation}\label{arg_rem_4}Arg(\hat R(s))=\varphi_L;\end{equation}
\par
As we will show later, this fact is fundamental in the existence of non-trivial zeros of the Riemann zeta function on the critical line.
\par
When changing the imaginary part of a complex number, the vectors of the vector system of the second approximate equation of the Riemann zeta function can occupy an any position in the entire range of angles $[0, 2\pi]$, in consequence of which they form a polyline with self-intersections (fig. \ref{fig:finite_vector_system}), which complicates the analysis of this vector system.
\par
\begin{figure}[ht]
\centering
\includegraphics[scale=0.6]{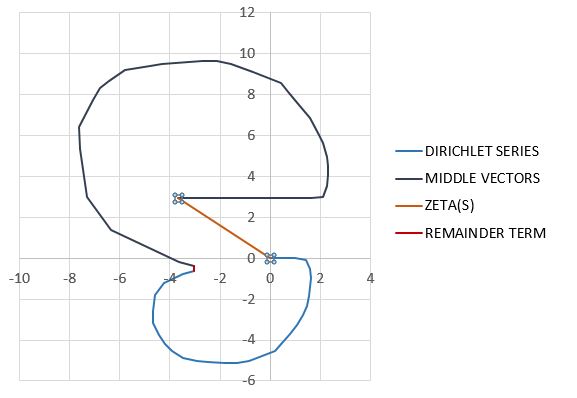}
\caption{Permutation of vectors of the Riemann spiral vector system, $s=0.25+5002.981i$}
\label{fig:finite_vector_system_perm}
\end{figure}

To obtain the polyline formed by the vectors $X_n$ and the middle vectors of the Riemann spiral $Y_n$, without self-intersections, the vectors can be ordered by a value of the angle, then they will form a polyline, which has no intersections (fig. \ref{fig:finite_vector_system_perm}).
\par
We determine the properties of the vector system of the second approximate equation of the Riemann zeta function in the permutation of vectors.
\par
As it is known from analytical geometry, the sum of vectors conforms the permutation law, i.e. it does not change when the vectors are permuted.
\par
Thus, the permutation of the vectors of the vector system of the second approximate equation of the Riemann zeta function does not affect a value of the Riemann zeta function.
\par
Since conformal symmetry, by definition, depends only on angles and does not depend on the actual position of the segments, then the permutation of the vectors of the vector system of the second approximate equation of the Riemann zeta function, conformal symmetry is kept because the angles between the lines formed by the vectors and axis of symmetry do not change.
\par
While, for mirror symmetry, the permutation of vectors forms a new pair of vertices and it is necessary to determine that they are on the same line perpendicular to the axis of symmetry and at the same distance from the axis of symmetry.
\begin{figure}[ht]
\centering
\includegraphics[scale=0.6]{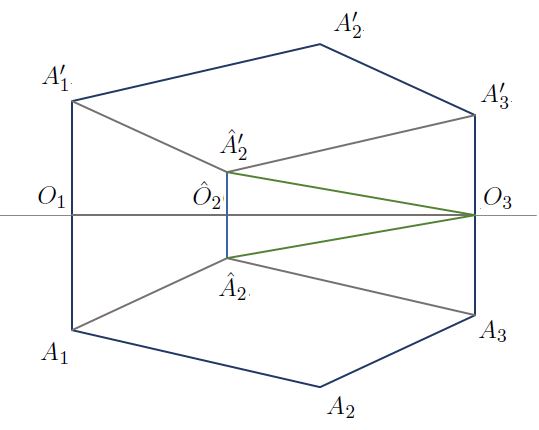}
\caption{Permutation vectors, new pair of vertices}
\label{fig:vector_permutation}
\end{figure}
\par
Consider (fig. \ref{fig:vector_permutation}) symmetrical polygon $S_1$ formed by the vertices $A_1A_2A_3$ $A_3'A_2'A_1'$.
\par
The axis of symmetry $O_1O_3$ divides segments $A_1'A_1$, $A_2'A_2$ and $A_3'A_3$ into equal segments, because vertices $A_1'$ and $A_1$, $A_2'$ and $A_2$, $A_3'$ and $A_3$ is mirror symmetrical relative to the axis of symmetry $O_1O_3$.
\par
Change the edges $A_1'A_2'$, $A_2'A_3'$ and $A_1A_2$, $A_2A_3$ of their places, we get two new vertices $\hat A_2'$ and $\hat A_2$.
\par
We connect vertices $\hat A_2'$, $O_3$ and $\hat A_2$, $O_3$, we get the triangles $T_1'$ and $T_1$ are formed respectively the vertices of $\hat A_2'A_3'O_3$ and $\hat A_2A_3O_3$.
\par
We show that the triangles $T_1'$ and $T_1$ are equal.
\par
Parallelograms $A_1'A_2'A_3'\hat A_2'$ and $A_1A_2A_3\hat A_2$ are equal by build, therefore, angles $A_2'A_3'\hat A_2'$ and $A_2A_3\hat A_2$ are equal;
\par
In the source polygon $S_1$, the angles $A_2'A_3'O_3$ and $A_2A_3O_3$ are equal, hence angles $\hat A_2'A_3'O_3$ and $\hat A_2A_3O_3$ are equal as the difference of equal angles;
\par
Then the triangles $T_1'$ and $T_2$ are equal by the equality of the two edges $\hat A_2'A_3'$, $\hat A_2A_3$ and $A_3'O_3$, $A_3O_3$ and the angle between them $\hat A_2'A_3'O_3$ and $\hat A_2A_3O_3$;
\par
Connect vertices $\hat A_2'$ and $\hat A_2$, we get the intersection of $\hat O_2$ segments $\hat A_2'\hat A_2$ and the axis of symmetry $O_1O_3$.
\par
The angles $\hat A_2'O_3A_3'$ and $\hat A_2O_3A_3$ are equal as the corresponding angles of equal triangles $T_1'$ and $T_1$, hence angles $\hat A_2'O_3\hat O_2$ and $\hat A_2O_3\hat O_2$ are equal;
\par
Then the triangle $\hat A_2'O_3\hat A_2$ is isosceles and the segment $\hat O_2O_3$ is its height, because the bisector in an equilateral triangle is its height.
\par
Therefore, the vertices $\hat A_2'$ and $\hat A_2$ lie on the line perpendicular to the axis of symmetry $O_1O_3$ and they have the same distance from the axis of symmetry $O_1O_3$.
\par
Thus, when the corresponding edges of the symmetric polygon are permuted, the mirror symmetry is kept (Lemma 2). $\square$
\par
Now we define the property of the vector system of the second approximate equation of the Riemann zeta function when $\sigma=1/2$.
\par
\begin{figure}[ht]
\centering
\includegraphics[scale=0.6]{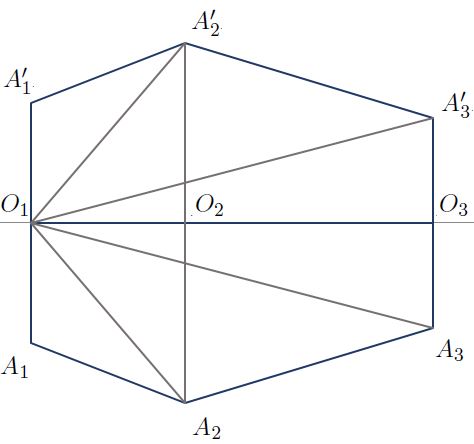}
\caption{Equal vectors possessing conformal symmetry}
\label{fig:conformal_vectors}
\end{figure}

Consider (fig. \ref{fig:conformal_vectors}) polygon, formed by the vertices $A_1A_2A_3A_3'A_2'A_1' $, whose edges $A_1'A_2'$, $A_2'A_3'$ and $A_1A_2$, $A_2A_3$ are equal, have conformal symmetry with relation to the axis of symmetry $O_1O_3$ and vertices $A_1'$ and $A_1$ is mirror symmetrical relative to the axis of symmetry $O_1O_3$.
\par
We show that the vertices $A_2'$, $A_3'$ and $A_2$, $A_3$ respectively are mirror symmetric about the axis of symmetry $O_1O_3$ (Lemma 3).
\par
The axis of symmetry $O_1O_3$ divides the segment $A_1'A_1$ into equal segments and segment $A_1'A_1$ is perpendicular to the axis of symmetry $O_1O_3$.
\par
The angles $A_1A_2A_3$ and $A_1'A_2'A_3'$ are equal by Lemma 1, since edges $A_1A_2$, $A_2A_3$ and $A_1'A_2'$, $A_2'A_3'$ respectively have conformal symmetry.
\par
Construct segments $O_1A_2'$, $O_1A_3'$ and $O_1A_2$, $O_1A_3$.
\par
Triangles $O_1A_1'A_2'$ and $O_1A_1A_2$ are equal by the equality of the two sides $O_1A_1'$, $A_1'A_2'$ and $O_1A_1$, $A_1A_2$ respectively and the angle between them $O_1A_1'A_2'$ and $O_1A_1A_2$.
\par
The angles $O_1A_2'A_3'$ and $O_1A_2A_3$ are equal as parts of the equal angles $A_1'A_2'A_3'$ and $A_1A_2A_3$ since angles $A_1'A_2'O_1$ and $A_1A_2O_1$ are equal as the corresponding angles of equal triangles.
\par
Triangles $O_1A_2'A_3'$ and $O_1A_2A_3$ are equal by the equality of the two sides $O_1A_2'$, $A_2'A_3'$ and $O_1A_2$, $A_2A_3$ respectively and the angle between them $O_1A_2'A_3'$ and $O_1A_2A_3$.
\par
The angles $O_3O_1A_2'$ and $O_3O_1A_2$ are equal as parts os the equal angles $A_1'O_1O_3$ and $A_1O_1O_3$ since angles $A_1'O_1A_2'$ and $A_1O_1A_2$ are equal as the corresponding angles of equal triangles.
\par
The angles $O_3O_1A_3'$ and $O_3O_1A_3$ are equal as parts of the equal angles $A_1'O_1O_3$ and $A_1O_1O_3$ since angles $A_1'O_1A_2'$, $A_1O_1A_2$ and $A_2'O_1A_3'$, $A_2O_1A_3$ respectively are equal as the corresponding angles of equal triangles.
\par
Triangles $A_3'O_1A_3$ and $A_2'O_1A_2$ are isosceles and segments $O_1O_2$ and $O_1O_3$ are respectively their height, hence the segments $O_1O_2$ and $O_1O_3$ divide segments $A_2'A_2$ and $A_3 A_3'$ into equal parts and segments $A_2'A_2$ and $A_3 A_3'$ are perpendilular axis of symmetry $O_1O_3$. $\square$
\par
\begin{figure}[ht]
\centering
\includegraphics[scale=0.6]{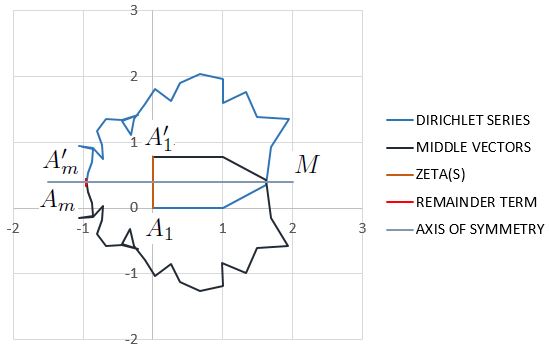}
\caption{Mirror symmetry of the Riemann spiral vector system, $s=0.5+5002.981i$}
\label{fig:mirror_symmetry}
\end{figure}

We construct a polyline formed by the vectors of the second approximate equation of the Riemann zeta function when $\sigma=1/2$ (fig. \ref{fig:mirror_symmetry}).
\par
In accordance with (\ref{x_vect}) and (\ref{y_vect}) for $X_n$ and $Y_n$ respectively all segments $A_1'A_2', A_2'A_3', ... \\A_{m-2}'A_{m-1}', A_{m-1}'A_m'$ and $A_1A_2, A_2A_3, ... A_{m-2}A_{m-1}, A_{m-1}A_m$ when $\sigma=1/2$ are equal.
\par
We draw the axis of symmetry $M$ of the vector system of the second approximate equation of the Riemann zeta function through the middle of the segment formed by the vector $\hat R(s)$ of the remainder term when $\sigma=1/2$.
\par
We obtain two mirror-symmetric vertices $A_m ' $ and $A_m$ besause when $\sigma=1/2$ the vector $\hat R(s)$ of the remainder term is perpendicular to the axis of symmetry $M$ of the vector system of the second approximate equation of the Riemann zeta function.
\par
\textit{According to Lemma 3, the vector system of the second approximate equation of the Riemann zeta function when $\sigma=1/2$ has mirror symmetry}
\footnote{We later show that the mirror symmetry of the vector system of the second approximate equation of the Riemann zeta function is also determined by the argument of the Riemann zeta function when $\sigma=1/2$.}.
\par
Corollary 1. The vector $A_1'A_1$ corresponds to the vector of a value of the Riemann zeta function when $\sigma=1/2$, therefore, when $\sigma=1/2$, the argument of the Riemann zeta function up to the sign corresponds to the direction of the normal $L$ to the axis of symmetry of the vector system of the second approximate equation of the Riemann zeta function.
\par
Corollary 2. Consider the projection of the vector system of the second approximate equation of the Riemann zeta function when $\sigma=1/2$ on the axis of symmetry $M$ of this vector system (fig. \ref{fig:mirror_symmetry}).
\par
Segments $A_1'A_1$ and $A_m'A_m$ is perpendicular to the axis of symmetry of $M$, hence, their projection on axis of symmetry $M$ equal to zero.
\par
The projections of the vectors $X_n$ and $Y_n$ on the axis of symmetry $M$ are equal in modulus and opposite sign, hence
\begin{equation}\label{x_n_y_n_m}(\sum_{n=1}^{m}{X_n})_M+(\sum_{n=1}^{m}{Y_n})_M=0;\end{equation}
\par
Therefore
\begin{equation}\label{zeta_app_m}(\sum_{n=1}^{m}{X_n})_M+(\sum_{n=1}^{m}{Y_n})_M+\hat R(s)_M=0;\end{equation}
\par
and
\begin{equation}\label{zeta_m}\hat \zeta(s)_M=0;\end{equation}
\par
Corollary 3. Consider the projection of the vector system of the second approximate equation of the Riemann zeta function when $\sigma=1/2$ on the normal $L$ to the axis of symmetry $M$ of this vector system (fig. \ref{fig:mirror_symmetry}).
\par
According to the rules of vector summation
\begin{equation}\label{zeta_app_l}\hat \zeta(s)=\hat \zeta(s)_L=(\sum_{n=1}^{m}{X_n})_L+(\sum_{n=1}^{m}{Y_n})_L+\hat R(s)_L;\end{equation}
\par
\textit{This expression, as we will show later, called the Riemann-Siegel formula, is used to find the non-trivial zeros of the Riemann zeta function on the critical line.}
\par
It should be noted that, while, vector $X_1$ remains fixed relative to the axes $x=Re(s)$ and $y=Im(s)$, relative to the normal of $L$ to the axis of symmetry and the axis of symmetry $M$ of the vector system of the second approximate equation of the Riemann zeta function vectors $X_n$ and $Y_n$ are rotated, in accordance with (\ref{x_vect}) and (\ref{y_vect}) towards each other with equal speeds and the angles of the vector of the remainder term $R(s)$ remains fixed when $\sigma=1/2$.
\par
This rotation of the vectors $X_n$ and $Y_n$ leads to the cyclic behavior of the projection $\hat\zeta(s)_L$ of the Riemann zeta function on the normal $L$ to the axis of symmetry of the vector system of the second approximate equation of the Riemann zeta function on \textit{the critical line,} i.e. $\hat\zeta(s)_L=\zeta(1/2+it)_L$ alternately takes the maximum positive and maximum negative value, therefore, $\hat\zeta(s)_L$ cyclically takes a value of zero.
\par
Projections $\zeta(s)_L$ and $\zeta(s)_M$ of the Riemann zeta function respectively on the normal $L$ to the axis of symmetry and on the axis of symmetry $M$ of the vector system of the second approximate equation of the Riemann zeta function in the critical strip when $\sigma\ne 1/2$ have the same cyclic behavior, i.e. $\zeta(s)_L$ and $\zeta(s)_M$ alternately take the maximum positive and maximum negative value, hence $\zeta(s)_L$ and $\zeta(s)_M$ is cyclically takes a value of zero, and, as we will show later, if $\zeta(s)_L$ or $\zeta(s)_M$ are set to zero for any value of $\sigma+it$, then they take a value of zero for value $1-\sigma+it$ also.
\subsection{Non-trivial zeros of the Riemann zeta function}
Using the vector equation of the Riemann zeta function (\ref{zeta_eq_2_vect}), we can obtain the vector equation of the non-trivial zeros of the Riemann zeta function.
\begin{equation}\label{zeta_eq_2_zero}\sum_{n=1}^{m}{X_n(s)}+\sum_{n=1}^{m}{Y_n(s)}+R(s)=0;\end{equation}
\par
Denote the sums of vectors $X_n$ and $Y_n$:
\begin{equation}\label{l_1_l_2}L_1=\sum_{n=1}^{m}{X_n(s)}; L_2=\sum_{n=1}^{m}{Y_n(s)};\end{equation}
\par
The vectors $L_1$ and $L_2$ are invariants of the vector system of the second approximate equation of the Riemann zeta function, since they do not depend on the order of the vectors $X_n$ and $Y_n$, nor on their quantity.
\par
We can now determine the geometric condition of the non-trivial zeros of the Riemann zeta function:
\begin{equation}\label{zeta_zero_vect}L_1+L_2+R=0;\end{equation}
\par
This condition means that when the Riemann zeta function takes a value of non-trivial zero when $\sigma=1/2$, the vectors $L_1$, $L_2$ and $R$ form \textit{an isosceles} triangle, because when $\sigma=1/2$ $|L_1|=|L_2|$, and when $\sigma\ne 1/2$, if the Riemann zeta function takes a value of non-trivial zero, these vectors must form a triangle of \textit{general form,} because when $\sigma\ne 1/2$ $|L_1|\ne |L_2|$.
\par
According to Hardy's theorem \cite{HA2}, the Riemann zeta function has an infinite number of zeros when $\sigma=1/2$.
\par
\textit{This fact is confirmed by the projection values $\hat\zeta(s)_L$ and $\hat\zeta(s)_M$ of the Riemann zeta function respectively on the normal $L$ to the axis of symmetry and on the axis of symmetry $M$ of the vector system of the second approximate equation of the Riemann zeta function on the critical line.}
\par
The Riemann zeta function takes a value of non-trivial zero when $\sigma=1/2$ every time when the projection $\hat\zeta(s)_L$ takes a value of zero, because the projection $\hat\zeta(s)_M$ when $\sigma=1/2$, according to the mirror symmetry of the vector system of the second approximate equation of the Riemann zeta function on the critical line are equel to zero identicaly.
\par
\textit{Geometric meaning of non-trivial zeros of the Riemann zeta function means that when the Riemann zeta function takes a value of non-trivial zero when $\sigma=1/2$, the vector system of the second approximate equation of the Riemann zeta function, as a consequence of the mirror symmetry of this vector system on the critical line, in accordance with Lemma 3, forms a symmetric polygon, since the sum of the vectors forming the polygon is equal to zero.}
\par
Using the vector system of the second approximate equation of the Riemann zeta function, we can also explain the geometric meaning of the Riemann-Siegel function \cite{SI,GA,TI}, which is used to compute the non-trivial zeros of the Riemann zeta function on the critical line:
\begin{equation}\label{zeta_ri_si}Z(t)=e^{\theta i}\zeta(\frac{1}{2}+it);\end{equation}
\par
where
\begin{equation}\label{teta}e^{\theta i}=\Big(\chi(\frac{1}{2}+it)\Big)^{\frac{1}{2}};\end{equation}
\par
According to the rules of multiplication of complex numbers, the Riemann-Siegel function determines the projection $\hat\zeta(s)_L$ of the Riemann zeta function on the normal $L$ to the axis of symmetry of the vector system of the second approximate equation of the Riemann zeta function, since
\begin{equation}\label{teta2}\theta=\frac{Arg(\chi(\frac{1}{2}+it))}{2}=\varphi_L;\end{equation}
\par
Which corresponds to the results of our research, i.e. $\hat \zeta(s)=\hat \zeta(s)_L$, since when $\sigma=1/2$ in accordance with the mirror symmetry of the vector system of the second approximate equation of the Riemann zeta function on the critical line $\hat \zeta(s)_M=0$.
\par
In conclusion of the research of the vector system of the second approximate equation of the Riemann zeta function, we establish another fundamental fact.
\par
The first middle vector $Y_1$ of the Riemann spiral rotates relative to the axes $x=Re(s)$ and $y=Im(s)$ around a fixed first vector $X_1$ of the Riemann spiral, and in accordance with the argument of the CHI function (\ref{chi_eq_app}) does $N(t)$ complete rotations:
\begin{equation}\label{x_1_num}N(t)=\frac{|Arg(\chi(s))|}{2\pi};\end{equation}
\par
Later we show that when $\sigma=1/2$ the first middle vector $Y_1$ of the Riemann spiral passes through the zero of the complex plane average once for each complete rotation
\footnote{The research of the vector system of the second approximate equation of the Riemann zeta function on the critical line shows that if the first middle vector $Y_1$ of the Riemann spiral for any complete rotation around the fixed first vector $X_1$ of the Riemann spiral never passes through the zero of the complex plane, then for another complete rotation it passes through the zero value of the complex plane twice.}
since, in accordance with the mirror symmetry of the vector system of the second approximate equation of the Riemann zeta function on the critical line, the end of first middle vector $Y_1$ makes reciprocating motions along the normal $L$ to the axis of symmetry of this vector system, since this normal passes through the end of the first vector $X_1$ of the Riemann spiral, which is at the zero of the complex plane.
\par
\textit{Thus, we can determine the number of non-trivial zeros of the Riemann zeta function on the critical line (as opposed to the Riemann-von Mangoldt formula, which determines the number of non-trivial zeros of the Riemann zeta function in the critical strip) via the number of complete rotations of the first middle vector $Y_1$ of the Riemann spiral:}
\begin{equation}\label{zeta_zero_num}N_0(t)=\Bigg[\frac{|Arg(\chi(s))-\alpha_2|}{2\pi}\Bigg]+2;\end{equation}
\par
where $\alpha_2$ argument to the second base point \footnote{base point is a value of a complex variable in which the first middle vector $Y_1$ occupies the position opposite to the first vector $X_1$ of the Riemann spiral.} of the first middle vector $Y_1$ of the Riemann spiral.
\section{Variants of confirmation of the Riemann hypothesis}
Before proceeding to the variants of confirmation of the Riemann hypothesis based on the analysis of the vector system of the second approximate equation of the Riemann zeta function, we consider two traditional approaches.
\par
The authors of the first approach estimate the proportion of non-trivial zeros $k$ and the proportion of simple zeros $k^*$ on the critical line compared to $N(T)$ the number of non-trivial zeros of the Riemann zeta function in the critical strip.
\begin{equation}\label{k}k=\lim_{T\to\infty}\inf\frac{N_0(T)}{N(T)};\end{equation}
\begin{equation}\label{k_s}k^*=\lim_{T\to\infty}\inf\frac{N_{0s}(T)}{N(T)};\end{equation}
\par
In the critical strip $N(T)$, the number of non-trivial zeros of the Riemann zeta function is determined by the Riemann-von Mangoldt formula \cite{TI}:
\begin{equation}\label{n_mn}N(T)=\frac{T}{2\pi}(\log{\frac{T}{2\pi}}-1)+\frac{7}{8}+S(T)+\mathcal{O}(\frac{1}{T});\end{equation}
\begin{equation}\label{s_mn}S(T)=\frac{1}{\pi}Arg(\zeta(\frac{1}{2}+iT))=\mathcal{O}(\log T), T\to\infty;\end{equation}
\par
This expression \cite{BU} is used to determine the proportion of non-trivial zeros $k$ and the proportion of simple zeros $k^*$ on the critical line:
\begin{equation}\label{k_2}k\ge 1-\frac{1}{R}\log\Big(\frac{1}{T}\int_{1}^{T}|V\psi(\sigma_0+it)|^2dt\Big)+o(1);\end{equation}
\begin{equation}\label{k_s_2}k^*\ge 1-\frac{1}{R}\log\Big(\frac{1}{T}\int_{2}^{T}|V\psi(\sigma_0+it)|^2dt\Big)+o(1);\end{equation}
\par
where $V(s)$ is some function whose number of zeros is the same as the number of zeros $\zeta(s)$ in a contour bounded by a rectangle (i. e. not on the critical line):
\begin{equation}\label{rect}\frac{1}{2}<\sigma<1; 0<t<T;\end{equation}
\par
$\psi(s)$ some \glqq mollifier\grqq \ function that has no zeros and compensates for the change of $|V(s)|$.
\par
R is some positive real number
\par
and
\begin{equation}\label{sigma_0}\sigma_0=\frac{1}{2}-\frac{R}{\log T};\end{equation}
\par
In recent papers \cite{BU, FE} the function $V(s)$ is used in the form of Levenson \cite{LE}.
\begin{equation}\label{v}V(s)=Q(-\frac{1}{\log T}\frac{d}{dt})\zeta(s);\end{equation}
\par
where $Q(x)$ is a real polynomial with $Q(0)=1$ and $Q'(x)=Q'(1-x)$.
\par
In this approach, the authors use different kinds of polynomials $Q(x)$, \glqq mollifier\grqq \ functions $\psi (s)$, as well as different methods of approximation and evaluation of the integral
\begin{equation}\label{int}\int_{1}^{T}|V\psi(\sigma_0+it)|^2dt;\end{equation}
\par
In the paper \cite{BU} in 2011 it is proved that
\begin{equation}\label{k_bu}k\ge .4105; k^*\ge .4058\end{equation}
\par
In the parallel paper \cite{FE} also in 2011 it is proved that
\begin{equation}\label{k_fe}k\ge .4128;\end{equation}
\par
To understand the dynamics of the results in this direction, we compare them with the paper \cite{CO} in which in 1989.
\begin{equation}\label{k_co}k\ge .4088; k^*\ge .4013\end{equation}
\par
\textit{The complexity of this approach lies in the fact that different methods of approximation of the integral (\ref{int}) allow to obtain an insufficiently accurate result.}
\par
We hope, the authors of \cite{PR} in 2019 found a way to accurately determine the number of non-trivial zeros in the critical strip and on the critical line, because they claim to have obtained the result $k=1$.
\par
The second traditional direction of confirmation of the Riemann hypothesis is connected with direct verification of non-trivial zeros of the Riemann zeta function.
\par
As we have already mentioned, the Riemann-Siegel formula (\ref{zeta_ri_si}) is used for this.
\par
One of the recent papers \cite{GO}, which besides computing $10^{13}$ of the first non-trivial zeros of the Riemann zeta function on the critical line offers a statistical analysis of these zeros and an improved approximation method of the Riemann-Siegel formula, was published in 2004.
\par
\textit{The complexity of this approach lies in the fact that it is impossible to calculate all the non-trivial zeros of the Riemann zeta function, therefore, by this method of confirmation the Riemann hypothesis it can be refuted rather than confirm.}
\par
It should be noted that there is \textit{a contradiction} between the results of the first and second method of confirmation of the Riemann hypothesis.
\par
The determination of the proportion of non-trivial zeros on the critical line is in no way related to \textit{the specific interval} of the imaginary part of a complex number, the authors of the first approach do not try to show that there is a sufficiently large interval where the Riemann hypothesis is true in accordance with the second method of confirmation the Riemann Hypothesis.
\par
In other words, by determining the proportion of non-trivial zeros of the Riemann zeta function on the critical line, the authors of the first direction of confirmation of the Riemann hypothesis \textit{indirectly indicate that part of the non-trivial zeros of the Riemann zeta function lies not on the critical line at any interval}, i.e. even where these zeros are already verified by the second method.
\par
The site \cite{WA} collects some unsuccessful attempts to prove the Riemann hypothesis, some of them are given with a detailed error analysis.
\par
In his presentation, Peter Sarnak \cite{SA}, in the course of analyzing the different approaches to proving the Riemann hypothesis, mentions that about three papers a week are submitted for consideration annually.
\par
However, the situation with the proof of the Riemann hypothesis remains uncertain and we tend to agree with Pete Clark \cite{KL}:
\par
\textit{So far as I know, there is no approach to the Riemann Hypothesis which has been fleshed out far enough to get an even moderately skeptical expert to back it, with any odds whatsoever. I think this situation should be contrasted with that of Fermat's Last Theorem [FLT]: a lot of number theorists, had they known in say 1990 that Wiles was working on FLT via Taniyama-Shimura, would have found that plausible and encouraging.}
\par
Based on the research results described in the second section of our paper, we propose several new approaches to confirm the Riemann hypothesis based on the properties of the vector system of the second approximate equation of the Riemann zeta function:
\par
1) the first method is based on determining the exact number (\ref{zeta_zero_num}) of non-trivial zeros of the Riemann zeta function on the critical line;
\par
2) the second method of confirmation of the Riemann hypothesis is based on the analysis of projections $\zeta(s)_L$ and $\zeta(s)_M$ of the Riemann zeta function respectively on the normal $L$ to the axis of symmetry and on the axis of symmetry $M$ of the vector system of the second approximate equation of the Riemann zeta function in the critical strip;
\par
3) the third method is based on the vector condition (\ref{zeta_zero_vect}) of non-trivial zeros of the Riemann zeta function.
\par
Determining the exact number (\ref{zeta_zero_num}) of non-trivial zeros of the Riemann zeta function on the critical line is based on several facts:
\par
1) the vector system of the second approximate equation of the Riemann zeta function when $\sigma=1/2$ has mirror symmetry (\ref{zeta_m});
\par
2) the vector system of the second approximate equation of the Riemann zeta function rotates around the end of the first vector $X_1$ of the Riemann spiral (\ref{phi_m}) in a fixed coordinate system of the complex plane;
\par
3) vectors $X_n$ and the middle vectors $Y_n$ of the Riemann spirals rotate in opposite directions (\ref{x_vect}) and (\ref{y_vect}) in the moving coordinate system formed by a normal $L$ to the axis of symmetry and axis of symmetry $M$ of the vector system of the second approximate equation of the Riemann zeta function.
\par
The first middle vector $Y_1$ of the Riemann spiral rotates with the vector system of the second approximate equation of the Riemann zeta function around the end of the first vector $X_1$ of the Riemann spiral (\ref{phi_m}) in the fixed coordinate system of the complex plane, so the first middle vector $Y_1$ of the Riemann spiral, in accordance with (\ref{chi_eq_app}), periodically passes \textit{base point} where it takes a position opposite to the first vector $X_1$ of the Riemann spiral:
\begin{equation}\label{base_point}Arg(\chi(\frac{1}{2}+it_k))=(2k-1)\pi;\end{equation}
\par
The argument of \textit{the base points} of the first middle vector $Y_1$ of the Riemann spiral differs from the argument of the Gram points \cite{GR} by $\pi/2$:
\begin{equation}\label{gram_point}\theta=\frac{Arg(\chi(\frac{1}{2}+it_n))}{2}=(n-1)\pi;\end{equation}
\par
\textit{We consider that the complete rotation of the first middle vector $Y_1$ of the Riemann spiral is the rotation from any base point to the next base point.}
\begin{figure}[ht]
\centering
\includegraphics[scale=0.6]{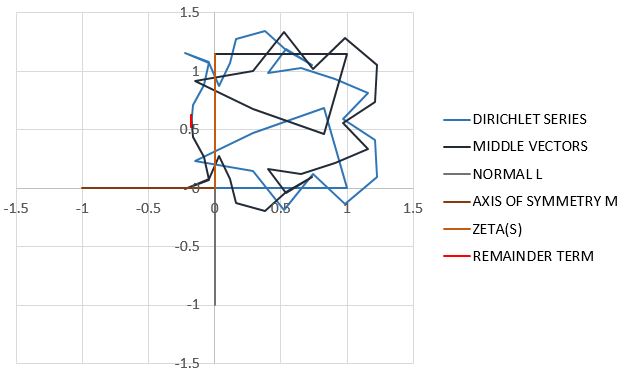}
\caption{Base point \#4520, the first middle vector is above the real axis, $s=0.5+5001.099505i$}
\label{fig:gram_point_4520}
\end{figure}
\par
\begin{figure}[ht]
\centering
\includegraphics[scale=0.6]{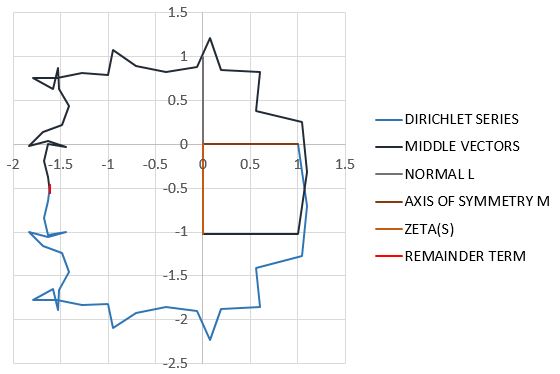}
\caption{Base point \#4525, the first middle vector is below the real axis, $s=0.5+5005.8024855i$}
\label{fig:gram_point_4525}
\end{figure}

In accordance with the mirror symmetry of the vector system of the second approximate equation of the Riemann zeta function when $\sigma=1/2$, at the base points the end of the first middle vector $Y_1$ of the Riemann spiral is on the imaginary axis of the complex plane, so the vector can occupy one of two positions above (fig. \ref{fig:gram_point_4520}) or bottom (fig. \ref{fig:gram_point_4525}) of the real axis of the complex plane.
\par
Since the vectors $X_n$ and the middle vectors $Y_n$ of the Riemann spiral rotate in opposite directions (\ref{x_vect}, \ref{y_vect}) in the moving coordinate system formed by the normal $L$ to the axis of symmetry and the axis of symmetry $M$ of the vector system of the second approximate equation of the Riemann zeta function, in the fixed coordinate system of the complex plane, the first middle vector $Y_1$ of the Riemann spiral, when it is below the real axis of the complex plane, rotates towards the first vector $X_1$ Riemann Spirals and opposite when the first middle vector $Y_1$ of the Riemann spiral is on top of the real axis of the complex plane, it rotates away from the first vector $X_1$ of the Riemann spiral.
\par
\begin{figure}[ht]
\centering
\includegraphics[scale=0.6]{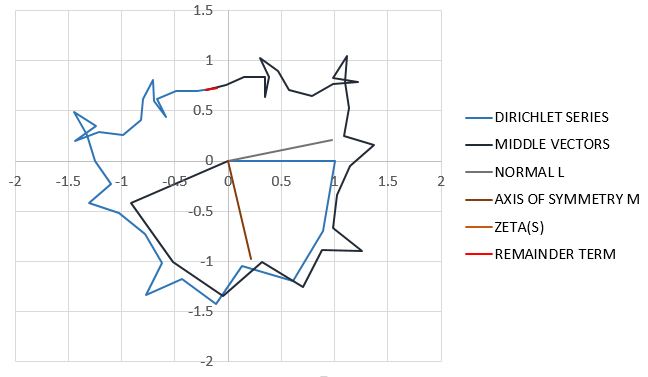}
\caption{Non-trivial zero of the Riemann zeta function \#4525, after the base point, $s=0.5+5006.208381106i$}
\label{fig:root_4525}
\end{figure}

If the first middle vector $Y_1$ of the Riemann spiral at the base point rotates towards the first vector $X_1$ of the Riemann spiral, then, in accordance with the mirror symmetry of the vector system of the second approximate equation of the Riemann zeta function when $\sigma=1/2$, until the rotation of the first middle vector $Y_1$ of the Riemann spiral is completed, the axis of symmetry $M$ of this vector system crosses the zero of the complex plane (fig. \ref{fig:root_4525}).
\par
\begin{figure}[ht]
\centering
\includegraphics[scale=0.6]{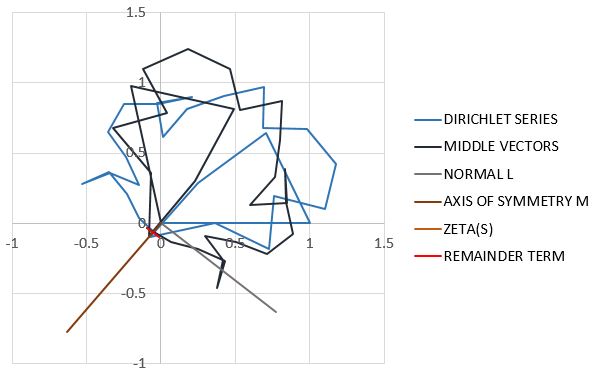}
\caption{Non-trivial zero of the Riemann zeta function \#4520, up to the base point, $s=0.5+5000.834381i$}
\label{fig:root_4520}
\end{figure}

Conversely, if the first middle vector $Y_1$ of the Riemann spiral at the base point rotates away from the first vector $X_1$ of the Riemann spiral, then from the beginning of the rotation of the first middle vector $Y_1$ of the Riemann spiral, the axis of symmetry $M$ of this vector system has already crossed the zero of the complex plane (fig. \ref{fig:root_4520}).
\par
In the result of the mirror symmetry of the vector system of the second approximate equation of the Riemann zeta function, when the axis of symmetry $M$ of this vector system crosses the zero of the complex plane, the end of the first middle vector $Y_1$ of the Riemann spiral also crosses the zero of the complex plane, so the vector system of the second approximate equation of the Riemann zeta function at this point forms a closed polyline.
\par
\textit{It is known from analytical geometry that the sum of vectors forming a closed polyline is equal to zero.}
\par
Therefore, when the first middle vector $Y_1$ of the Riemann spiral is at the base point, we can be sure that either before that point or after that point the Riemann zeta function takes a value of non-trivial zero.
\par
\textit{This fact allows us to conclude that one non-trivial zero of the Riemann zeta function corresponds to one complete rotation of the first middle $Y_1$ vector of the Riemann spiral.}
\par
In accordance with of different combinations of the location of the first middle vector $Y_1$ of the Riemann spiral relative to the real axis of the complex plane at adjacent base points, there may be a different number of non-trivial zeros of the Riemann zeta function in the intervals between the base points:
\par
a) one non-trivial zero if the first middle vector $Y_1$ of the Riemann spiral occupies the same positions at two serial base points;
\par
b) two zeros if at the first base point the first middle vector $Y_1$ of the Riemann spiral is below and at the next base point is above the real axis of the complex plane;
\par
c) no zero, if, on the contrary, at the first base point the first middle vector $Y_1$ of the Riemann spiral is above, and at the next base point, is below the real axis of the complex plane.
\par
Denote the types of base points:
\par
$a_1$ - if the first middle vector $Y_1$ of the Riemann spiral occupies a position above\footnote{If the Riemann zeta function can take a value of non-trivial zero at the base point, then the first middle vector $Y_1$ will occupy the position of the first vector $X_1$ of the Riemann spiral (we will not consider this state in detail), we assume that this position belongs to the base point type $a_1$.} the real axis of the complex plane;
\par
$a_2$ - if the first middle vector $Y_1$ of the Riemann spiral occupies the position below the real axis of the complex plane.
\par
Then we can define a sequence of base points of the same types, which correspond to the first type of interval, which has one non-trivial zero: $$A_1=a_1a_1;$$ $$A_2=a_2a_2;$$
\par
and sequence base points with different types, which correspond to the second and third interval type, respectively: $$B=a_2a_1;$$ $$C=a_1a_2;$$
\par
It is obvious that the sequence $A_1$ and $A_2$ can't follow each other, because $$a_1a_1a_2a_2=A_1 CA_2$$
\par
or $$a_2 a_2a_1 a_1=A_2 BA_1$$
\par
It is also clear that each other can not follow the sequence $B$, since $$a_2a_1 a_2 a_1=BC$$
\par
and each other can not follow sequence $C$, since $$a_1a_2a_1a_2=CBC$$
\par
Therefore, if at some interval between two base points of the Zeta-function of Riemann has no non-trivial zeros, i.e. is the interval of type $C$, then at another interval the Zeta-function of Riemann will have two non-trivial zero, it is the interval of type $B$, since these intervals appear every time when the type of base point changed, such as $$a_1a_2a_1=CB$$ or more long chain $$a_1a_2a_2a_2a_1a_1a_1a_1a_1a_2a_2a_1=CA_2A_2BA_1A_1A_1A_1CA_2B$$
\par
\textit{Thus, the total number of non-trivial zeros of the Riemann zeta function at critical line always corresponds to the number of base points or the number of complete rotations of the first middle $Y_1$ vector of the Riemann spiral around the end of the first vector $X_1$ of the Riemann spiral in the fixed coordinate system of the complex plane.}
\par
We used this property to obtain the expression (\ref{zeta_zero_num}) the number of non-trivial zeros of the Riemann zeta function on the critical line through the angle of the first middle vector $Y_1$ of the Riemann spiral.
\par
We substitute the exact expression (\ref{arg_chi}) argument CHI functions in (\ref{zeta_zero_num}):
\begin{equation}\label{zeta_zero_num_0}N_0(T)=\Bigg[\Big|\frac{T}{2\pi}(\log{\frac{T}{2\pi}}-1)-\frac{1}{8}+\frac{2\mu(T)-\alpha_2}{2\pi}\Big|\Bigg]+2;\end{equation}
\par
where $\mu(T)$ is the remainder term of the gamma function (\ref{mu}) when $\sigma=1/2$;
\par
$\alpha_2$ argument of the CHI function at the second base point.
\par
Comparing the expression for the number of non-trivial zeros in the critical strip (\ref{n_mn}) and the expression for the number of non-trivial zeros on the critical line (\ref{zeta_zero_num_0}), we see that these values \glqq match\grqq.
\par
We are in a paradoxical situation where we know the exact number of zeros on the critical line, because $\mu(T)\to 0$ at $T\to \infty$ and do not know the exact number of zeros in the critical strip, because in the most optimistic estimate \cite{TR}:
\begin{equation}\label{s_tr}|S(T)|<1.998+0.17\log(T); T>e;\end{equation}
\par
\textit{Thus, based on the last estimate (\ref{s_tr}) of the remainder term of the Riemann-von Mangoldt formula, we can say that \glqq almost all\grqq \ non-trivial zeros of the Riemann zeta function lie on the critical line.}
\par
It should be noted that for the final solution of the problem by this method it is not enough to show that
\begin{equation}\label{s_li}|S(T)|<\mathcal{O}(\frac{\log(T)}{\log\log(T)}); T\to \infty;\end{equation}
\par
Although that this result was obtained by Littlewood \cite{LI} provided that the Riemann hypothesis is true, since the expression $(2\mu(T)-\alpha_2)/2\pi$ has a limit at $T\to \infty$, and the expression $\log(T)/\log\log (T)$ has no such limit, although it grows very slowly.
\par
\textit{In other words, by comparing the number of non-trivial zeros of the Riemann zeta function in the critical strip and on the critical line, it is almost impossible to confirm the Riemann hypothesis.}
\par
Using the method of analyzing the vector system of the second approximate equation of the Riemann zeta function, we can offer a confirmation of the Riemann hypothesis \textit{from the contrary.}
\par
This approach is to show that Riemann zeta functions \textit{cannot have non-trivial zeros} when $\sigma\ne 1/2$.
\par
This problem is solved in different ways by the second and third methods of confirmation of the Riemann hypothesis, based on the properties of the vector system of the second approximate equation of the Riemann zeta function.
\par
We have already considered the dynamics of the first middle vector $Y_1$ of the Riemann spiral at the base points when $\sigma=1/2$, now consider the dynamics of the first middle vector $Y_1$ of the Riemann spiral at the base points when $\sigma\ne 1/2$.
\par
In accordance with (\ref{alpha1_ex}) base point when $\sigma\ne 1/2$ are in the neighborhood $\epsilon=\mathcal{O}(t^{-1})$ points when $\sigma=1/2$, i.e. they have practically the same value at $t\to \infty$.
\par
Therefore, at the base points when $\sigma\ne 1/2$ the first middle vector $Y_1$ of the Riemann spiral occupies a position opposite to the first vector $X_1$ of the Riemann spiral, while, in accordance with the violation of the mirror symmetry of the vector system of the second approximate equation of the Riemann zeta function, the end of the first middle vector $Y_1$ of the Riemann spiral \textit{cannot be} on the imaginary axis of the complex plane.
\par
\begin{figure}[ht]
\centering
\includegraphics[scale=0.6]{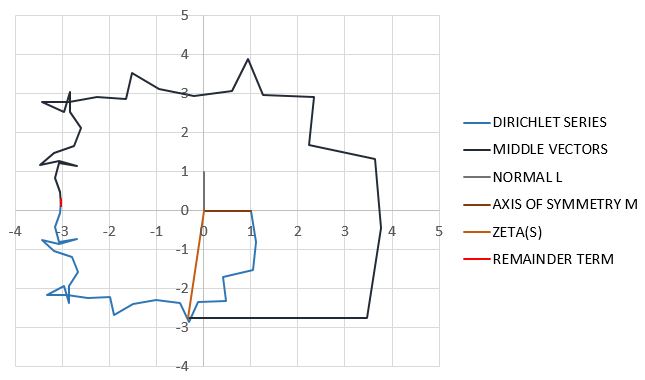}
\caption{Base point \#4525, the first middle vector of the Riemann spiral is on the left, $s=0.35+5005.8024855i$}
\label{fig:gram_point_4525_left}
\end{figure}
\begin{figure}[ht]
\centering
\includegraphics[scale=0.6]{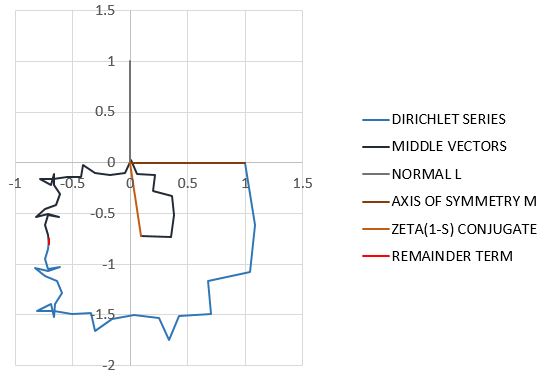}
\caption{Base point \#4525, the first middle vector of the Riemann spiral is on the right, $s=0.65+5005.8024855i$}
\label{fig:gram_point_4525_right}
\end{figure}

Therefore, when $\sigma<1/2$, the end of the first middle vector $Y_1$ of the Riemann spiral is on the left (fig. \ref{fig:gram_point_4525_left}), when $\sigma>1/2$ it is on the right (fig. \ref{fig:gram_point_4525_right}) from the imaginary axis of the complex plane, regardless of the position relative to the real axis of the complex plane.
\par
\begin{figure}[ht]
\centering
\includegraphics[scale=0.6]{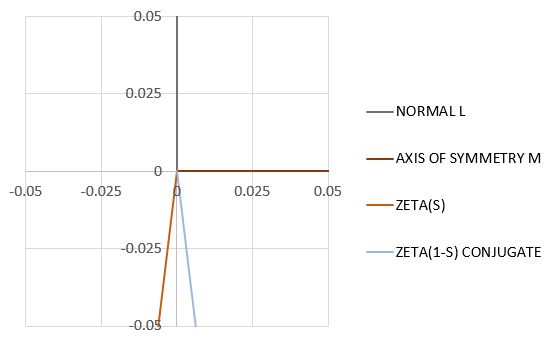}
\caption{The deviation of the vector of values of the Riemann zeta function in the base point \#4525, when $\sigma=0.35$ and $\sigma=0.65$}
\label{fig:gram_point_4525_two_angles}
\end{figure}

It is easy to notice that at the base point of the vector values of the Riemann zeta function at values $\sigma+it$ and $1 - \sigma+it$ deviate from the normal $L$ to the axis of symmetry of the vector system of the second approximate equation in different directions at \textit{the same angle} (fig. \ref{fig:gram_point_4525_two_angles}).
\par
This behavior of vectors of values of the Riemann zeta function for values of the complex variable symmetric about the critical line can be easily explained by the arithmetic of arguments of a complex numbers at the base point.
\begin{equation}\label{arg1}Arg(\zeta(s))_B=Arg(\chi(s))_B+Arg(\zeta(1-s))_B\end{equation}
\par
At the base point
\begin{equation}\label{arg2}Arg(\chi(s))_B=\pi;\end{equation}
\par
since
\begin{equation}\label{zeta_conj3}\zeta(1-\sigma+it)=\zeta(\overline{1-\sigma-it})=\zeta(\overline{1-s})=\overline{\zeta(1-s)};\end{equation}
\begin{equation}\label{arg3}Arg(\zeta(1-s))=-Arg(\overline{\zeta(1-s)})=-Arg(\zeta(1-\sigma+it));\end{equation}
\par
then at the base point
\begin{equation}\label{arg4}Arg(\zeta(\sigma+it))_B=\pi-Arg(\zeta(1-\sigma+it))_B;\end{equation}
\par
This ratio of arguments is kept in the moving coordinate system formed by the normal $L$ to the axis of symmetry and the axis of symmetry $M$ of the vector system of the second approximate equation of the Riemann zeta function \textit{for any values} of the complex variable symmetric about the critical line.
\begin{equation}\label{arg5}Arg(\zeta(\sigma+it))-\frac{Arg(\chi(\sigma+it))}{2}=-(Arg(\zeta(1-\sigma+it))-\frac{Arg(\chi(\sigma+it))}{2});\end{equation}
\par
The vector system of the second approximate equation of the Riemann zeta function rotates relative to the end of the first vector $X_1$ of the Riemann spiral (\ref{phi_m}), while, in accordance with (\ref{x_vect}, \ref{y_vect}), each next vector of this vector system rotates relative to the previous vector in the same direction as the entire vector system, therefore, the angle of twist of the vectors \textit{grows monotonically.}
\par
Thus, in accordance with (\ref{x_vect}, \ref{y_vect}) when $\sigma<1/2$, the angle of twist of vectors \textit{grows faster} than when $\sigma=1/2$, because the angles between the vectors are equal, but when $\sigma<1/2$, the modulus of each vector is greater than the modulus of the corresponding vector when $\sigma=1/2$, so the vector system when $\sigma<1/2$ is twisted \textit{at a greater angle} than in the case $\sigma=1/2$.
\par
While, when $\sigma>1/2$, the angle of twist of vectors \textit{grows slower} than when $\sigma=1/2$, because the angles between the vectors are equal to, but $\sigma>1/2$, the modulus of each vector is less than the modulus of the corresponding vector when $\sigma=1/2$, so the vector system when $\sigma>1/2$ is twisted \textit{at a smaller angle} than in the case $\sigma=1/2$.
\begin{figure}[ht]
\centering
\includegraphics[scale=0.6]{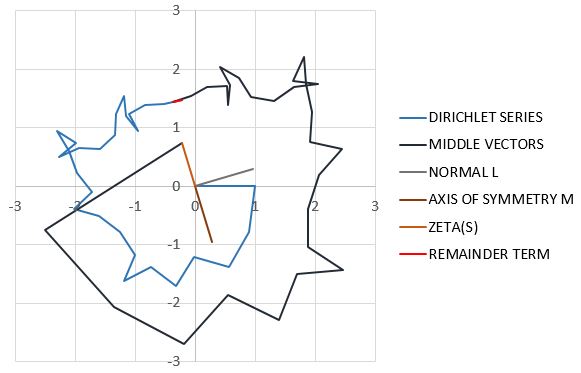}
\caption{Base point \#4525, the vector of values of the Riemann zeta function parallel to the axis of symmetry, $s=0.35+5006.186i$}
\label{fig:gram_point_4525_symmetry1}
\end{figure}
\par
\begin{figure}[ht]
\centering
\includegraphics[scale=0.6]{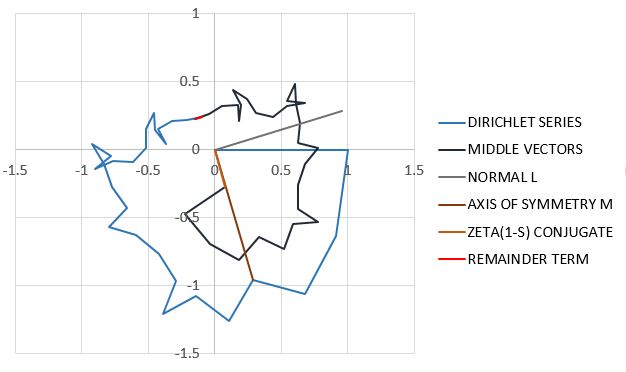}
\caption{Base point \#4525, the vector of values of the Riemann zeta function parallel to the axis of symmetry, $s=0.65+5006.186i$}
\label{fig:gram_point_4525_symmetry2}
\end{figure}
\begin{figure}[ht]
\centering
\includegraphics[scale=0.6]{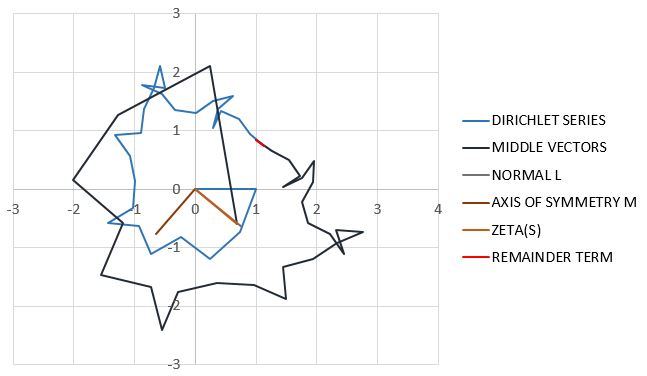}
\caption{Base point \#4525, the vector of values of the Riemann zeta function parallel to the normal to the axis of symmetry, $s=0.35+5006.484i$}
\label{fig:gram_point_4525_normal1}
\end{figure}
\begin{figure}[ht!]
\centering
\includegraphics[scale=0.6]{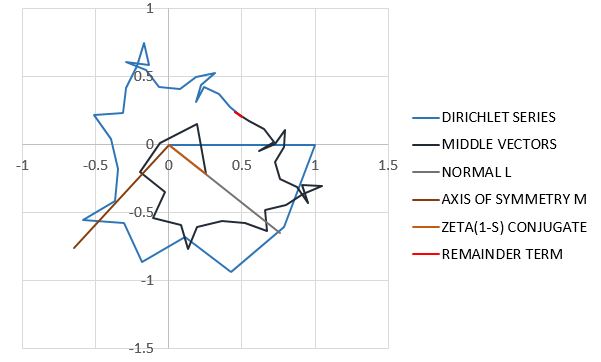}
\caption{Base point \#4525, the vector of values of the Riemann zeta function parallel to the normal to the axis of symmetry, $s=0.65+5006.484i$}
\label{fig:gram_point_4525_normal2}
\end{figure}

In accordance with the identified ratio of the arguments (\ref{arg5}) and a monotonic increase of the angle of twist of vectors, we can conclude that the vectors of value of the Riemann zeta function, when values of the complex variable are symmetric about the critical line, rotate in the moving coordinate system formed by a normal of $L$ to the axis of symmetry and axis of symmetry $M$ of the vector system of the second approximate equation of the Riemann zeta function in opposite directions with the same speed and therefore at different values of the imaginary part of a complex number has a special positions:
\par
a) directed in different sites, along the axis of symmetry $M$ (fig. \ref{fig:gram_point_4525_symmetry1}, \ref{fig:gram_point_4525_symmetry2});
\par
b) directed in the same site, along the normal $L$ to the axis of symmetry of the vector system of the second approximate equation of the Riemann zeta function (fig. \ref{fig:gram_point_4525_normal1}, \ref{fig:gram_point_4525_normal2}).
\par
While, when $\sigma=1/2$, in accordance with the mirror symmetry of the vector system of the second approximate equation of the Riemann zeta function when $\sigma=1/2$, the vector of value of the Riemann zeta function rotates so that it is always located along the normal $L$ to the axis of symmetry of this coordinate system (fig. \ref{fig:gram_point_4525}). 
\par
According to the position relative to the real axis of the complex plane at the base point when $\sigma\ne 1/2$, the first middle vector $Y_1$ of the Riemann spiral, when it is below the real axis of the complex plane, rotates towards the first vector $X_1$ of the Riemann spiral and, conversely, when the first middle vector $Y_1$ of the Riemann spiral is above the real axis of the complex plane, it rotates away from the first vector $X_1$ of the Riemann spiral.
\par
If the first middle vector $Y_1$ of the Riemann spiral at the base point when $\sigma\ne 1/2$ rotates towards the first vector $X_1$ of the Riemann spiral, then, in accordance with the conformal symmetry of the vector system of the second approximate equation of the Riemann zeta function when $\sigma\ne 1/2$, until the rotation of the first middle vector $Y_1$ of the Riemann spiral is completed, it will take a special position when the axis of symmetry $M$ of this vector system passes through both the zero of the complex plane and the end of the first middle vector $Y_1$ of the Riemann spiral.
\par
Conversely, if the first middle vector $Y_1$ of the Riemann spiral at the base point when $\sigma\ne 1/2$ rotates away from the first vector $X_1$ of the Riemann spiral, from the beginning of the rotation of the first middle vector $Y_1$ of the Riemann spiral, it already occupied a special position when the axis of symmetry $M$ of this vector system is passed through both the zero of the complex plane and the end of the first middle vector $Y_1$ of the Riemann spiral.
\par
This special position of the first middle vector $Y_1$ of the Riemann spiral corresponds to the special position of the vector of value of the Riemann zeta function when it locates along to the axis of symmetry $M$ (fig. \ref{fig:gram_point_4525_symmetry1}, \ref{fig:gram_point_4525_symmetry2}) of the vector system of the second approximate equation of the Riemann zeta function.
\par
According to the rules of summation of vectors in this special position of the first middle vector $Y_1$ of the Riemann spiral, the projection of the vector system of the second approximate equation of the Riemann zeta function on the normal $L$ to the axis of symmetry of this vector system is equal to zero, hence
\begin{equation}\label{zeta_l_y1}\zeta(s)_L=0;\end{equation}
\par
By increasing the imaginary part of a complex number, the first middle vector $Y_1$ of the Riemann spiral will move to another special position when the normal $L$ to the axis of symmetry of the vector system of the second approximate equation of the Riemann zeta function passes through both the zero of the complex plane and the end of the first middle vector $Y_1$ of the Riemann spiral at that moment it rotated by an angle $\pi/2$ from the first special position in the moving coordinate system formed by the normal $L$ to the axis of symmetry and the axis of symmetry $M$ of this vector system.
\par
Second special position of the first middle vector $Y_1$ of the Riemann spiral corresponds to the special position of the vector of value of the Riemann zeta function when it locates along the normal $L$ to the axis of symmetry (fig. \ref{fig:gram_point_4525_normal1}, \ref{fig:gram_point_4525_normal2}) of the vector system of the second approximate equation of the Riemann zeta function.
\par
According to the rules of summation of vectors in this special position of the first middle vector $Y_1$ of the Riemann spiral, the projection of the vector system of the second approximate equation of the Riemann zeta function on the axis of symmetry $M$ of this vector system is equal to zero, hence
\begin{equation}\label{zeta_m_y1}\zeta(s)_M=0;\end{equation}
\par
Thus, when we performed an additional analysis of the vector system of the second approximate equation of the Riemann zeta function in accordance with the identified ratio of the arguments (\ref{arg5}) and monotonic increase of the angle of twist of vectors we found that each base point corresponds to \textit{two special positions} of the first middle vector $Y_1$ of the Riemann spiral, in which when $\sigma\ne 1/2$ the projections of the vector system of the second approximate equation of the Riemann zeta function on the normal $L$ to the axis of symmetry and on the axis of symmetry $M$ of this vector system \textit{alternately} take a value of zero and when $\sigma=1/2$ one of which corresponds to the non-trivial zero of the Riemann zeta function.
\par
Now we need to make sure that the modulus of the Riemann zeta function cannot take a value of zero except for the special position of the first middle vector $Y_1$ of the Riemann spiral when $\sigma=1/2$.
\par
We will analyzing the projections of the vector system of the second approximate equation of the Riemann zeta function on the normal $L$ to the axis of symmetry and on the axis of symmetry $M$ of this vector system starting from the boundary of the critical strip, where the Riemann zeta function has no non-trivial zeros, i.e. for values $\sigma=0$ and $\sigma=1$.
\par
\begin{figure}[ht]
\centering
\includegraphics[scale=0.6]{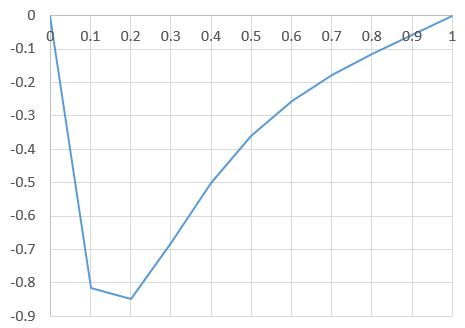}
\caption{Base point \#4525, projection of the vector system on the normal $L$ to the axis of symmetry when $\zeta(1+5006,09072i)_L=0$}
\label{fig:gram_point_4525_zeta_1_l}
\end{figure}

Construct the graphs of the projections of the vector system of the second approximate equation of the Riemann zeta function on the normal $L$ to the axis of symmetry and on the axis of symmetry $M$ of this vector system when $\zeta(1+it)_L=0$, in accordance with the ratio of the arguments (\ref{arg5}) of the vectors of values of the Riemann zeta function, when values of the complex variable symmetric about the critical line, we obtain another equality $\zeta(0+it)_L=0$ (fig. \ref{fig:gram_point_4525_zeta_1_l}).
\begin{figure}[ht]
\centering
\includegraphics[scale=0.6]{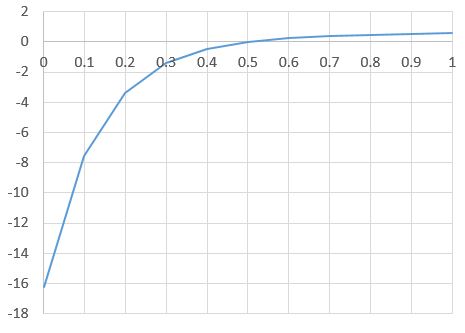}
\caption{Base point \#4525, projection of the vector system on the axis of symmetry $M$ when $\zeta(1+5006,09072i)_L=0$}
\label{fig:gram_point_4525_zeta_1_m}
\end{figure}
\par
Analysis of the projections of the vector system of the second approximate equation of the Riemann zeta function on the normal $L$ to the axis of symmetry and on the axis of symmetry $M$ of this vector system shows that in the interval $A_k$ of the imaginary part of a complex number from $t_1$: $\zeta(1+it_1)_L=0$ to $t_2$: $\zeta(1/2+it_2)_L=0$, any value of $t$ corresponds to a value $0<\sigma<1$: $\zeta(\sigma+it)_L=0$.
\par
In other words, in the interval $A_k$ for each value of the imaginary part of a complex number, the graph of the function $\zeta(\sigma+it)_L=0$ crosses the abscissa axis twice for the symmetric values $\sigma+it$ and $1-\sigma+it$ until when $\sigma=1/2$ it reaches a value of the non-trivial zero of the Riemann zeta function.
\par
While the graph of the function $\zeta(\sigma+it) _M=0$ at any value $t$ of the imaginary part of a complex number from the interval $A_k$ crosses the abscissa axis only once, at the point $\sigma=1/2$ (fig. \ref{fig:gram_point_4525_zeta_1_m}).
\par
\begin{figure}[ht]
\centering
\includegraphics[scale=0.6]{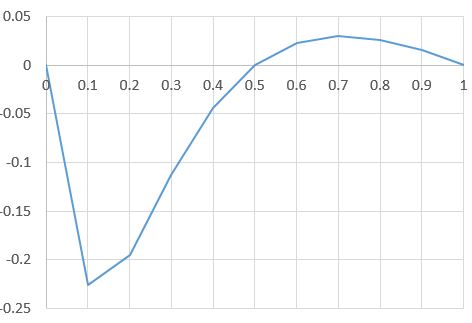}
\caption{Base point \#4525, projection of the vector system on the axis of symmetry $M$ when $\zeta(1+5006,4559i)_M=0$}
\label{fig:gram_point_4525_zeta_2_m}
\end{figure}

Now we construct the graphs of the projections of the vector system of the second approximate equation of the Riemann zeta function on the normal $L$ to the axis of symmetry and on the axis of symmetry $M$ of this vector system when $\zeta(1+it)_M=0$, in accordance with the ratio of the arguments (\ref{arg5}) of the vectors of value of the Riemann zeta function, when values of the complex variable symmetric about the critical line, we obtain another equality $\zeta(0+it)_M=0$ (fig. \ref{fig:gram_point_4525_zeta_2_m}).
\begin{figure}[ht]
\centering
\includegraphics[scale=0.6]{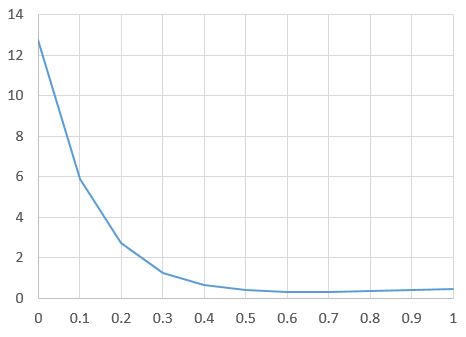}
\caption{Base point \#4525, projection of the vector system on the normal $L$ to the axis of symmetry when $\zeta(1+5006,4559i)_M=0$}
\label{fig:gram_point_4525_zeta_2_l}
\end{figure}
\par
Analysis of the projections of the vector system of the second approximate equation of the Riemann zeta function on the normal $L$ to the axis of symmetry and on the axis of symmetry $M$ of this vector system shows that in the interval $C_k$ of the imaginary part of a complex number from $t'_1$: $\zeta(1+it'_1)_M=0$ to $t'_2$: $\zeta(1/2+it'_2)_M=0$, any value of $t'$ corresponds to a value $0<\sigma<1$: $\zeta(\sigma+it')_M=0$.
\par
In other words, in the interval $C_k$ for each value of the imaginary part of a complex number, the graph of the function $\zeta(\sigma+it')_M=0$ crosses the abscissa axis three times for the symmetric values $\sigma+it'$, $1-\sigma+it'$ and at the point $\sigma=1/2$.
\par
While the graph of the function $\zeta(\sigma+it') _L=0$ at any value $t'$ of the imaginary part of a complex number from the interval $C_k$ never crosses the abscissa axis (fig. \ref{fig:gram_point_4525_zeta_2_l}).
\par
Therefore, in the interval $C_k$ for any value of the imaginary part of a complex number and any value of the real part of a complex number $0<\sigma<1$, the Riemann zeta function has no non-trivial zeros.
\begin{figure}[ht]
\centering
\includegraphics[scale=0.6]{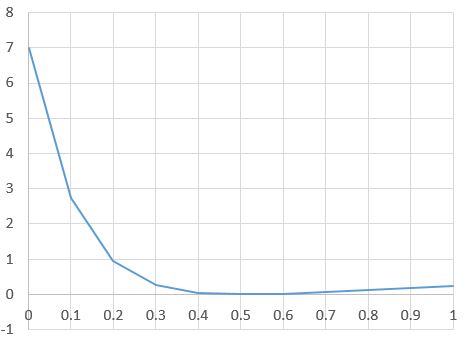}
\caption{Base point \#4525, projection of the vector system on the normal $L$ to the axis of symmetry when $\zeta(1/2+5006,208381106i)_L=0$}
\label{fig:gram_point_4525_zeta_3_l}
\end{figure}
\begin{figure}[ht!]
\centering
\includegraphics[scale=0.6]{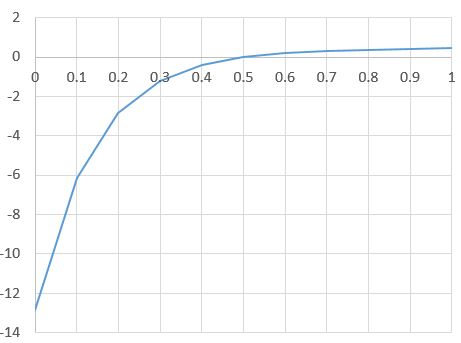}
\caption{Base point \#4525, projection of the vector system on the axis of symmetry $M$ when $\zeta(1/2+5006,208381106i)_L=0$}
\label{fig:gram_point_4525_zeta_3_m}
\end{figure}
\par
Analysis of the projections of the vector system of the second approximate equation of the Riemann zeta function on the normal $L$ to the axis of symmetry and on the axis of symmetry $M$ of this vector system in the interval of $B_k$ between intervals $A_k$ and $C_k$ and in the interval $D_k$ between intervals $C_k$ and $A_{k+1}$ shows that the projection at any value of the imaginary part of a complex number and any value of the real part of a complex number when $0<\sigma<1$ is not equal to zero, therefore, in the interval of $B_k$ and $D_k$ the Riemann zeta function has no non-trivial zeros.
\par
\begin{figure}[ht]
\centering
\includegraphics[scale=0.6]{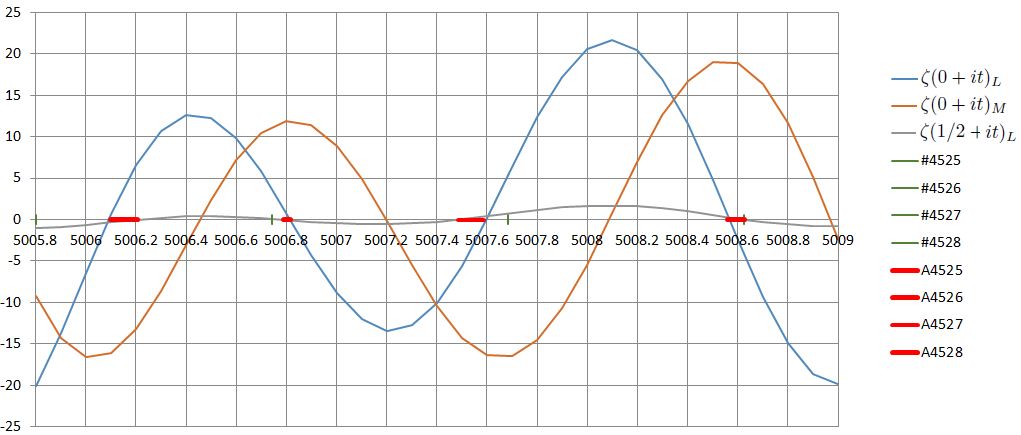}
\caption{Graphics of projections of the Riemann zeta function}
\label{fig:graphics_projections}
\end{figure}

It should be noted that the sign of the projection of the vector system on the axis of symmetry $M$ of the vector system of the second approximate equation of the Riemann zeta function changes at each interval $C_k$, i.e. it depends on the number of the base point (fig. \ref{fig:graphics_projections}).
\par
While the sign of the projection of the vector system on the normal $L$ to the axis of symmetry of the vector system of the second approximate equation of the Riemann zeta function depends on the location of the first middle vector $Y_1$ of the Riemann spiral at the base point (fig. \ref{fig:graphics_projections}), i.e. it changes at each interval $A_k$.
\par
Thus, when we performed an additional analysis of the vector system of the second approximate equation of the Riemann zeta function we found that each base point corresponds to \textit{four intervals} in which the projections of the vector system of the second approximate equation of the Riemann zeta function on the normal $L$ to the axis of symmetry and on the axis of symmetry $M$ of this vector system take certain values and only in one interval when $\sigma=1/2$ they can be zero at the same time (fig. \ref{fig:gram_point_4525_zeta_3_l} and \ref{fig:gram_point_4525_zeta_3_m}), in this moment the Riemann zeta function takes a value of non-trivial zero.
\par
Now that we know all possible variants of the ratio of the projection of vectors system of the second approximate equation of the Riemann zeta function on normal $L$ to the axis of symmetry and on the axis of symmetry $M$ of this vector system we found out \textit{what and why} can be these projections for any value of the real part of a complex number in the critical strip, and on the boundary of the critical strip, i.e. for values  $\sigma=0$ and $\sigma=1$, where the Riemann zeta function has no non-trivial zeros, as it proved Adamar and vallée Poussin.
\par
We have found out how these projections change when the imaginary part of the complex number changes and how these changes are related to the number of the base point and the position of the first middle vector $Y_1$ of the Riemann spiral at the base point, we have to answer two questions:
\par
1) Why should these projections have such ratios for any base point?
\par
2) Why can't there be any other reason for the modulus of the Riemann zeta function to go to zero when $\sigma\ne 1/2$?
\par
The first question has already been answered in the analysis of the projections of the vector system on the normal $L$ to the axis of symmetry and on the axis of symmetry $M$ of the vector system of the second approximate equation of the Riemann zeta function:
\par
1) The vector system of the second approximate equation of the Riemann zeta function rotates when the imaginary part of the complex number changes, which we can confirm this by the equation of the axis of symmetry (\ref{phi_m}) of this vector system;
\par
2) The vector system periodically passes the base points, where it is convenient to fix its special properties;
\par
3) The special properties of the vector system are determined by the position of the first middle vector $Y_1$ of the Riemann spiral at the base point relative to the real axis of the complex plane and the axes of the moving coordinate system formed by the normal $L$ to the axis of symmetry and the axis of symmetry $M$ of the vector system of the second approximate equation of the Riemann zeta function;
\par
4) The arguments of $\zeta(s)$ and $\overline{\zeta(1-s)}$ have axial symmetry around the normal $L$ to the axis of symmetry of the vector system of the second approximate equation of the Riemann zeta function;
\par
5) The vectors of value of $\zeta(s)$ and $\overline{\zeta(1-s)}$ rotate in different directions at the same speed in the moving coordinate system formed by the normal $L$ to the axis of symmetry and the axis of symmetry $M$ of the vector system of the second approximate equation of the Riemann zeta function;
\par
6) Projections of the vector system of the second approximate equation of the Riemann zeta function on the normal $L$ to the axis of symmetry and on the axis of symmetry $M$ of this vector system are determined by the vectors of value of $\zeta(s)$ and $\overline{\zeta(1-s)}$.
\par
\textit{Thus, the projections of the vector system of the second approximate equation of the Riemann zeta function on the normal $L$ to the axis of symmetry and on the axis of symmetry $M$ of this vector system have periodic properties (which we described earlier) with relation to the base points, so these properties of the projections of the vector system must be observed at any base point, and therefore for all values of the complex variable, where this vector system determines a value of the Riemann zeta function, i.e., on the entire complex plane except the real axis, where, as is known, the Riemann zeta function has only trivial zeros.}
\par
To answer the second question, why there can be no other reason for the modulus of the Riemann zeta function to zero when $\sigma\ne 1/2$, it is necessary to consider possible variants of such a zero transformation, for example:
\par
1) The interval $A_k$, where $\zeta(s)_L=0$, intersects with the interval $C_k$, where $\zeta(s)_M=0$, for any value of a complex variable when $\sigma\ne 1/2$;
\par
2) The condition $\zeta(s)_L=0$ and $\zeta(s)_M=0$ is satisfied in the interval $B_k$ or $D_k$ , where $\zeta(s)_L\ne 0$ and $\zeta(s)_M\ne 0$, for any value of the complex variable in the critical strip.
\par
\textit{It is obvious that the module $\zeta(s)$ can not arbitrarily be reduced to zero because this would require that the zero was reduced module of all vectors $X_n$ and $Y_n$ of the Riemann spiral, which contradicts (\ref{x_vect}, \ref{y_vect}).}
\par
Somebody can think of other reasons why the modulus of the Riemann zeta function can transform to zero when $\sigma\ne 1/2$, we hope that upon careful examination they will all be refuted, because the identified properties of the projections of the vector system of the second approximate equation of the Riemann zeta function on the normal $L$ to the axis of symmetry and on the axis of symmetry $M$ of this system show the direction of confirmation why the Riemann zeta function may not have non-trivial zeros when $\sigma\ne 1/2$.
\par
Now we consider another method of confirmation the Riemann hypothesis, based on the properties of the vector system of the second approximate equation of the Riemann zeta function.
\par
We need to find out in which cases the invariants $L_1$ and $L_2$ (\ref{l_1_l_2}) of this vector system and the vector of the remainder term $R$ of the second approximate equation of the Riemann zeta function can form a triangle (\ref{zeta_zero_vect}).
\par
\begin{figure}[ht]
\centering
\includegraphics[scale=0.6]{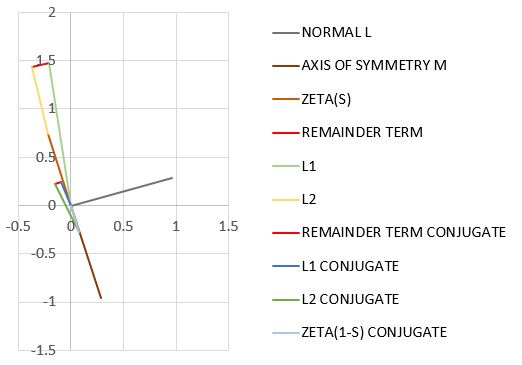}
\caption{Base point \#4525, the vector of values of the Riemann zeta function parallel to the axis of symmetry, $s=0.35+5006.186i$ and $s=0.65+5006.186i$}
\label{fig:gram_point_4525_symmetry_inv}
\end{figure}
\begin{figure}[ht]
\centering
\includegraphics[scale=0.6]{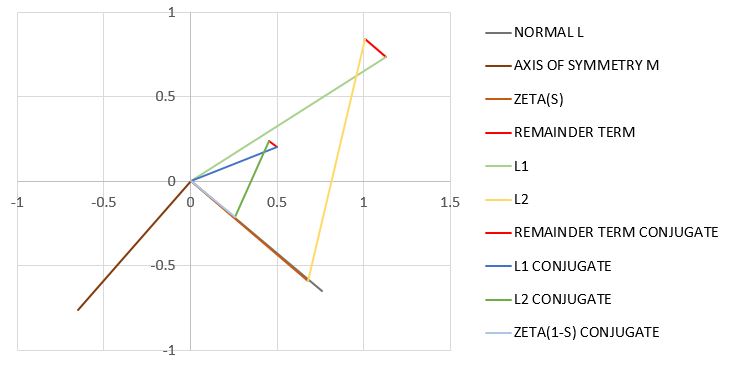}
\caption{Base point \#4525, the vector of values of the Riemann zeta function parallel to the normal to the axis of symmetry, $s=0.35+5006.484i$ and $s=0.65+5006.484i$}
\label{fig:gram_point_4525_normal_inv}
\end{figure}

Consider the invariants $L_1$ and $L_2$ at the special points of the first middle vector $Y_1$ of the Riemann spiral, when the axis of symmetry $M$ of the vector system of the second approximate equation of the Riemann zeta function passes through both the zero of the complex plane and the end of the first middle vector $Y_1$ of the Riemann spiral (fig. \ref{fig:gram_point_4525_symmetry_inv}), and when the normal $L$ to the axis of symmetry of the vector system of the second approximate equation of the Riemann zeta function passes through both the zero of the complex plane and the end of the first middle vector $Y_1$ of the Riemann spiral (fig. \ref{fig:gram_point_4525_normal_inv}).
\par
\begin{figure}[ht]
\centering
\includegraphics[scale=0.6]{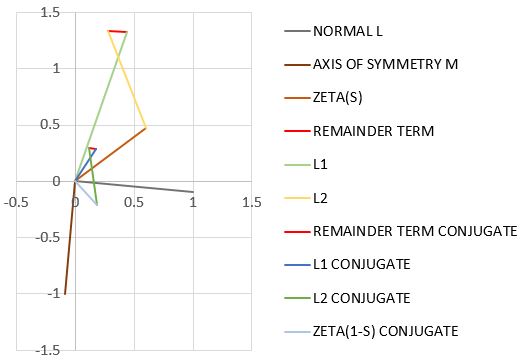}
\caption{Base point \#4525, intermediate position of the vector of values of the Riemann zeta function, $s=0.35+5006.186i$ and $s=0.65+5006.186i$}
\label{fig:gram_point_4525_intermediate_inv}
\end{figure}

In all other cases, in accordance with the continuity of values of the Riemann zeta function, the invariants $L_1$ and $L_2$ will occupy different intermediate positions (fig. \ref{fig:gram_point_4525_intermediate_inv}).
\par
\begin{figure}[ht]
\centering
\includegraphics[scale=0.6]{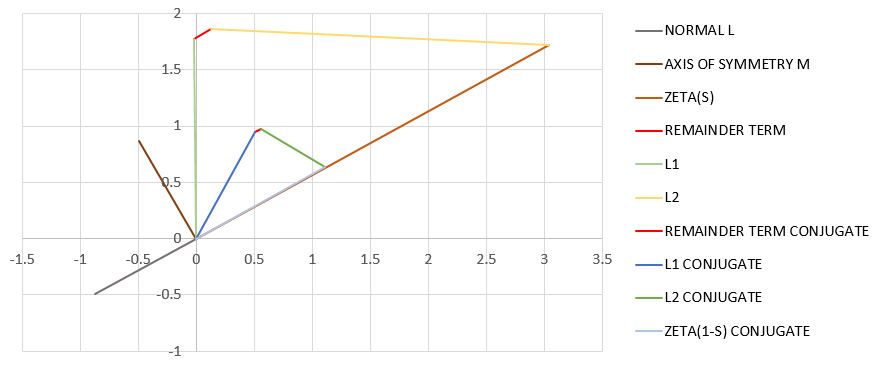}
\caption{Base point \#4520, the vector of values of the Riemann zeta function parallel to the normal to the axis of symmetry, $s=0.35+5001.415i$ and $s=0.65+5001.415i$}
\label{fig:gram_point_4520_normal_inv_opposite}
\end{figure}

When the normal $L$ to the axis of symmetry of the vector system of the second approximate equation of the Riemann zeta function passes through both the zero of the complex plane and the end of the first middle vector $Y_1$ of the Riemann spiral (fig. \ref{fig:gram_point_4520_normal_inv_opposite}) invariants $L_1$ and $L_2$ can occupy a position close to the trapezoid, but in any case, obviously, can not form a triangle.
\par
\textit{It is obvious that the invariants $L_1$ and $L_2$ can occupy the position closest to the triangle only when the axis of symmetry $M$ of the vector system of the second approximate equation of the Riemann zeta function passes through both the zero of the complex plane and the end of the first middle vector $Y_1$ of the Riemann spiral (fig. \ref{fig:gram_point_4525_symmetry_inv}), so we will perform further analysis of the invariants $L_1$ and $L_2$ in this special position of the first middle vector $Y_1$ of the Riemann spiral.}
\par
\begin{figure}[ht]
\centering
\includegraphics[scale=0.6]{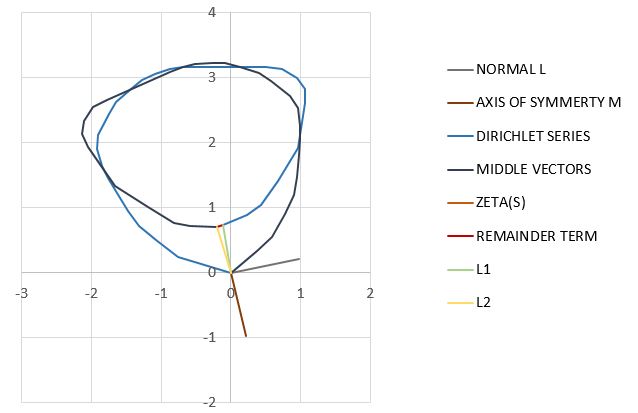}
\caption{Base point \#4525, non-trivial zero of the Riemann zeta function, $s=0.5+5006.208381106i$}
\label{fig:gram_point_4525_symmetry_inv_zero}
\end{figure}

It should be noted that in this position of the axis of symmetry $M$ of the vector system of the second approximate equation of the Riemann zeta function when $\sigma=1/2$ in accordance with the mirror symmetry, this vector system forms a closed polyline, and therefore the Riemann zeta function takes a value of non-trivial zero (fig. \ref{fig:gram_point_4525_symmetry_inv_zero}).
\par
In the ordinate of the non-trivial zero of the Riemann zeta function, the projection modulus of the invariant $L_1$ on the axis of symmetry $M$ of the vector system of the second approximate equation of the Riemann zeta function increases when $\sigma$ decreases and, conversely, it decreases when $\sigma$ increases.
\par
The projection of the invariant $L_2$ on the axis of symmetry $M$ of the vector system of the second approximate equation of the Riemann zeta function changes according to the sign of the projection of its gradient on the axis of symmetry $M$:
\begin{equation}\label{grad_l_2}grad_M L_2=\sum_{n=1}^{m}\Big(\frac{\partial Y_n}{\partial\sigma}\Big)_M;\end{equation}
\par
\begin{figure}[ht]
\centering
\includegraphics[scale=0.6]{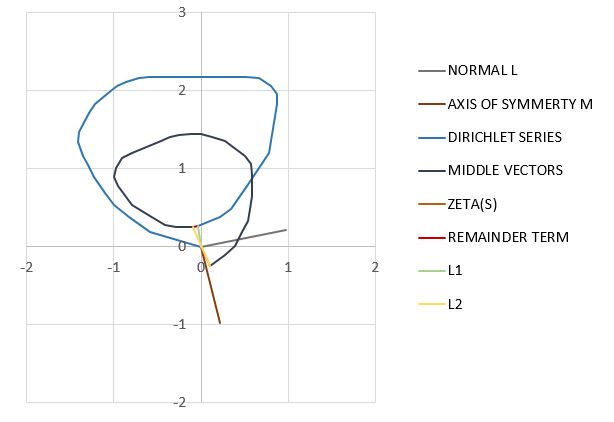}
\caption{Base point \#4525, the ordinate of the non-trivial zero of the Riemann zeta function, positive gradient, $s=0.65+5006.208381106i$}
\label{fig:gram_point_4525_symmetry_3_inv}
\end{figure}
\begin{figure}[ht!]
\centering
\includegraphics[scale=0.6]{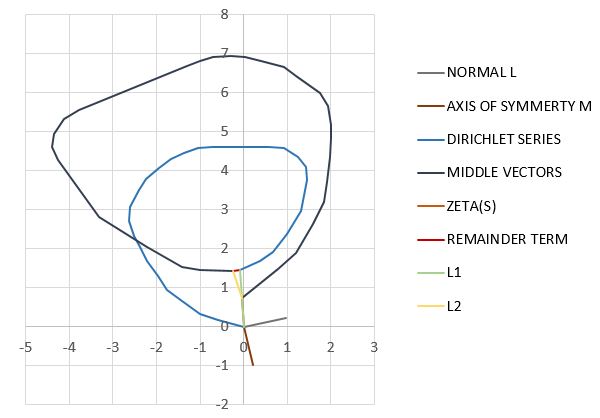}
\caption{Base point \#4525, the ordinate of the non-trivial zero of the Riemann zeta function, positive gradient, $s=0.35+5006.208381106i$}
\label{fig:gram_point_4525_symmetry_2_inv}
\end{figure}
\begin{figure}[ht]
\centering
\includegraphics[scale=0.6]{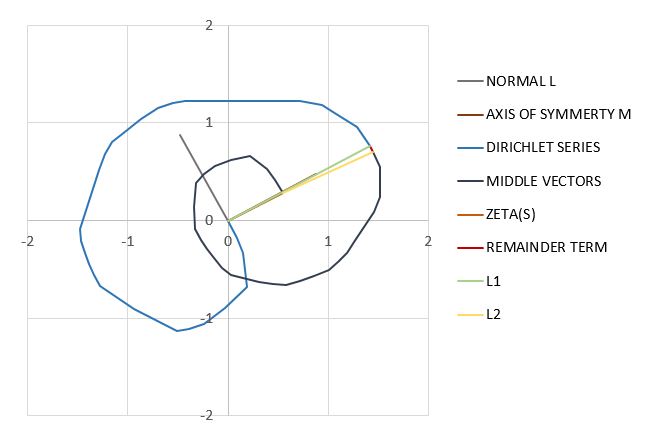}
\caption{Base point \#4521, the ordinate of the non-trivial zero of the Riemann zeta function, negative gradient, $s=0.65+5001.889773627i$}
\label{fig:gram_point_4525_symmetry_5_inv}
\end{figure}
\begin{figure}[ht!]
\centering
\includegraphics[scale=0.6]{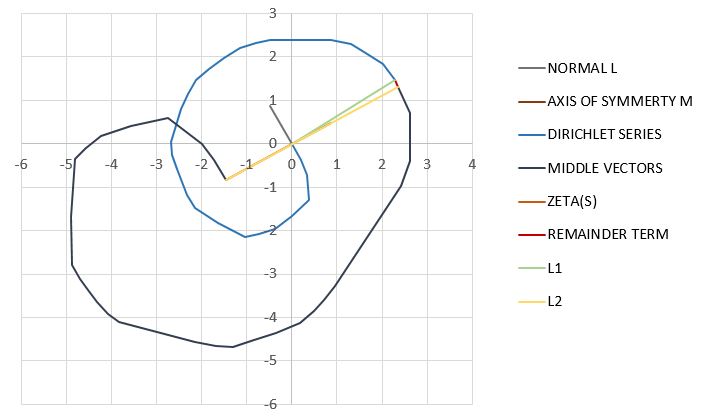}
\caption{Base point \#4521, the ordinate of the non-trivial zero of the Riemann zeta function, negative gradient, $s=0.35+5001.889773627i$}
\label{fig:gram_point_4525_symmetry_4_inv}
\end{figure}

If the sign of the projection of the invariant $L_2$ on the axis of symmetry $M$ of the vector system of the second approximate equation of the Riemann zeta function is equal to the sign of the projection of its gradient on the axis of symmetry $M$, then the projection of the invariant $L_2$ on the axis of symmetry $M$ increases when $\sigma$ increases (fig. \ref{fig:gram_point_4525_symmetry_3_inv}) respectively while $\sigma$ decreases it decreases also (fig. \ref{fig:gram_point_4525_symmetry_2_inv}) and, conversely, if the sign of the projection of the invariant $L_2$ on the axis of symmetry $M$ is not equal to the sign of the projection of its gradient on the axis of symmetry $M$, then when $\sigma$ increases the projection of the invariant $L_2$ on the axis of symmetry $M$ decreases (fig. \ref{fig:gram_point_4525_symmetry_5_inv}) respectively while $\sigma$ decreases it increases (fig. \ref{fig:gram_point_4525_symmetry_4_inv}).
\par
The sign of the projection of the gradient of the invariant $L_2$ on the axis of symmetry $M$ depends on the distribution of angles and modulus of the middle vectors, which becomes obvious if we arrange the middle vectors in increasing order of their angles (fig. \ref{fig:gram_point_4525_symmetry_inv_zero}).
\par
Researches show that when $\sigma\ne 1/2$ at the point when $\zeta(s)_M=0$ the sign of the projection of the gradient of the invariant $L_2$ on the axis of symmetry $M$ is kept.
\par
Now, when we consider the invariants $L_1$ and $L_2$ at the point where the axis of symmetry $M$ of the vector system of the second approximate equation of the Riemann zeta function passes through both the zero of the complex plane and the end of the first middle vector $Y_1$ of the Riemann spiral (fig. \ref{fig:gram_point_4525_symmetry_inv}), it is sufficient to consider the sum of projections of vectors $L_1$, $L_2$ and R on the axis of symmetry $M$, because at this point the sum of projections of vectors $L_1$, $L_2$ and R on the normal $L$ to the axis of symmetry is equal to zero.
\par
Consider separately the sum of projections of invariants $L_1$ and $L_2$ on the axis of symmetry $M$ of the vector system of the second approximate equation of the Riemann zeta function:
\begin{equation}\label{delta_l}\Delta L=\sum_{n=1}^{m}{(X_n(s))_M}+\sum_{n=1}^{m}{(Y_n(s))_M};\end{equation}
\par
and the projection of vector R of the remainder term of the second approximate equation of the Riemann zeta function:
\begin{equation}\label{delta_r}\Delta R=R\sin(\Delta\varphi_R);\end{equation}
\par
where $\Delta\varphi_R$ is the deviation of the vector of the remainder term R of the second approximate equation of the Riemann zeta function from the normal $L$ to the axis of symmetry of the vector system of the second approximate equation of the Riemann zeta function.
\par
Then the vector condition (\ref{zeta_zero_vect}) of the non-trivial zero of the Riemann zeta function can be rewritten as follows:
\begin{equation}\label{zeta_zero_module}|\Delta L|=|\Delta R|;\end{equation}
\par
We already know that in accordance with the mirror symmetry of the vector system of the second approximate equation of the Riemann zeta function when $\sigma=1/2$ at the point when the axis of symmetry $M$ of this vector system passes through both the zero of the complex plane and the end of the first middle vector $Y_1$ of the Riemann spiral (fig. \ref{fig:gram_point_4525_symmetry_inv_zero}), the Riemann zeta function takes a value of non-trivial zero.
\par
Thus $\Delta L=0$ and $\Delta R=0$.
\par
When $\sigma\ne 1/2$, the mirror symmetry of the vector system of the second approximate equation of the Riemann zeta function is broken, in other words, $\Delta L\ne 0$ and $\Delta R\ne 0$, hence the change in $\Delta L$, if the Riemann zeta function can take a value of non-trivial zero, must be compensated by the change in $\Delta R$.
\par
\begin{figure}[ht]
\centering
\includegraphics[scale=0.6]{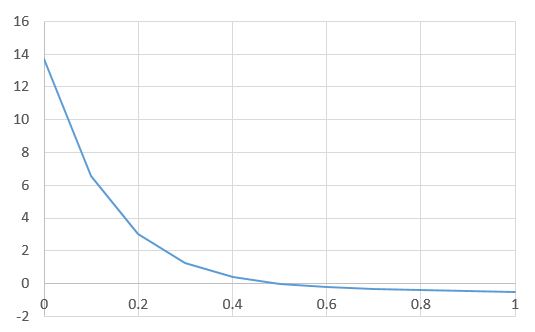}
\caption{Base point \#4525, sum of projections of invariants $L_1$ and $L_2$, ordinate $\zeta(0.35+5006.186i)_L=0$}
\label{fig:projections_sum_l1_l2}
\end{figure}
\begin{figure}[ht]
\centering
\includegraphics[scale=0.6]{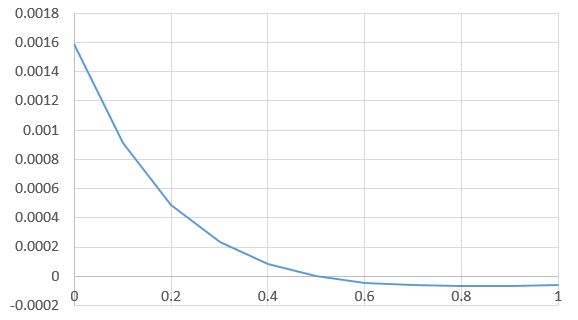}
\caption{Base point \#4525, projection of the remainder term $R$, ordinate $\zeta(0.35+5006.186i)_L=0$}
\label{fig:projection_r}
\end{figure}

Consider the dependence of the sum of projections of the invariants $L_1$ and $L_2$ (fig. \ref{fig:projections_sum_l1_l2}) and the projection of the vector $R$ of the remainder term (fig. \ref{fig:projection_r}) on the axis of symmetry $M$ of the vector system of the second approximate equation of the Riemann zeta function from the real part of a complex number at the point when the axis of symmetry $M$ of this vector system passes through both the zero of the complex plane and the end of the first middle vector $Y_1$ of the Riemann spiral.
\par
It is obvious that these functions are equal only at the point $\sigma=1/2$, i.e. on the critical line.
\par
\begin{figure}[ht]
\centering
\includegraphics[scale=0.6]{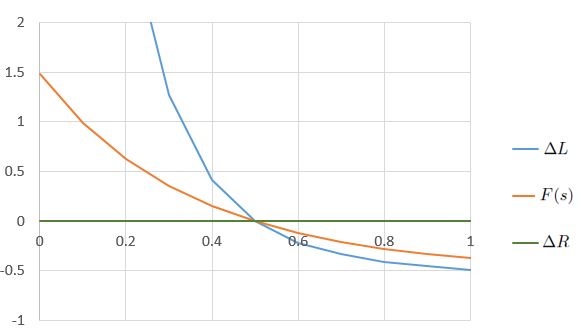}
\caption{Base point \#4525, boundary function $F(s)$, ordinate $\zeta(0.35+5006.186i)_L=0$, $A=0.5$}
\label{fig:boundary_function}
\end{figure}

Consider the \textit{boundary function} that separates values $\Delta L$ and $\Delta R$ (fig. \ref{fig:boundary_function}), at the point where the axis of symmetry $M$ of the vector system of the second approximate equation of the Riemann zeta function passes through both the zero of the complex plane and the end of the first middle vector $Y_1$ of the Riemann spiral:
\begin{equation}\label{boundary_function}F(s)=A \Big(\frac{\sum_{n=1}^{m}{|Y_n(s)|}}{\sum_{n=1}^{m}{|X_n(s)|}}-1\Big);\end{equation}
\par
where $A$ is some constant greater than zero.
\par
It is obvious that the sum of projections of invariants $L_1$ and $L_2$ (fig. \ref{fig:projections_sum_l1_l2}) on the axis of symmetry $M$ of the vector system of the second approximate equation of the Riemann zeta function when $\sigma=0$ is determined by a value of the projection $\zeta(0+it)_M$ (fig. \ref{fig:graphics_projections}), because in this case, $|\Delta L|\gg|\Delta R|$ (fig. \ref{fig:projection_r}), therefore, the sum of projections of invariants $L_1$ and $L_2$ on the axis of symmetry $M$ at all points except $\sigma=1/2$ modulo more than values of the boundary function, while values the projection of the vector  $R$ of the remainder term on the axis of symmetry $M$ of the vector system of the second approximate equation of the Riemann zeta function in all points except $\sigma=1/2$ modulo smaller than values of the boundary functions.
\par
Now we need to find out the dependence of the angle of deviation $\Delta\varphi_R$ of vector remainder member of the second approximate equation of the Riemann zeta function from normal $L$ to the axis of symmetry of the vector system of the second approximate equation of the Riemann zeta function, because a value of this angle determines a value of $\Delta R$.
\par
Although that Riemann recorded the remainder term (\ref{zeta_eq_2_zi}) of the second approximate equation of the Riemann zeta function explicitly, Riemann and other authors use the argument of the remainder term only for the particular case of (\ref{arg_rem_2}) when $\sigma=1/2$, so we use the explicit expression for the CHI functions (\ref{chi_eq_ex3}) and will calculate the argument of the remainder term using the exact values of the Riemann zeta function \cite{ZF}:
\begin{equation}\label{delta_varphi_r}\Delta\varphi_R=\frac{1}{2}Arg(\chi(s))-Arg(R(s))=\frac{1}{2}\arccos(\frac{Re(\chi(s))}{|\chi(s)|})-\arccos(\frac{Re(R(s))}{|R(s)|});\end{equation}
\par
where
\begin{equation}\label{r_vect}R(s)=\zeta(s)-\sum_{n=1}^{m}{X_n(s)}-\sum_{n=1}^{m}{Y_n(s)};\end{equation}
\begin{equation}\label{r_module}|R(s)|=\sqrt{Re(R(s))^2+Im(R(s))^2};\end{equation}
\par
\begin{figure}[ht]
\centering
\includegraphics[scale=0.6]{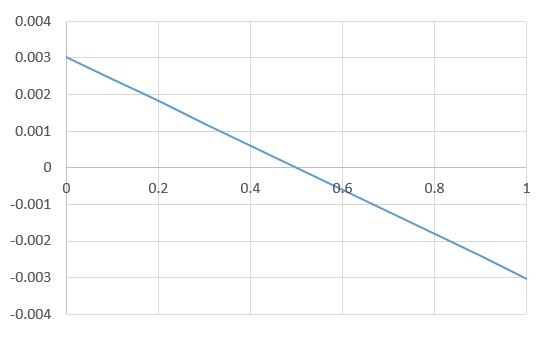}
\caption{Base point \#4525, the deviation angle of vector $\Delta\varphi_R$ of remainder term, rad, ordinate $\zeta(0.35+5006.186i)_L=0$}
\label{fig:delta_varphi_r_real}
\end{figure}

The calculations show a linear relationship (fig. \ref{fig:delta_varphi_r_real}) of the angle of deviation $\Delta\varphi_R$ of the  vector of the remainder term of the second approximate equation of the Riemann zeta function from normal $L$ to the axis of symmetry of this vector system from the real part of a complex number.
\begin{figure}[ht!]
\centering
\includegraphics[scale=0.6]{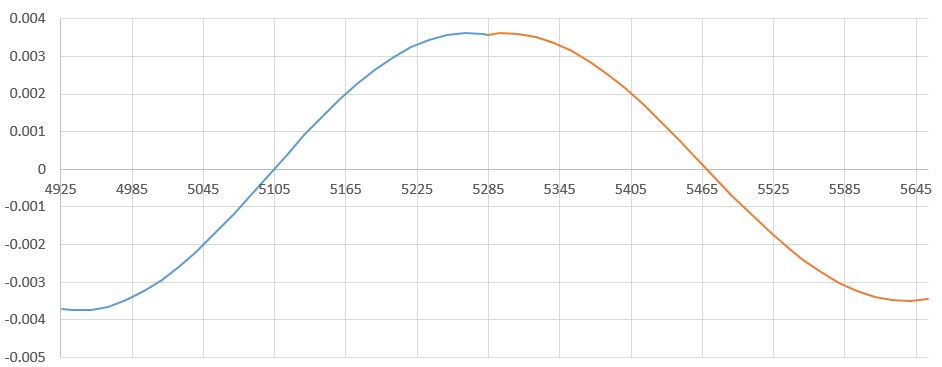}
\caption{The deviation angle of vector $\Delta\varphi_R$ of remainder term, rad (from the imaginary part of a complex number)}
\label{fig:delta_varphi_r_complex}
\end{figure}
\begin{figure}[ht!]
\centering
\includegraphics[scale=0.6]{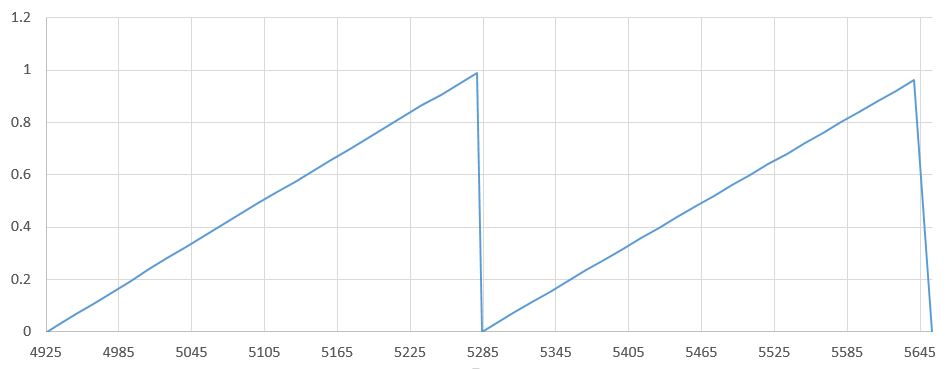}
\caption{Fractional part of $\sqrt{t/2\pi}$ (from imaginary part of complex number)}
\label{fig:frac_part_m_complex}
\end{figure}
\begin{figure}[ht!]
\centering
\includegraphics[scale=0.6]{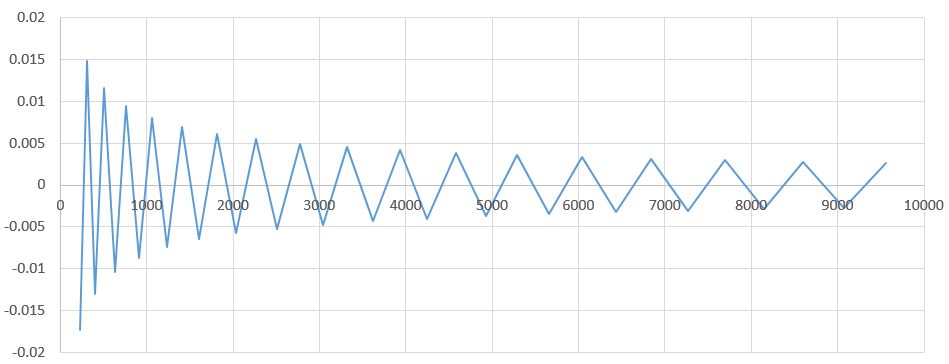}
\caption{The maximum derivation angle of the vector $\Delta\varphi_R$ of the remainder term, rad (from the imaginary part of a complex number)}
\label{fig:delta_varphi_r_max_complex}
\end{figure}
\par
\textit{This dependence is confirmed for any values of the imaginary part of a complex number, since it is determined by the form of the polyline formed by the vector system of the second approximate equation of the Riemann zeta function, for symmetric values $\sigma+it$ (fig. \ref{fig:gram_point_4525_left}) and $1-\sigma+it$ (fig. \ref{fig:gram_point_4525_right}).}
\par
It can be shown that these polylines are mirror congruent, hence all the angles of these polylines are mirror congruent, including the arguments of the remainder terms.
\par
The angle of deviation $\Delta\varphi_R$ of the vector of the remainder term of the second approximate equation of the Riemann zeta function from normal $L$ to the axis of symmetry of the vector system of the second approximate equation dzeta-functions of Riemann has a periodic dependence (fig. \ref{fig:delta_varphi_r_complex}) from the imaginary part of a complex number, with a period equal to an interval (\ref{m_interval}), where fractional part of the expression $\sqrt{t/2\pi}$, which varies from 0 to 1 (fig. \ref{fig:frac_part_m_complex}), has the periodic dependence also.
\par
The angle of deviation $\Delta\varphi_R$ of the vector of the remainder term of the second approximate equation of the Riemann zeta function from normal $L$ to the axis of symmetry of this vector system has the maximum value at the boundaries of the intervals (\ref{m_interval}), and the absolute value of the maximum variance asymptotically decreases (fig. \ref{fig:delta_varphi_r_max_complex}) with the growth of the imaginary part of a complex number, which corresponds to the evaluation of the remainder term of the gamma function:
\begin{equation}\label{mu_limit}\mu(s)\to 0, t\to\infty;\end{equation}
\par
\textit{Thus, at the point where the axis of symmetry $M$ of the vector system of the second approximate equation of the Riemann zeta function passes through both the zero of the complex plane and the end of the first middle vector $Y_1$ of the Riemann spiral, the condition (\ref{zeta_zero_module}) can be true only when $\sigma=1/2$, since $\sigma\ne 1/2$ $|\Delta L|>|\Delta R|$ on any interval (\ref{m_interval}).}
\par
Somebody can think of different options, why expression (\ref{boundary_function}) cannot be \textit{bounding function}, but we hope that they will all be refuted.
\section{Summary}
Was Riemann going to speak at the Berlin Academy of Sciences on the occasion of his election as a corresponding member to present a proof of the asymptotic law of the distribution of prime numbers?
\par
Only Riemann himself could have answered this question, but we are inclined to assume that he did not intend to.
\par
Neither before nor after (Riemann died seven years later) he did not return to the subject publicly.
\par
Riemann did not publish any paper in progress, so the paper on the analysis of the formula, which is now called the Riemann-Siegel formula was published by Siegel based on Riemann's notes, but this paper is rather a development of the analytical analysis of the Riemann zeta function, rather than a proof of the asymptotic distribution law of prime numbers.
\par
The speech was rather about the analytical function of the complex variable.
\par
Riemann showed how bypassing the difficult arguments about the convergence of the series, which defines the function, we can get its analytical continuation using the residue theorem.
\par
Riemann also used a feature of a function of a complex variable for which zeros are its singular points, which carry basic information about the function.
\par
It was in connection with this feature of the complex variable function that the Riemann hypothesis appeared - the hypothesis of the distribution of zeros of the Riemann zeta function.
\par
Strange is also the absence of any mention of the representation of complex numbers by points on the plane, although Riemann in his report uses the rotation of the zeta function by the angle $Arg(\chi(s))/2$, which is certainly an operation on complex numbers as points on the plane.
\par
Such relation to complex numbers seems even more strange in light of the fact that Riemann was a student of Gauss, who was one of the first to introduce the representation of complex numbers by points on the plane.
\par
\textit{In other words, we tend to assume that the zeta function was chosen to show how we can solve the problem by methods of functions of a complex variable, Riemann is not important the problem itself, it is important approach and methods that gives the theory of functions of a complex variable as the apotheosis of the theory of analytic functions.}
\par
Although that a theorem of the distribution of prime numbers is proved analytically, i.e. using the analytical Riemann zeta function, later it was found a proof that uses the functions of a real variable, i.e. an elementary proof.
\par
Thus, the role of the Riemann zeta function has shifted towards regularization, the so-called methods of generalized summation of divergent series.
\par
The Riemann hypothesis seems to belong to such problems, since the generalized Riemann hypothesis deals with an entire class of Dirichlet L-functions.
\par
Zeta function regularization is also used in physics, particularly in quantum field theory.
\par
We prefer to conclude that the main role of the Riemann zeta function is in the understanding of generalized methods for summing asymptotic divergent series and as a consequence of constructing an analytic continuation of the functions of a complex variable.
\par
Based on that on numerous forums \glqq seriously\grqq \ is discussed infinite sum of natural numbers:
\begin{equation}\label{123}1+2+3+ = -\frac{1}{12};\end{equation}
\par
only a few understand the essence of the generalized summation or regularization of divergent series.
\par
Here mathematics encounters philosophy, namely with \textit{the law of unity of opposites.}
\par
The essence of generalized summation is \textit{regularization} (that is why the second name of the method is regularization) - it means that a divergent series cannot exist without its \textit{opposite} - convergent series.
\par
In other words, there is only \textit{one series}, but it converges in one region and diverges in another.
\par
The most natural such series is the Dirichlet series:
\begin{equation}\label{dirichlet}\sum_{n=1}^{\infty}\frac{1}{n^s};\end{equation}
\par
which in the real form was researched by Euler (he first raised the question of the need for the concept of generalized summation), and in the complex form it was considered by Riemann.
\par
Riemaann left the question of generalized summation, skillfully replacing the divergent Dirichlet series by already regularized integral of gamma functions of a complex variable:
\begin{equation}\label{gamma_int}\Gamma (z)={\frac {1}{e^{i2\pi z}-1}}\int \limits _{L}\!t^{\,z-1}e^{-t}\,dt;\end{equation}
\par
The essence of unity \glqq summation\grqq \ of an infinite series is in the definition of \glqq sum\grqq \ this series in the region where this series converges and in the region where \textit{the same} series diverges.
\par
The misconception begins in the definition of \textit{infinite sum.}
In the case of a convergent series, it only seems to us that we can find the \glqq sum\grqq \ of this series.
In fact, we find \textit{limit of partial sums} because in accordance with the convergence of the series, such a limit exists (by definition).
\par
Euler first formulated the need for a different concept of the word \glqq sum\grqq \ in the application to the divergent series (or rather to the series in the region where it diverges), he explained this by the practical need to attach some value divergent series.
\par
We now know that this value is found as \textit{limit of generalized partial sums} of an infinite series.
\par
And the main condition that is imposed on the method of obtaining generalized partial sums (except that there must be a limit) is \textit{regularity}, i.e. the limit of generalized partial sums in the region where the series converges must be equal to the limit of partial sums of this series.
\par
This understanding, as Hardy observed, came only with the development of the theory of the function of a complex variable, namely the notion of \textit{analytic continuation}, which is closely linked to the infinite series that defines the analytic function, and is also closely linked to the fact that this infinite series converges in one region and diverges in another.
\par
An analytic continuation of a function of a complex variable, if it is possible, is unique and this fact (which has a regorous proof) is possible only if there exists a limit of generalized partial sums of an infinite series by which the analytic function is defined, in the region where this series diverges.
\par
In the theory of generalized summation of divergent series, it is also rigorously proved that if the limit of generalized partial sums exists for two different regularization methods (obtaining generalized partial sums), then it has the same value.
\par
\textit{It is this correspondence of different methods of obtaining generalized partial sums and analytic continuation of the function of a complex variable that Hardy had in mind.}
\par
The Dirichlet series $\sum\limits_{n=1}^{\infty}\frac{1}{n^s}$ in complex form defines the Riemann zeta function.
\par
As is known, this series diverges in the critical strip, where the Riemann zeta function has non-trivial zeros.
\par
Hence all non-trivial zeros of the Riemann zeta function are \textit{limit} of generalized partial sums of the Dirichlet series, while the trivial zeros of the Riemann zeta function define odd Bernoulli numbers, which are all zero.
\par
In the Riemann zeta function theory to regularize the Dirichlet series the Euler-MacLaren formula
\footnote{unfortunately, the fact that the Euler-McLaren summation formula is used for the Riemann zeta function as a method of generalized summation of divergent series is not mentioned in all textbooks.}
is traditionally used, as we mentioned earlier, this formula is used if the partial sums of a divergent series are suitable for calculating the generalized sum of that divergent series.
\par
The Euler-MacLaren generalized summation formula allows us to move from an infinite sum to an improper integral, i.e. to the limit of generalized partial sums of the Dirichlet series.
\par
The geometric analysis of partial sums of the Dirichlet series, which defines the Riemann zeta function, allowed us to make an important conclusion that a generalized summation of infinite series is possible if this series diverges asymptotically.
\par
And then we get into the essence of the definition of analytic functions through \textit{asymptotic} infinite series, and a value of the function equal to a value of the asymptote in any case or when the series asymptotically converge or when a series diverges asymptotically, while in order to find the limit of this asymptote when the series converges, use partial sums of this series and when the series diverges, to find the limit of the asymptote, we must use the regular generalized partial sum.
\par
And then we can obtain an analytic continuation of the function given by the \textit{asymptotic infinite series.}
\par
Therefore, the result obtained at the very beginning of the research, namely the use of an alternative method of generalized summation of Cesaro, which was obtained geometrically, is as important as the description of the various options to confirm the Riemann hypothesis.
\section{Conclusion}
\par
We believe that the method of geometric analysis of Dirichlet series, based on the representation of complex numbers by points on the plane, will complement the set of tools of analytical number theory.
\par
The method described in this paper allows us to get into the essence of the function of a complex variable, to identify regularities that explain any value of this function (including in the region where the series that difine this function diverges) and most importantly it gives \textit{an idea of the exact value of zero} as the sum of vectors that form a closed polyline.
\par
After going through the analysis of the vector system of the second approximate equation of the Riemann zeta function, we can formulate the results of the second method of confirmation the Riemann hypothesis without using this vector system, since we only needed it to find the key points that indicate that the Riemann zeta function \textit{cannot have non-trivial zeros in the critical strip, except for the critical line.}
\par
Obviously we can move from the fixed coordinate system formed by the axes $x=Re(s)$ and $y=Im(s)$, to the moving coordinate system formed by axes $L$ with angle $\varphi_L=Arg(\chi(s))/2$ and by axes $M$ with angle $\varphi_M=(Arg(\chi(s))+\pi)/2$ passing through the zero of the complex plane.
\par
Then from the functional equation of the Riemann zeta function:
\begin{equation}\label{zeta_func_eq3}\zeta(s)=\chi(s)\zeta(1-s); \end{equation}
\par
and equality arguments of functions:
\begin{equation}\label{arg6}Arg(\zeta(1-s))=-Arg(\zeta(\overline{1-s})); \end{equation}
\par
using the arithmetic of the arguments of complex numbers, when we rotate the vector of a value of Riemann zeta function by the angle $Arg(\chi(s))/2$ in the negative direction, we obtain:
\begin{equation}\label{zeta_arg}Arg(\zeta(s))-\frac{Arg(\chi(s))}{2}=-(Arg(\zeta(\overline{1-s}))-\frac{Arg(\chi(s))}{2});\end{equation}
\par
Consequently, in the moving coordinate system formed by the axes $L$ and $M$, the angles of the vectors of value of $\zeta(s)$ and $\zeta(\overline{1-s})$ are symmetric about the axis $L$, thus, in accordance with the symmetry of the angles, the vector of value of $\hat\zeta(s)=\zeta (1/2+it)$ is always directed along the axis $L$.
\par
In other words, in the moving coordinate system formed by the axes $L$ and $M$, the vector of value of $ \hat\zeta(s)=\zeta(1/2+it)$ remains fixed and only changes its modulus and sign, while the vectors of value of $\zeta(s)$ and $\zeta(\overline{1-s})$ when $\sigma\ne 1/2$ rotate in this moving coordinate system in different directions with the same speed.
\par
Therefore:
\par
a)  the projection $\hat\zeta(s) _M=\zeta(1/2+it)_M$ is always zero;
\par
b)  the projections $\zeta(s)_L$ and $\zeta(\overline{1-s})_L$ are periodically equal to zero when the vectors of value of $\zeta(s)$ and $\zeta(\overline{1-s})$ both locate along the axis $M$ they have opposite directions and at this point the projections $\zeta(s)_M$ and $\zeta(\overline{1-s})_M$ are not equal to zero, because $\zeta(s)_L$ is an odd harmonic function, and $\zeta(s)_M$ is an even harmonic functionthat conjugate to $\zeta(s)_L$;
\par
c) the projections $\zeta(s)_M$ and $\zeta(\overline{1-s})_M$ are periodically equal to zero when the vectors of value of $\zeta(s)$ and $\zeta(\overline{1-s})$ both locate along the axis $L$ while they have the same direction and at this point the projections $\zeta(s)_L$ and $\zeta(\overline{1-s})_L$ are not equal to zero, since $\zeta(s)_L$ is an odd harmonic function, and $\zeta(s)_M$ is an even harmonic functionthat conjugate to $\zeta(s)_L$;
\par
\begin{figure}[ht!]
\centering
\includegraphics[scale=0.6]{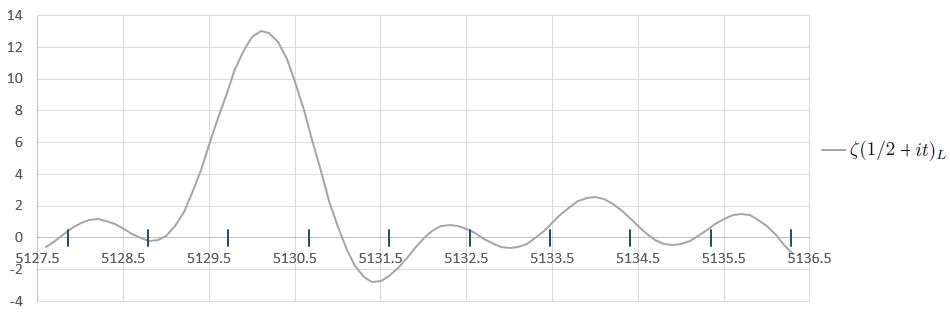}
\caption{Graphics of projections of the Riemann zeta function, $\sigma=1/2$}
\label{fig:graphics_projections_2}
\end{figure}

So only $\hat\zeta(s)_M$ and $\hat\zeta(s)_L$ can be equal to zero at the same time, because $\hat\zeta(s)_L=\zeta(1/2+it)_L$ is an odd harmonic function, and $\hat\zeta(s)_M=\zeta(1/2+it)_M$ is a harmonic function \textit{identically equal to zero} (fig. \ref{fig:graphics_projections_2}).
\par
\begin{figure}[ht!]
\centering
\includegraphics[scale=0.6]{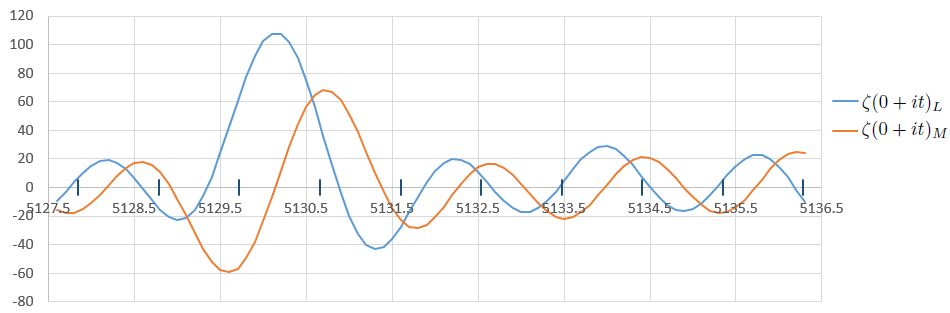}
\caption{Graphics of projections of the Riemann zeta function, $\sigma=0$}
\label{fig:graphics_projections_3}
\end{figure}

\textit{Therefore, when $\sigma\ne 1/2$, the Riemann zeta function cannot be zero, because when $\sigma\ne 1/2$, the projections $\zeta(s)_L$ and $\zeta(s)_M$ cannot be equal to zero at the same time} (fig. \ref{fig:graphics_projections_3}).
\par
There may be a Lemma on conjugate harmonic functions identically nonzero, but we have not found it, just like the Lemma on the symmetric polygon we proved earlier (Lemma 3).
\par
Therefore, to confirm our conclusions, we are forced to return to the vector system of the second approximate equation of the Riemann zeta function.
\par
Vertical marks on the graphs (fig. \ref{fig:graphics_projections_2} and \ref{fig:graphics_projections_3}) is \textit{the base point}, corresponding to solution of the equation:
\begin{equation}\label{}Arg(\chi(\frac{1}{2}+it_k))=(2k-1)\pi;\end{equation}
\par
Using the mirror symmetry property of the vector system of the second approximate equation of the Riemann zeta function when $\sigma=1/2$, we previously defined two types of base points:
\par
$a_1$ - if the first middle vector of the Riemann spiral at the base point is above or along the real axis of the complex plane (the second position corresponds to the non-trivial zero of the Riemann zeta function at the base point is a likely event), in this case the non-trivial zero of the Riemann zeta function that corresponds to this base point is located between this and the previous base point;
\par
$a_2$ - if the first middle vector of the Riemann spiral at the base point is below the real axis of the complex plane, then the non-trivial zero of the Riemann zeta function that corresponds to this base point is located between this and the next base point;
\par
as well as four types of intervals between base points of different types:
\par
$A_1=a_1a_1$ and $A_2=a_2a_2$ - intervals of this kind contain one non-trivial zero of the Riemann zeta function;
\par
$B=a_2a_1$ - interval of this kind contains two non-trivial zeros of the zeta function of Riemann;
\par
$C=a_1a_2$ - interval of this kind does not contain any non-trivial zero of the Riemann zeta function.
\par
We also found that there cannot be the following combinations of intervals:
\par
$A_1A_2$, $A_2A_1$, $BB$ and $CC$;
\par
Hence we have \textit{a fixed set} of combinations of intervals:
\par
$A_1A_1=a_1a_1a_1$, $A_1C=a_1a_1a_2$, $CB=a_1a_2a_1$, $CA_2=a_1a_2a_2$, $BA_1=a_2a_1a_1$, $BC=a_2a_1a_2$, $A_2B=a_2a_2a_1$, $A_2A_2=a_2a_2a_2$;
\par
It is obvious that the sequence of base points which are represented in the graphs (fig. \ref{fig:graphics_projections_2} and \ref{fig:graphics_projections_3}), contains all the possible combinations of intervals:
\par
$a_2a_2a_2a_1a_1a_1a_2a_1a_2a_2=A_2A_2BA_1A_1BCBA_2$;
\par
In accordance with the properties of the vector system of the second approximate equation of the Riemann zeta function, we can determine the sign of the function $\zeta(\sigma+it)_L$ and the sign of its first derivative at each base point.
\par
This requires:
\par
1) determine the type of base point by the position of the first middle vector Riemann spiral relative to the real axis of the complex plane at the base point;
\par
2) determine the direction of the normal $L$ to the axis of symmetry of the vector system of the second approximate equation of the Riemann zeta function at the base point;
\par
Then
\par
A) If the first middle vector Riemann sprial is below the real axis of the complex plane, the function $\zeta (\sigma+it)_L$ will have the sign opposite to the direction of the normal $L$ to the axis of symmetry;
\par
B) If the first middle vector of the Riemann sprial is above the real axis of the complex plane, the function $\zeta (\sigma+it)_L$ will have a sign corresponding to the direction of the normal $L$ to the axis of symmetry;
\par
C) The sign of the first derivative of the function $\zeta(\sigma+it)_L$ always has a sign corresponding to the direction of the normal $L$ to the axis of symmetry;
\par
These rules are executed for any combination of intervals from \textit{a fixed set,} so they are executed for any combination of intervals that can occur.
\par
In the moving coordinate system formed by axes $L$ and $M$ when $\sigma\ne 1/2$, the vector of values of the Riemann zeta function  must rotate on the angle $\pi/2$ from the position corresponding to $\zeta(\sigma+it)_L=0$ to the position corresponding to $\zeta(\sigma+it)_M=0$ as well as from the position corresponding to $\zeta(\sigma+it)_M=0$ to position corresponding to $\zeta(\sigma+it)_L=0$, hence \textit{all the zeros of the function $\zeta(\sigma+it)_M$ when $\sigma\ne 1/2$, lie between the zeros of the function $\zeta(\sigma+it)_L$.}
\par
Corollary 1. The function $\zeta(1/2+it)_L$ has an infinite number of zeros (this statement corresponds to Hardy's theorem \cite{HA2} on an infinite number of non-trivial zeros of the Riemann zeta function on the critical line).
\par
Corollary 2. The number of non-trivial zeros of the Riemann zeta function on the critical line corresponds to the number of base points:
\begin{equation}\label{}N_0(T)=\Bigg[\Big|\frac{T}{2\pi}(\log{\frac{T}{2\pi}}-1)-\frac{1}{8}+\frac{2\mu(T)-\alpha_2}{2\pi}\Big|\Bigg]+2;\end{equation}
\par
where $\mu(T)$ is the remainder term of the gamma function (\ref{mu}) when $\sigma=1/2$;
\par
$\alpha_2$ argument of the CHI function at the second base point.
\par
Corollary 3. The function $\zeta(\sigma+it)_L$ when $\sigma\ne 1/2$ has an infinite number of zeros.
\par
Corollary 4. The function $\zeta(\sigma+it)_M$ when $\sigma\ne 1/2$ has an infinite number of zeros.
\par
Corollary 5. Zeros of function $\zeta(\sigma+it)_L$ and functions $\zeta(\sigma+it)_M$ when $\sigma\ne 1/2$ are not the same, because all zeros of function $\zeta(\sigma+it)_M$ when $\sigma\ne 1/2$ lie between zeros of function $\zeta(\sigma+it)_L$.
\par
Based on the obtained results, we believe that the methods of confirmation of the Riemann hypothesis based on the properties of the vector system of the second approximate equation of the Riemann zeta function will soon lead to its proof.
\section{Gratitudes}
Special thanks to Professor July Dubensky, who allowed to speak at the seminar and instilled confidence in the continuation of the research. Colleagues and friends who listened to the first results and supported throughout the research. My wife, who helped to prepare the presentation at the seminar, supported and believed in success.
\par
The organizers of the conference of Chebyshev collection in Tula, inclusion in the list of participants of this conference allowed to prepare the first version of the paper, but a exclusion from the list of participants mobilized and allowed to obtain new important results of the research.

\bibliographystyle{unsrt}  


\end{document}